\documentclass{amsart}
\usepackage{amssymb}
\usepackage{amsmath}
\usepackage{latexsym}
\usepackage{amsfonts}
\usepackage{amscd}
\usepackage{amsthm}
\newcommand{\gt}[1]{\mathfrak{#1}}

\newcommand{\mc}[1]{\mathcal{#1}}

\newcommand{\RR}{{\mathbb R}}%Reals
\newcommand{\CC}{{\mathbb C}}%Complex
\newcommand{\ZZ}{{\mathbb Z}}%Integers
\newcommand{\ZZh}{\ZZ+\tfrac{1}{2}}%Integers plus half
\newcommand{\hZZ}{\tfrac{1}{2}\ZZ}%half integers

\newcommand{\QQ}{{\mathbb Q}}%Rationals
\newcommand{\hh}{{\mathbf h}}%hyperbolic space
\newcommand{\FF}{{\mathbb F}}%Field
\newcommand{\ii}{{\bf i}}%Sqrt of -1
\newcommand{\lab}{{\langle}}    %Left angle brackets
\newcommand{\rab}{{\rangle}}    %Right angle brackets

\newcommand{\w}{\omega}     %cube root of 1
\newcommand{\wc}{\bar{\omega}}  %conj cube root of 1

\newcommand{\fset}{\Sigma}    %a finite set

\newcommand{\End}{\operatorname{End}}
\newcommand{\Aut}{\operatorname{Aut}}
\newcommand{\Inn}{\operatorname{Inn}}
\newcommand{\Gal}{\operatorname{Gal}}
\newcommand{\Id}{\operatorname{Id}}

\newcommand{\GL}{\operatorname{\textsl{GL}}}      %GL group
\newcommand{\PSL}{\operatorname{\textsl{PSL}}}    %PSL group
\newcommand{\SL}{\operatorname{\textsl{SL}}}      %SL group
\newcommand{\Sp}{\operatorname{\textsl{Spin}}}    %Spin group
\newcommand{\SO}{\operatorname{\textsl{SO}}}    %SO group

\newcommand{\Co}{\operatorname{\textsl{Co}}}    %Conway 1
\newcommand{\LL}{\Lambda}     %Leech lattice

\newcommand{\Cl}{{\rm Cliff}}   %Clifford algebra
\newcommand{\Cm}{{\rm CM}}      %Clifford algebra module
\newcommand{\vfn}{{V^{f\natural}}}      %\cdot 1 SVOA
\newcommand{\afn}{{A^{f\natural}}}    %Fermionic \cdot 1 SVOA

\newcommand{\vn}{V^{\natural}}     %Moonshine Module

\newcommand{\ogi}{\theta}   %``involution'' in \og
\newcommand{\ogz}{\gt{z}}

\newcommand{\vac}{{\bf 1}}      %Vacuum
\newcommand{\cas}{{\bf \omega}} %Conformal element
\newcommand{\scas}{{\bf \tau}}  %Super conformal element

\newtheorem{thm}{Theorem}[section]
\newtheorem{cor}[thm]{Corollary}
\newtheorem{lem}[thm]{Lemma}
\newtheorem{prop}[thm]{Proposition}

\theoremstyle{definition}

\theoremstyle{remark}
\newtheorem{rmk}[thm]{Remark}

\numberwithin{equation}{subsection}
\theoremstyle{plain}
\newtheorem*{sthm}{Theorem}

\begin{document}

\title{Super-moonshine for Conway's largest sporadic group
%     \footnote{{\it MSC2000}
%          Primary 17B69; Secondary 20D08.
%               }
     }

\author{John F. Duncan
%          \footnote{
%          Department of Mathematics,
%          Yale University,
%          P.O. Box 208283,
%          New Haven,
%          CT 06520-8283,
%          U.S.A.}
%          {}\footnote{
%          email: {\tt john.duncan@yale.edu}}
%          %homepage: {\tt
%          %http://math.yale.edu/\~{}jfd22/} }
%               }
%\address{Department of Mathematics, Yale University,
%    New Haven, Connecticut 06520}
%\email{john.duncan@yale.edu}
}


\date{September 14, 2006}

%\maketitle

\begin{abstract}
We study a self-dual $N=1$ super vertex operator algebra and prove
that the full symmetry group is Conway's largest sporadic simple
group. We verify a uniqueness result which is analogous to that
conjectured to characterize the Moonshine vertex operator algebra.
The action of the automorphism group is sufficiently transparent
that one can derive explicit expressions for all the McKay--Thompson
series. A corollary of the construction is that the perfect double
cover of the Conway group may be characterized as a point-stabilizer
in a spin module for the Spin group associated to a $24$ dimensional
Euclidean space.
\end{abstract}

\maketitle

\section{Introduction}

The preeminent example of the structure of vertex operator algebra
(VOA) is the Moonshine VOA denoted $\vn$ which was first constructed
in \cite{FLM} and whose full automorphism group is the Monster
sporadic group.

Following \cite{HohnPhD} we say that a VOA is nice when it is
$C_2$-cofinite and satisfies a certain natural grading condition
(see \S\ref{sec:SVOAstruc:SVOAs}), and we make a similar definition
for super vertex operator algebras (SVOAs). We say that a VOA is
rational when all of its modules are completely reducible (see
\S\ref{sec:SVOAMods}). Then conjecturally $\vn$ may be characterized
among nice rational VOAs by the following properties:
\begin{itemize}
\item   self-dual
\item   rank $24$
\item   no small elements
\end{itemize}
where a self-dual VOA is one that has no non-trivial irreducible
modules other than itself, and we write ``no small elements'' to
mean no non-trivial vectors in the degree one subspace, since in a
nice VOA this is the $L(0)$-homogeneous subspace with smallest
degree that can be trivial. In this note we study what may be viewed
as a super analogue of $\vn$. More specifically, we study an object
$\afn$ characterized among nice rational $N=1$ SVOAs by the
following properties.
\begin{itemize}
\item   self-dual
\item   rank $12$
\item   no small elements
\end{itemize}
An $N=1$ SVOA is an SVOA which admits a representation of the
Neveu--Schwarz superalgebra, and now ``no small elements'' means no
non-trivial vectors with degree ${1}/{2}$. (We define an SVOA to be
rational when its even sub-VOA is rational, and a self-dual SVOA is
an SVOA with no irreducible modules other than itself.)

We find that the full automorphism group of $\afn$ is Conway's
largest sporadic group $\Co_1$. Thus by considering the graded
traces of elements of $\Aut(\afn)$ acting on $\afn$, that is, the
McKay--Thompson series, one can associate modular functions to the
conjugacy classes of $\Co_1$, and we obtain moonshine for $\Co_1$.
This is directly analogous to the moonshine which exists for the
Monster simple group, was first observed in \cite{ConNorMM} and is
to some extent explained by the existence of $\vn$.

The main results of this paper are the following three theorems.
\begin{sthm}[\ref{ThmCnstafn}]
The space $\afn$ admits a structure of self-dual rational $N=1$
SVOA.
\end{sthm}

\begin{sthm}[\ref{ThmSymms}]
The automorphism group of the $N=1$ SVOA structure on $\afn$ is
isomorphic to Conway's largest sporadic group, $\Co_1$.
\end{sthm}

\begin{sthm}[\ref{ThmUniq}]
The $N=1$ SVOA $\afn$ is characterized among nice rational $N=1$
SVOAs by the properties: self-dual, rank $12$, trivial degree $1/2$
subspace.
\end{sthm}

The earliest evidence in the mathematical literature that there
might be an object such as $\afn$ was given in \cite{FLMBerk} where
it was suggested to study a $\ZZ/2$-orbifold of a certain SVOA
$V_L^f=A_L\otimes V_L$ associated to the $E_8$ lattice. Here $L$ is
a lattice of $E_8$ type, $V_L$ denotes the VOA associated to $L$,
and $A_L$ denotes the Clifford module SVOA associated to the space
$\CC\otimes_{\ZZ}L$. Since the lattice $L=E_8$ is self-dual, $V^f_L$
is a self-dual SVOA, and one finds that the graded character of
$V^f_L$ satisfies
\begin{gather}
        {\sf tr}|_{V^f_L}q^{L(0)-c/24}
        =\frac{\theta_{E_8}(\tau)\eta(\tau)^8}
        {\eta(\tau/2)^{8}\eta(2\tau)^{8}}
            =q^{-1/2}(1+8q^{1/2}+276q+2048q^{3/2}+\ldots)
\end{gather}
and is a Hauptmodul for a certain genus zero subgroup of the Modular
group $\bar{\Gamma}=\PSL(2,\ZZ)$. One can check that, but for the
constant term, these coefficients exhibit moonshine phenomena for
$\Co_1$. For example, we have that $276$ is the dimension of an
irreducible module for $\Co_1$, and $2048=1+276+1771$ is a possible
decomposition of the degree ${3}/{2}$ subspace into irreducibles for
$\Co_1$. The only problem being that there is no irreducible
representation of $\Co_1$ of dimension $8$, and the space
corresponding to the constant term in the character of $V^f_L$ would
have to be a sum of trivial modules were it a $\Co_1$-module at all.
As observed in \cite{FLMBerk}, orbifolding $V^f_L$ by a suitable
lift of $-1$ on $L$, one obtains a space $\vfn$ with ${\sf
tr}|_{\vfn} q^{L(0)}={\sf tr}|_{V^f_L} q^{L(0)}-8$; that is, a space
with the correct character for $\Co_1$.
\begin{gather}
        {\sf tr}|_{\vfn}q^{L(0)-c/24}
        =\frac{\theta_{E_8}(\tau)\eta(\tau)^8}
        {\eta(\tau/2)^{8}\eta(2\tau)^{8}}-8
            =q^{-1/2}(1+276q+2048q^{3/2}+\ldots)
\end{gather}
Also one finds that $\vfn$ admits a reasonably transparent action by
a group of the shape $2^{1+8}.(W_{E_8}'/\{\pm 1\})$ which is the same
as that of a certain involution centralizer and maximal subgroup in
$\Co_1$. (Here $W_{E_8}'$ denotes the derived subgroup of the Weyl
group of type $E_8$.)

We note here that the existence of $\vfn$ has certainly been known
to Richard E. Borcherds for some time. Also, an action of $\Co_1$ on
the SVOA underlying what we will call $\afn$ was considered in
\cite{BorRybMMIII}.

The space $\vfn$ may be described as
\begin{gather}
    \vfn=(V^f_L)^0\oplus (V^f_L)^0_{\ogi}
\end{gather}
where $(V^f_L)_{\ogi}$ denotes a $\ogi$-twisted $V^f_L$-module and
$\ogi=\ogi_f\otimes \ogi_b$ is an involution on $V^f_L$ obtained by
letting $\ogi_f$ be the parity involution on $A_L$, and letting
$\ogi_b$ be a lift to $\Aut(V_L)$ of the $-1$ symmetry on $L$. The
$\ogi$-twisted module $(V^f_L)_{\ogi}$ may be described as a tensor
product of twisted modules $(V^f_L)_{\ogi}=(A_L)_{\ogi_f}\otimes
(V_L)_{\ogi_b}$, and $(V^f_L)^0$ and $(V^f_L)_{\ogi}^0$ denote
$\ogi$-fixed points. Now one is in a situation very similar to that
which gives rise to the Moonshine VOA $\vn$ via the Leech lattice
VOA as is carried out in \cite{FLM}, and one can hope that a similar
approach as that used in \cite{FLM} would yield the desired result:
an SVOA structure on $\vfn$ with an action by $\Co_1$.

\medskip

We find it convenient to pursue a different approach. We construct
an SVOA $\afn$ whose graded character coincides with that of $\vfn$,
and we do so using a purely fermionic construction; that is, using
Clifford module SVOAs. It turns out that this fermionic construction
enables one to analyze the symmetries in quite a transparent way. We
find that there is a specific vector in the degree ${3}/{2}$
subspace of $\afn$ such that the corresponding point stabilizer in
the full group of SVOA automorphisms of $\afn$ is precisely the
group $\Co_1$. Furthermore, this vector naturally gives rise to a
representation of the Neveu--Schwarz superalgebra on $\afn$, and
thus $\Co_1$ is realized as the full group of automorphisms of a
particular $N=1$ SVOA structure on $\afn$. These results are
presented in \S\ref{sec:strucafn}.

%\medskip

In \S\ref{SecUniq} we consider an analogue of the uniqueness
conjecture for $\vn$. It turns out that the $N=1$ SVOA structure on
$\afn$ is unique in the sense that the vectors in $(\afn)_{3/2}$
that give rise to a representation of the Neveu--Schwarz
superalgebra on $\afn$ form a single orbit under the action of
$\Sp_{24}(\RR)$ on $\afn$. Making use of modular invariance for
$n$-point functions on a self-dual SVOA, due to \cite{ZhuModInv} and
\cite{HohnPhD}, and utilizing also some ideas of \cite{DonMasEfctCC}
and \cite{DonMasHlmVOA}, one can show that the SVOA structure
underlying $\afn$ is characterized up to some technical conditions
by the by now familiar properties:
\begin{itemize}
\item   self-dual
\item   rank $12$
\item   no small elements
\end{itemize}
Combining these results, we find that $\afn$ is characterized among
$N=1$ SVOAs by the above three properties. This uniqueness result
accentuates the analogy between $\vfn$ and certain other celebrated
algebraic structures: the Golay code, the Leech lattice, and
(conjecturally) the Moonshine VOA $\vn$.

%\medskip

The homogeneous subspace of $\afn$ of degree $3/2$ may be identified
with a half-spin module over $\Sp_{24}(\RR)$. The main idea behind
the construction of $\afn$ is to realize this spin module in such a
way that the vector giving rise to the $N=1$ structure on $\afn$ is
as obvious as possible. Essentially, we achieve this by using the
Golay code to construct a particular idempotent in the Clifford
algebra of a $24$ dimensional vector space. The automorphism group
of $\afn$ turns out to be the quotient by center of the subgroup of
$\Sp_{24}(\RR)$ fixing this idempotent. Thus a curious corollary of
the construction is that the group $\Co_0$ (the perfect double cover
of $\Co_1$) may be characterized as a point stabilizer in a spin
module over $\Sp_{24}(\RR)$.

The Golay code is characterized among length $24$ doubly-even linear
binary codes (see \S\ref{Sec:Notation}) by the conditions of
self-duality and having no ``small elements'' (that is, no weight
$4$ codewords). The Golay code is an important ingredient in the
construction of $\afn$, and these defining properties yield direct
influence upon the structure of $\afn$. For example, the condition
``no small elements'' allows one to conclude that $\Aut(\afn)$ is
finite, and the uniqueness of the Golay code ultimately entails the
uniqueness of the $N=1$ structure on $\afn$. The uniqueness of the
$N=1$ structure in turn allows one to formulate the following
characterization of $\Co_0$.
\begin{sthm}[\ref{Thm:Co0Chrztn}]
Let $M$ be a spin module for $\Sp_{24}(\RR)$ %with invariant
%bilinear form denoted $\lab\cdot\,,\cdot\rab$,
and let $t\in M$ such that $\lab xt,t\rab=0$ whenever
$x\in\Sp_{24}(\RR)$ is an involution with ${\rm tr}|_{24}x=16$. Then
the subgroup of $\Sp_{24}(\RR)$ fixing $t$ is isomorphic to $\Co_0$.
\end{sthm}
It is perhaps interesting to note that the Leech lattice, which has
automorphism group $\Co_0$ and which furnishes a popular definition
of this group \cite{ConCnstCo0}, does not figure directly in our
construction of $\afn$. Although the Golay code does play a
prominent role, the uniqueness of $\afn$, or alternatively the above
Theorem~\ref{Thm:Co0Chrztn}, provide definitions for the group
$\Co_0$ relying neither on the Leech lattice nor the Golay code.

\medskip

The SVOA $\vfn$ constructed from the $E_8$ lattice admits a
natural structure of $N=1$ SVOA in analogy with the way in which a
usual lattice VOA is naturally equipped with a Virasoro element.
Thus a corollary of the uniqueness result for $\afn$ is that
$\afn$ is isomorphic to the $N=1$ SVOA $\vfn$ discussed above. In
the penultimate section \S\ref{LatConst} we consider the
construction of $\vfn$ in more detail, and we indicate how to
construct an explicit isomorphism with $\afn$.

In the final section \S\ref{sec:MTseries} we consider the
McKay--Thompson series associated to elements of $\Co_1$ acting on
$\afn$.
%each of $\vfn$ and $\afn$. In the former case one can treat any
%element of the extraspecial type involution centralizer in $\Co_1$,
%and in the later case one can derive explicit formulae for any
%element of $\Co_1$ in terms of its associated frame shape. In
%particular, one can derive explicit expressions for the
One can derive explicit expressions for each of the McKay--Thompson
series associated to elements of $\Co_1$ in terms of the frame
shapes of the corresponding preimages in $\Co_0$. These expressions
are recorded in Theorem~\ref{ThmChars}.

%%%%%%%%%%%%%%%%%%%%%%%%%%%%%%%%%%%%%%%%%%%%%%%%%%%%%%%%%

\subsection{Notation}\label{Sec:Notation}

If $M$ is a vector space over $\mathbb{F}$ and $\mathbb{E}$ is a
field containing $\FF$, we write $_{\mathbb{E}}M =\mathbb{E}
\otimes_{\mathbb{F}}M$ for the vector space over $\mathbb{E}$
obtained by extension of scalars. For the remainder we shall use
$\FF$ to denote either $\RR$ or $\CC$. We choose a square root of
$-1$ in $\CC$ and denote it by $\ii$. For $q$ a prime power $\FF_q$
shall denote a field with $q$ elements. For $G$ a finite group we
write $\FF G$ for the group algebra of $G$ over $\FF$. %If $\{a_i\}$
%is a set of elements in some algebraic object $A$, then we write
%$\langle a_i\rangle$ for the sub-object of $A$ generated by the
%$a_i$.
\medskip

For $\fset$ a finite set, we denote the power set of $\fset$ by
$\mc{P}(\fset)$. The set operation of symmetric difference (which we
denote by $+$) equips $\mc{P}(\fset)$ with a structure of
$\FF_2$-vector space, and with this structure in mind, we sometimes
write $\FF_2^{\fset}$ in place of $\mc{P}(\fset)$. Suppose that
$\fset$ has $N$ elements. The space $\FF_2^{\fset}$ comes equipped
with a function $w:\FF_2^{\fset}\to \{0,1,\ldots,N\}$ called {\em
weight}, which assigns to an element $\gamma\in\FF_2^{\fset}$ the
cardinality of the corresponding element of $\mc{P}(\fset)$. An
$\FF_2$-subspace of $\FF_2^{\fset}$ is called a {\em linear binary
code of length $N$}. A linear binary code $\mc{C}$ is called {\em
even} if $2|w(C)$ for all $C\in\mc{C}$, and is called {\em
doubly-even} if $4|w(C)$ for all $C\in\mc{C}$. For
$\mc{C}<\FF_2^{\fset}$ a linear binary code, the {\em co-code}
$\mc{C}^*$ is the space $\FF_2^{\fset}/\mc{C}$. We write $X\mapsto
\bar{X}$ for the canonical map $\FF_2^{\fset}\to\mc{C}^*$. The
weight function on $\FF_2^{\fset}$ induces a function $w^*$ on
$\mc{C}^*$ called {\em co-weight}, which assigns to
$\bar{X}\in\mc{C}^*$ the minimum weight amongst all preimages of
$\bar{X}$ in $\FF_2^{\fset}$. Once a choice of code $\mc{C}$ has
been made, it will be convenient to regard $w^*$ as a function on
$\FF_2^{\fset}$ by setting $w^*(X):=w^*(\bar{X})$ for
$X\in\FF_2^{\fset}$.
%The reader may refer to \cite{CoS93} for much
%information about codes.

\medskip

The Virasoro algebra is the universal central extension of the Lie
algebra of polynomial vector fields on the circle (see
\S\ref{sec:PreN=1SVOAs} for an algebraic formulation). In the case
that a vector space $M$ admits an action by the Virasoro algebra,
and the action of $L(0)$ is diagonalizable, we write
$M=\coprod_nM_n$ where $M_n=\{v\in M\mid L(0)v=nv\}$. We call $M_n$
the homogeneous subspace of degree $n$, and we write ${\rm
deg}(u)=n$ for $u\in M_n$.

%\medskip

When $M$ is a super vector space, we write $M=M_{\bar{0}}\oplus
M_{\bar{1}}$ for the superspace decomposition. For $u\in M$ we set
$|u|=k$ when $u$ is $\ZZ/2$-homogeneous and $u\in M_{\bar{k}}$ for
$k\in\{0,1\}$. The dual space $M^*$ has a natural superspace
structure such that $M^*_{\gamma}=(M_{\gamma})^*$ for $\gamma\in
\{\bar{0},\bar{1}\}$.

There are various vector spaces throughout the paper that admit an
action by a linear automorphism of order two denoted $\ogi$.
Suppose that $M$ is such a space. Then we write $M^k$ for the
$\ogi$-eigenspace %of $U$
with eigenvalue $(-1)^k$. Note that $M$ may be a super vector space,
and $M=M^0\oplus M^1$ may or may not coincide with the superspace
grading on $M$.

\medskip

We denote by $D_z$ the operator on formal Laurent series which is
formal differentiation in the variable $z$, so that if $f(z)=\sum
f_rz^{-r-1}\in V\{z\}$ is a formal Laurent series with coefficients
in some space $V$, we have $D_zf(z)=\sum(-r)f_{r-1}z^{-r-1}$. For
$m$ a positive integer, we set $D_z^{(m)}=\tfrac{1}{m!}D_z^m$.

%\medskip

As is customary, we use $\eta(\tau)$ to denote the Dedekind eta
function.
\begin{gather}\label{Dedetafun}
        \eta(\tau)=q^{1/24}\prod_{n= 1}^{\infty}(1-q^n)
\end{gather}
Here $q=e^{2\pi\ii\tau}$ and $\tau$ is a variable in the upper half
plane $\hh=\{\sigma+\ii t\mid t>0\}$.

\medskip

We write $\Co_1$ for an abstract group isomorphic to Conway's
largest sporadic group. We write $\Co_0$ for an abstract group
isomorphic to the perfect double cover of $\Co_1$. It is well known
that $\Co_0$ is isomorphic to the automorphism group of the Leech
lattice \cite{ConCnstCo0}, and in particular, admits an orthogonal
representation of degree $24$ writable over $\ZZ$.

\medskip

The most specialized notations arise in \S\ref{sec:cliffalgs}. We
include here a list of them, with the relevant subsections indicated
in brackets. %They are grouped roughly according to similarity of
%appearance, rather than by order of appearance, so that the list may
%be easier to search through, whenever the need might arise.
\begin{small}
\begin{list}{}{     \itemsep -3pt
                    \topsep 3pt
                         }
\item[$\gt{u}$]  A real or complex vector space of even dimension with
non-degenerate bilinear form, assumed to be positive definite in the
real case (\S\ref{sec:cliffalgs:struc}).

\item[$\{e_i\}_{i\in\fset}$]     A basis for $\gt{u}$,
orthonormal in the sense that $\lab e_i,e_j\rab=\delta_{ij}$ for
$i,j\in\fset$ (\S\ref{sec:cliffalgs:struc}).

\item[$e_I$]   We write $e_I$ for
$e_{i_1}\cdots e_{i_k}\in\Cl(\gt{u})$ when $I=\{i_1,\ldots,i_k\}$ is
a subset of $\fset$ and $i_1<\cdots<i_k$
(\S\ref{sec:cliffalgs:struc}).

\item[$g(\cdot)$]   We write $g\mapsto
g(\cdot)$ for the natural homomorphism $\Sp(\gt{u})\to\SO(\gt{u})$.
Regarding $g\in\Sp(\gt{u})$ as an element of $\Cl(\gt{u})^{\times}$
we have $g(u)=gug^{-1}$ in $\Cl(\gt{u})$ for $u\in\gt{u}$. More
generally, we write $g(x)$ for $gxg^{-1}$ when $x$ is any element of
$\Cl(\gt{u})$.

\item[$e_I(\cdot)$] When $I$ is even, $e_I\in\Cl(\gt{u})$ lies also
in $\Sp(\gt{u})$, and $e_I(a)$ denotes $e_Ia{e_I}^{-1}$ when
$a\in\Cl(\gt{u})$.

\item[$\alpha$] The main anti-automorphism of a Clifford algebra
(\S\ref{sec:cliffalgs:struc}).

\item[$\ogz$]  We denote $e_{\fset}\in\Sp(\gt{u})$ also by $\ogz$
(\S\ref{sec:cliffalgs:spin}).

\item[$\ogi$]  The map which is $-\Id$ on $\gt{u}$, or the parity
involution on $\Cl(\gt{u})$ (\S\ref{sec:cliffalgs:struc}), or the
parity involution on $A(\gt{u})_{\Theta}$
(\S\ref{sec:cliffalgs:SVOAs}).

\item[$1_E$]   A vector in $\Cm(\gt{u})_E$ such that $x1_E=1_E$ for
$x$ in $E$ (\S\ref{sec:cliffalgs:mods}).

\item[$\vac_E$]     The vector corresponding to $1_E$ under the
identification between $\Cm(\gt{u})_{E}$ and
$(A(\gt{u})_{\ogi,E})_{N/8}$ when $\gt{u}$ has dimension $2N$
(\S\ref{sec:cliffalgs:SVOAs}).

\item[$\fset$] A finite ordered set indexing an orthonormal basis
for $\gt{u}$ (\S\ref{sec:cliffalgs:struc}).

\item[$\mc{E}$]     A label for the basis
$\{e_i\}_{i\in\fset}$ (\S\ref{sec:cliffalgs:struc}).

\item[$E$]     A subgroup of $\Cl(\gt{u})^{\times}$ homogeneous with
respect to the $\FF_2^{\mc{E}}$ grading on $\Cl(\gt{u})$
(\S\ref{sec:cliffalgs:mods}).

\item[$A(\gt{u})$]  The Clifford module SVOA associated to the
vector space $\gt{u}$ (\S\ref{sec:cliffalgs:SVOAs}).

\item[$A(\gt{u})_{\ogi}$]     The canonically $\ogi$-twisted module
over $A(\gt{u})$ (\S\ref{sec:cliffalgs:SVOAs}).

\item[$A(\gt{u})_{\ogi,E}$]   The $\ogi$-twisted module
$A(\gt{u})_{\ogi}$ realized in such a way that the subspace of
minimal degree is identified with $\Cm(\gt{u})_E$
(\S\ref{sec:cliffalgs:SVOAs}).

\item[$A(\gt{u})_{\Theta}$]   The direct sum of $A(\gt{u})$-modules
$A(\gt{u})\oplus A(\gt{u})_{\ogi}$ (\S\ref{sec:cliffalgs:SVOAs}).

\item[$\mc{C}(E)$]  The linear binary code on $\fset$ consisting of
elements $I$ in $\FF_2^{\fset}$ for which $E$ has non-trivial
intersection with $\FF e_I\subset\Cl(\gt{u})$
(\S\ref{sec:cliffalgs:mods}).

\item[$\Cm(\gt{u})_E$]   The module over $\Cl(\gt{u})$ induced from
a trivial module over $E$ (\S\ref{sec:cliffalgs:mods}).

\item[$\Cl(\gt{u})$]     The Clifford algebra associated to the vector space
$\gt{u}$ (\S\ref{sec:cliffalgs:struc}).

\item[$\Sp(\gt{u})$]     The spin group associated to the vector space
$\gt{u}$ (\S\ref{sec:cliffalgs:spin}).

\item[$\lab\cdot\,,\cdot\rab$]     A non-degenerate symmetric
bilinear form on $\gt{u}$ or on $\Cl(\gt{u})$
(\S\ref{sec:cliffalgs:struc}), or on the $\Cl(\gt{u})$-module
$\Cm(\gt{u})_E$ (\S\ref{sec:cliffalgs:mods}). In the case that
$\gt{u}$ is real all of these forms will be positive definite.

\item[$\lab\cdot|\cdot\rab$]  A non-degenerate symmetric bilinear
form on $A(\gt{u})_{\Theta}$ (\S\ref{sec:cliffalgs:SVOAs}).
\end{list}
\end{small}
From \S\ref{sec:strucafn} onwards we restrict to the case that
\begin{small}
\begin{list}{}{     \itemsep -3pt
                    \topsep 3pt
                         }
\item[$\gt{u}=\gt{l}$] is a $24$ dimensional real vector space
with positive definite symmetric bilinear form
$\lab\cdot\,,\cdot\rab$;
\item[$\fset=\Omega$] is an ordered set with $24$ elements
indexing an orthonormal basis $\mc{E}=\{e_i\}_{i\in\Omega}$ for
$\gt{l}$;
\item[$E=G$] is a subgroup of $\Sp(\gt{l})$ such that
$\mc{C}(G)$ is a copy of the Golay code $\mc{G}$;
\end{list}
\end{small}
since this is the situation that will be relevant for the
construction of $\afn$.

\section{SVOA structure}\label{sec:PreSVOAs}
%\subsubsection{SVOAs}

\subsection{SVOAs}\label{sec:SVOAstruc:SVOAs}

Suppose that $U=U_{\bar{0}}\oplus U_{\bar{1}}$ is a super vector
space over $\FF$. %Define $s:\ZZ/2\times\ZZ/2\to\CCm$ to be the
%bimultiplicative function such that $s(\bar{1},\bar{1})=-1$ and
%all other values are $1$. We write $s(u,v)$ for $s(\gamma,\delta)$
%when $u\in U_{\gamma}$ and $v\in U_{\delta}$ are
%$\ZZ/2$-homogeneous elements in $U$.
For an {\em SVOA structure} on $U$ we require the following data.
\begin{itemize}
\item{\em Vertex operators:}   a map $U\otimes U\to
U((z))$ denoted $u\otimes v\mapsto Y(u,z)v$ such that, when we write
$Y(u,z)v=\sum_{\ZZ}u_nvz^{-n-1}$, we have $u_nv\in
U_{\gamma+\delta}$ when $u\in U_{\gamma}$ and $v\in U_{\delta}$, and
such that $Y(u,z)v=0$ for all $v\in U$ implies $u=0$.
\item{\em Vacuum:}   a distinguished vector $\vac\in U_{\bar{0}}$
such that $Y(\vac,z)u=u$ for $u\in U$, and $Y(u,z)\vac|_{z=0}=u$.
\item{\em Conformal element:}   a distinguished vector $\cas\in
U_{\bar{0}}$ such that the operators $L(n)=\cas_{n+1}$ furnish a
representation of the Virasoro algebra on $U$(c.f.
\S\ref{sec:PreN=1SVOAs}).
\end{itemize}
This data furnishes an SVOA structure on $U$ just when the
following axioms are satisfied.
\begin{enumerate}
\item{\em Translation:}   for $u\in U$ we have
$[L(-1),Y(u,z)] =D_zY(u,z)$.
\item{\em Jacobi Identity:}   the following Jacobi
identity is satisfied for $\ZZ/2$ homogeneous $u,v\in U$.
\begin{gather}
    \begin{split}
        &
        z_0^{-1}\delta\left(\frac{z_1-z_2}{z_0}\right)
        Y(u,z_1)Y(v,z_2)\\
    &-(-1)^{|u||v|}%s(u,v)
        z_0^{-1}\delta\left(\frac{z_2-z_1}{-z_0}\right)
        Y(v,z_2)Y(u,z_1)
        \\
    &=
        z_2^{-1}\delta\left(\frac{z_1-z_0}{z_2}\right)
        Y(Y(u,z_0)v,z_2)
    \end{split}
\end{gather}
when $u\in U_{\gamma}$ and $v\in U_{\delta}$.
\item{\em $L(0)$-grading:}  the action of $L(0)$ on $U$ is
diagonalizable with rational eigenvalues bounded below, by $-N$
say, and thus defines a $\QQ_{>-N}$-grading $U=\coprod_n U_n$ on
$U$. This grading is such that the $L(0)$-homogeneous subspaces
$U_n=\{u\in U\mid L(0)u=nu\}$ are finite dimensional.
\end{enumerate}
In the case that $U=U_{\bar{0}}$ we are speaking of ordinary VOAs.

An SVOA $U$ is said to be {\em $C_2$-cofinite} in the case that
$U_{-2}U=\{u_{-2}v\mid u,v\in U\}$ has finite codimension in $U$.
Following \cite{HohnPhD} we say that an SVOA $U$ is {\em nice} when
it is $C_2$-cofinite, the eigenvalues of $L(0)$ are non-negative and
contained in $\hZZ$, and the degree zero subspace $U_0$ is one
dimensional and spanned by the vacuum vector. All the SVOAs we
consider in this paper will be nice SVOAs.

By definition the coefficients of $Y(\cas,z)$ define a
representation of the Virasoro algebra on $U$ (c.f.
\S\ref{sec:PreN=1SVOAs}). Let $c\in\CC$ be such that the central
element of the Virasoro algebra acts as $c{\rm Id}$ on $U$. Then $c$
is called the {\em rank} of $U$, and we denote it by ${\rm
rank}(U)$.

\subsection{SVOA Modules}\label{sec:SVOAMods}

For an SVOA module over $U$ we require a $\ZZ/2$-graded vector space
$M=M_{\bar{0}}\oplus M_{\bar{1}}$ and a map $U\otimes M\to M((z))$
denoted $u\otimes v\mapsto Y^M(u,z)v$ such that, when we write
$Y^M(u,z)v = \sum_{\ZZ} u^M_nvz^{-n-1}$, we have $u^M_nv \in
M_{\gamma+\delta}$ when $u\in U_{\gamma}$ and $v\in M_{\delta}$. The
pair $(M,Y^M)$ is called an {\em admissible $U$-module} just when
the following axioms are satisfied.
\begin{enumerate}
\item{\em Vacuum:} the operator $Y^M(\vac,z)$ is the identity on $M$.
\item{\em Jacobi Identity:}   for $\ZZ/2$-homogeneous $u,v\in U$
we have
\begin{gather}
    \begin{split}
        &
        z_0^{-1}\delta\left(\frac{z_1-z_2}{z_0}\right)
        Y^M(u,z_1)Y^M(v,z_2)\\
    &-(-1)^{|u||v|}%s(u,v)
        z_0^{-1}\delta\left(\frac{z_2-z_1}{-z_0}\right)
        Y^M(v,z_2)Y^M(u,z_1)
        \\
    &=
        z_2^{-1}\delta\left(\frac{z_1-z_0}{z_2}\right)
        Y^M(Y(u,z_0)v,z_2)
    \end{split}
\end{gather}
where $u\in U_{\gamma}$, $v\in U_{\delta}$.
\item{\em Grading:}  The space $M$ carries a $\hZZ$-grading
$M=\coprod_rM(r)$ bounded from below such that $u_nM(r)\subset
M(m+r-n-1)$ for $u\in U_m$.
\end{enumerate}
We say that an admissible $U$-module $(M,Y^M)$ is an {\em ordinary
$U$-module} if there is some $h\in \CC$ such that the $\hZZ$-grading
$M=\coprod_r M(r)$ satisfies $M(r)=\{m\in M\mid L(0)m=(r+h)m\}$, and
the spaces $M(r)$ are finite dimensional for each $r\in\hZZ$. In
this case we set $M_{h+r}=M(r)$ and write $M=\coprod_{r} M_{h+r}$.

All the SVOA modules that we consider will be ordinary modules. We
understand that unless otherwise qualified, the term module shall
mean ordinary module.

\begin{rmk}
The significance of the notion of admissible module is the
result of \cite{DonLiMasTwRepsVOAs} that if every admissible
module over a VOA is completely reducible, then there are finitely
many irreducible admissible modules up to equivalence, and every
finitely generated admissible module is an ordinary module.
\end{rmk}

A VOA is called {\em rational} if each of its admissible modules are
completely reducible. We define an SVOA to be {\em rational} if its
even sub-VOA is rational. An SVOA $U$ is called {\em simple} if it
is irreducible as a module over itself. We say that a rational SVOA
is {\em self-dual} if it is simple, and has no irreducible modules
other than itself.

\subsection{Intertwiners}

Let $V$ be a VOA, and let $(M^i,Y^i)$, $(M^j,Y^j)$ and $(M^k,Y^k)$
be three $V$-modules. Suppose given a map $M^i\otimes M^j\to
M^k\{z\}$, and employ the notation $u\otimes v\mapsto
Y^{ij}_k(u,z)v=\sum_{n\in\QQ}u_nvz^{-n-1}$. We assume that for any
$u\in M^i$ and $v\in M^j$ we have $u_nv=0$ for $n$ sufficiently
large, and we assume that $Y^{ij}_k(u,z)=0$ only for $u=0$. Then the
map $Y^{ij}_k$ is called an {\em intertwining operator of type
$\binom{k}{i\;j}$} just when the following axioms are satisfied.
\begin{enumerate}
\item{\em Translation:}   for $u\in M^i$ we have
$Y^{ij}_k(L^i(-1)u,z) =D_zY^{ij}_k(u,z)$.
\item{\em Jacobi Identity:}   for $u\in V$, $v\in M^i$ and $w\in
M^j$ we have
\begin{gather}
    \begin{split}
        &z_0^{-1}\delta\left(\frac{z_1-z_2}{z_0}\right)
        Y^k(u,z_1)Y^{ij}_k(v,z_2)w\\
    &-z_0^{-1}\delta\left(\frac{z_2-z_1}{-z_0}\right)
        Y^{ij}_k(v,z_2)Y^j(u,z_1)w
        \\
    &=
        z_2^{-1}\delta\left(\frac{z_1-z_0}{z_2}\right)
        Y^{ij}_k(Y^i(u,z_0)v,z_2)w
    \end{split}
\end{gather}
\end{enumerate}

\subsection{Twisted SVOA modules}\label{sec:SVOATwMods}

Any SVOA $U$ say, admits an order two involution which is the
identity on $U_{\bar{0}}$ and acts as $-1$ on $U_{\bar{1}}$. We
refer to this involution as the {\em canonical involution} on $U$,
and denote it by $\sigma$.
%The notion of canonically twisted $U$-module is formulated as
%follows.

For a structure of {\em canonically twisted} or {\em
$\sigma$-twisted} $U$-module on a vector space $M$ we require a map
$Y^M:U\otimes M\to M((z^{1/2}))$ such that when we write
$Y^M(u,z)=\sum_{n}u^M_nz^{-n-1}$ for $\ZZ/2$-homogeneous $u\in U$,
then $u^M_n=0$ unless $u\in U_{\bar{k}}$ and $n\in
\ZZ+\tfrac{k}{2}$. The pair $(M,Y^M)$ is called an {\em admissible
canonically twisted $U$-module} when it satisfies just the same
axioms as for untwisted admissible $U$-modules except that we modify
the Jacobi identity axiom and the grading condition as follows.
\begin{enumerate}
\item[*]{\em Twisted Jacobi identity:}
For $\ZZ/2$-homogeneous $u,v\in U$ we require that
\begin{gather}
    \begin{split}
        &z_0^{-1}\delta\left(\frac{z_1-z_2}{z_0}\right)
        Y^M(u,z_1)Y^M(v,z_2)\\
    &-(-1)^{|u||v|}%s(u,v)
        z_0^{-1}\delta\left(\frac{z_2-z_1}{-z_0}\right)
        Y^M(v,z_2)Y^M(u,z_1)
        \\
    &=z_2^{-1}\delta\left(\frac{z_1-z_0}{z_2}\right)
    \left(\frac{z_1-z_0}{z_2}\right)^{-k/2}
        Y^M(Y(u,z_0)v,z_2)
        %,\qquad\text{ when $u\in U_{\bar{k}}$.}
    \end{split}
\end{gather}
when $u\in U_{\bar{k}}$.
\item[*]{\em Twisted grading:}
The space $M$ carries a $\ZZ$-grading $M=\coprod_rM(r)$ bounded from
below such that $u_nM(r)\subset M(m+r-n-1)$ for $u\in U_m$.
\end{enumerate}
We say that an admissible canonically twisted $U$-module $(M,Y^M)$
is an {\em ordinary canonically twisted $U$-module} if there is some
$h\in \CC$ such that the $\ZZ$-grading $M=\coprod_r M(r)$ satisfies
$M(r)=\{m\in M\mid L(0)m=(r+h)m\}$, and the spaces $M(r)$ are finite
dimensional for each $r\in\ZZ$. In this case we set $M_{h+r}=M(r)$
and write $M=\coprod_{r} M_{h+r}$.

As in \S\ref{sec:SVOAMods}, we convene that unless otherwise
qualified, the term ``canonically twisted module'' shall mean
``ordinary canonically twisted module''.

A canonically twisted module $(M,Y^M)$ over an SVOA $V$ is called
{\em $\sigma$-stable} if it admits an action by $\sigma$ compatible
with that on $V$, so that we have $\sigma Y^M(u,z)v=Y^M(\sigma
u,z)\sigma v$ for $u\in V$, $v\in M$.

An important result we will make use of is the following.
\begin{thm}[\cite{DonZhaMdltyOrbSVOA}]\label{thm:sdVOAhastwmod}
If $V$ is a self-dual rational $C_2$-cofinite SVOA then $V$ has a
unique irreducible $\sigma$-stable $\sigma$-twisted module.
\end{thm}
\begin{rmk}
Recall from \S\ref{sec:SVOAMods} that we define an SVOA to be
rational in case it's even sub-VOA is rational. This is a stronger
condition than the notion of rationality in
\cite{DonZhaMdltyOrbSVOA}, and an SVOA that is rational in our sense
is both rational and $\sigma$-rational in the sense of
\cite{DonZhaMdltyOrbSVOA}.
\end{rmk}

For $V$ satisfying the hypotheses of Theorem~\ref{thm:sdVOAhastwmod}
we will write $(V_{\sigma},Y_{\sigma})$ for the unique irreducible
canonically twisted $V$-module this theorem guarantees.

\subsection{$N=1$ SVOAs}\label{sec:PreN=1SVOAs}

The Neveu--Schwarz superalgebra is the Lie superalgebra spanned by
the symbols $L_m$, $G_{m+1/2}$, and ${\bf{c}}$, for $m\in\ZZ$, and
subject to the following relations \cite{KacWanSVOAs}.
\begin{gather}\label{SVirRelf}
    [L_m,L_n]=(m-n)L_{m+n}+\frac{m^3-m}{12}\delta_{m+n,0}\bf{c},\\
    \left[ G_{m+1/2},L_n\right]=
        \left(m+\frac{1}{2}-\frac{n}{2}\right)G_{m+n+1/2},\\
    \left\{G_{m+1/2},G_{n-1/2}\right\}=2L_{m+n}+
         \frac{m^2+m}{3}\delta_{m+n,0}\bf{c},\\%\quad
    \left[L_m,\bf{c}\right]=\left[G_{m+1/2},\bf{c}\right]=0.
\end{gather}
Note that this algebra is generated by the $G_{m+1/2}$ for
$m\in\ZZ$. The subalgebra generated by the $L_n$ is the Virasoro
algebra.

Suppose that $U$ is an SVOA. We say that $U$ is an $N=1$ SVOA if
there is an element $\scas\in U_{3/2}$ such that the operators
$G(n+\tfrac{1}{2})=\tau_{n+1}$ generate a representation of the
Neveu--Schwarz superalgebra on $U$. We refer to such an element
$\tau$ as a {\em superconformal element} for $U$. In general, a
superconformal element for an SVOA $U$ is not unique and may not
even exist, but for any superconformal $\tau$ we have that
$\tfrac{1}{2}\tau_{0}\tau =\tfrac{1}{2}G(-\tfrac{1}{2})\tau$ is a
conformal element. That is to say, the coefficients of
$Y(\tfrac{1}{2}G(-\tfrac{1}{2})\tau,z)$ generate a representation of
the Virasoro algebra on $U$. Given an SVOA $U$, we will always
assume that any superconformal element for $U$ is chosen so that
$\tfrac{1}{2}G(-\tfrac{1}{2})\tau=\cas$.

Note that the commutation relations for the Neveu--Schwarz
superalgebra of central charge $c$ are equivalent to the following
operator product expansions \cite{DixGinHarBB}.
\begin{gather}
    Y(\cas,z_1)Y(\cas,z_2)\sim\frac{c/2}{(z_1-z_2)^{4}}+
    \frac{2Y(\cas,z_2)}{(z_1-z_2)^{2}}+
    \frac{D_{z_2}Y(\cas,z_2)}{(z_1-z_2)}\\
    Y(\cas,z_1)Y(\scas,z_2)\sim
    \frac{3}{2}\frac{Y(\scas,z_2)}{(z_1-z_2)^{2}}+
    \frac{D_{z_2}Y(\scas,z_2)}{(z_1-z_2)}\\
    Y(\scas,z_1)Y(\scas,z_2)\sim\frac{2c/3}{(z_1-z_2)^{3}}+
    \frac{2Y(\cas,z_2)}{(z_1-z_2)}
\end{gather}
Consequently, we have the following
\begin{prop}\label{prop:SConfCriterion}
Suppose $U$ is an SVOA with conformal element $\cas$ and central
charge $c$. Then $\tau\in (U)_{3/2}$ is a superconformal element for
$U$ so long as $\tau_{2}\tau=\tfrac{2}{3}c\vac$, $\tau_{1}\tau=0$
and $\tau_{0}\tau={2}\cas$.
\end{prop}
Given an $N=1$ SVOA $U$ with superconformal element $\tau$ and
conformal element $\omega=\tfrac{1}{2}\tau_0\tau$, we write
$\Aut_{SVOA}(U)$ for the group of automorphisms of the SVOA
structure on $U$, and we write $\Aut(U)$ for group of automorphisms
of the $N=1$ SVOA structure on $U$. That is, for $U$ an $N=1$ SVOA,
$\Aut(U)$ denotes the subgroup of $\Aut_{SVOA}(U)$ comprised of
automorphisms that fix $\tau$.

\subsection{Adjoint operators}\label{Sec:AdjOps}

Suppose that $U$ is a nice SVOA.  For $M$ a module over $U$, let
$M'$ denote the restricted dual of $M$, and let
$\langle\cdot\,,\cdot\rangle_{M}$ be the natural pairing $M'\times
M\to \CC$. We define the adjoint vertex operators $Y':U\otimes M'\to
M'\{z\}$ by requiring that
\begin{gather}
    \langle Y'(u,z)w',w\rangle_M
    =(-1)^n\left\langle w',Y(e^{zL(1)}
        z^{-2L(0)}u,z^{-1})w\right\rangle_M
\end{gather}
for $u\in U_{n-1/2}\oplus U_n$ with $n\in \ZZ$, where $w'\in M'$
and $w\in M$. As in \cite{FHL} we have
\begin{prop}
The map $Y'$ equips $M'$ with a structure of $U$-module.
\end{prop}
Let $\langle\cdot| \cdot\rangle:U\otimes U\to\CC$ be a bilinear
form on $U$ such that $\langle U_m\mid U_n\rangle\subset\{0\}$
unless $m=n$. Then there is a unique grading preserving linear map
$\varphi:U\to U'$ determined by the formula
$\langle\varphi(u),v\rangle_U=\langle u\mid  v\rangle$, and
$\varphi$ is a $U$-module equivalence if and only if
\begin{gather}
    \left\langle Y(u,z)w_1\mid w_2\right\rangle
    %=\langle w_1\mid Y'(u,z)w_2\rangle
    =\left\langle w_1\mid Y(e^{zL(1)}
        z^{-2L(0)}u,z^{-1})w_2\right\rangle
\end{gather}
for all $u$, $w_1$ and $w_2$ in $U$. When this identity is
satisfied we say that the bilinear form $\langle\cdot|
\cdot\rangle$ is an invariant form for $U$.
\begin{prop}\label{InvFrmIsSymm}
Suppose that $\langle\cdot| \cdot\rangle$ is an invariant form for
$U$. Then it is symmetric.
\end{prop}
Just as in the untwisted case, we can define the adjoint
canonically twisted vertex operators. For $(M,Y)$ a canonically
twisted $U$ module, we define operators $Y':U\otimes M'\to
M'\{z\}$ by requiring that
\begin{gather}
    \langle Y'(u,z)w',w\rangle_{M}
    =(-1)^n\langle w',Y(e^{zL(1)}
        z^{-2L(0)}u,z^{-1})w\rangle_{M}
\end{gather}
for $u\in U_{n-1/2}\oplus U_n$ with $n\in\ZZ$, where $w'\in M'$,
$w\in M$, and $\langle\cdot\,,\cdot\rangle_{M}$ is the natural
pairing $M'\times M\to\CC$. As in the untwisted case, we have
\begin{prop}
The map $Y'$ equips $M'$ with a structure of canonically twisted
$U$-module.
\end{prop}

\subsection{Lattice SVOAs}\label{sec:SVOAstruc:LattSVOAs}

Let $L$ be a positive definite integral lattice, and recall that
$\FF$ denotes $\RR$ or $\CC$. Then the following standard
procedure associates to $L$, an SVOA defined over $\FF$ and we
shall denote it by $_{\FF}V_L$. We refer the reader to \cite{FLM}
for more details.

Let $_{\FF}\gt{h}=\FF\otimes_{\ZZ} L$ and let $_{\FF}\hat{\gt{h}}$
denote the Heisenberg Lie algebra described by
\begin{gather}
    _{\FF}\hat{\gt{h}}=\coprod_{n\in\ZZ}
         \,_{\FF}\gt{h}\otimes t^n\oplus \FF c,\quad
        [h(m),h'(n)]=m\langle h,h'\rangle\delta_{m+n,0},\quad
        [h(m),c]=0.
\end{gather}
We denote by $_{\FF}\hat{\gt{b}}$ and $_{\FF}\hat{\gt{b}}'$, the
(commutative) subalgebras of $_{\FF}\hat{\gt{h}}$ given by
\begin{gather}
    _{\FF}\hat{\gt{b}}=\coprod_{n\in\ZZ_{\geq 0}}\,
        _{\FF}\gt{h}\otimes t^n\oplus \FF c,\qquad
    _{\FF}\hat{\gt{b}}'=\coprod_{n\in\ZZ_{< 0}}\,
        _{\FF}\gt{h}\otimes t^n.
\end{gather}
Let $\hat{L}$ be the unique up to equivalence extension of $L$ by
a group $\langle\kappa\rangle$ with generator $\kappa$ of order
two, such that the commutators in $\hat{L}$ satisfy
\begin{gather}
    aba^{-1}b^{-1}=\kappa^{\langle\bar{a},\bar{b}\rangle+
        \langle\bar{a},\bar{a}\rangle
         \langle\bar{b},\bar{b}\rangle}
\end{gather}
where $a\mapsto \bar{a}$ denotes the natural homomorphism
$\hat{L}\to L$. We have the following short exact sequence.
\begin{gather}
    1\to\langle\kappa\mid\kappa^2=1\rangle\to\hat{L}
        \xrightarrow{-}L\to 1
\end{gather}
Let $\FF\{L\}$ denote the $\hat{L}$-module obtain by factoring
$\FF\hat{L}$ by the subalgebra generated by $\kappa+1$. We write
$\iota(a)$ for the image of $a\in\hat{L}$ in $\FF\{L\}$ under the
composition of maps $\hat{L}\hookrightarrow\FF\hat{L}\to\FF\{L\}$.
The space $\FF\{L\}$ is again an algebra, and is linearly
isomorphic to $\FF L$. The algebra $\FF\{L\}$ may be equipped with
a structure of $_{\FF}\hat{\gt{b}}$-module as follows.
\begin{gather}
    h(m)\cdot \iota(a)=0\;\;\text{for $m>0$},\qquad
    h(0)\cdot \iota(a)=\langle h,\bar{a}\rangle \iota(a),\qquad
    c\cdot \iota(a)=\iota(a).
\end{gather}
As an $_{\FF}\hat{\gt{h}}$-module, $_{\FF}V_L$ is defined to be
that induced from the $_{\FF}\hat{\gt{b}}$-module structure on
$\FF\{L\}$.
\begin{gather}
    _{\FF}V_L=U(_{\FF}\hat{\gt{h}})
        \otimes_{U(_{\FF}\hat{\gt{b}})}\FF\{L\}
\end{gather}
We identify $\FF\{L\}$ with the subspace $1\otimes \FF\{L\}$ of
$_{\FF}V_L$, and we set $\vac=1\otimes \iota(1)$. There is a
natural isomorphism of $_{\FF}\hat{\gt{b}}'$-modules,
$_{\FF}V_L\simeq S(_{\FF}\hat{\gt{b}}')\otimes\FF\{ L\}$.

Now we define the vertex operators on $_{\FF}V_L$. Let
$h\in\,_{\FF}\gt{h}$. We define generating functions $h(z)$ and
$l(h,z)$ of operators on $_{\FF}V_L$ by
\begin{gather}
    h(z)=\sum_{n\in\ZZ}h(n)z^{-n-1},\qquad
    l(h,z)=\sum_{n\in\ZZ,\; n\neq 0}\frac{h(n)}{-n}z^{-n}
\end{gather}
Then for $h\in\, _{\FF}\gt{h}$ and $a\in\hat{L}$, the vertex
operators associated to $\iota(a)=1\otimes \iota(a)$ and
$h(-n-1)=h(-n-1)\otimes \iota(1)$ are given by
\begin{gather}
    Y(h(-n-1),z)=D_z^{(n)}h(z),\qquad
    Y(\iota(a),z)=:\exp\left(l(\bar{a},z)\right)az^{\bar{a}(0)}:
\end{gather}
respectively, where the colons denote the Bosonic normal ordering:
that all operators $h(m)$ with $m\geq 0$ be commuted to the right
of all other operators. The remaining vertex operators are
determined by the requirement that $Y$ be a linear map, and that
\begin{gather}
    Y(u_{-1}v,z)=:Y(u,z)Y(v,z):\quad\text{for
    $u,v\in\,_{\FF}V_L$.}
\end{gather}
If $\{h_i\}$ is a basis for $_{\FF}\gt{h}$ and $\{h'_i\}$ is the
dual basis, let $\cas=\tfrac{1}{2}\sum h_i(-1)h_i'(-1)$. Then
$\cas$ is independent of the choice of basis, and we have the
following
\begin{thm}[\cite{FLM}]
The quadruple $(_{\FF}V_L,Y,\vac,\cas)$ is an SVOA and the rank of
$_{\FF}V_L$ coincides with the rank of $L$.
\end{thm}
We refer to $_{\FF}V_L$ as the SVOA over $\FF$ associated to $L$ via
the standard construction. The superspace decomposition of
$_{\FF}V_L$ is given by $_{\FF}V_L=\,_{\FF}V_{L_0}\oplus
\,_{\FF}V_{L_1}$ where $L_0$ is the sublattice of $L$ consisting of
elements with even norm squared, and $L_1$ is the (unique) coset of
$L_0$ in $L$, consisting of elements with odd norm squared.

\subsubsection{The involution $\ogi$ for $_{\FF}V_L$}

The lattice $L$ admits a non-trivial involution that acts by
$\alpha\mapsto-\alpha$ for $\alpha\in L$. This involution lifts
naturally to automorphisms of $_{\FF}\gt{h}$ and
$_{\FF}\hat{\gt{h}}$ and hence to $U(_{\FF}\hat{\gt{h}})$. We
denote it by $\ogi$. We may define an automorphism of $\FF\{L\}$
by
\begin{gather}
    \iota(a)\mapsto (-1)^{\langle\bar{a},\bar{a}\rangle/2+
        \langle\bar{a},\bar{a}\rangle^2/2}\iota(a^{-1})
\end{gather}
for $a\in\hat{L}$. We denote it also by $\ogi$. Recalling that
$_{\FF}V_L$ was constructed as
\begin{gather}
    _{\FF}V_L=U(_{\FF}\hat{\gt{h}})
        \otimes_{_{\FF}\hat{\gt{b}}}\FF\{L\}
\end{gather}
we may define an automorphism of $_{\FF}V_L$ by $\ogi\otimes\ogi$
where the $\ogi$ on the left is that defined for the left tensor
factor of $_{\FF}V_L$, and the one on the right is that just defined
on the right tensor factor. Since all these automorphisms may be
regarded as lifts of $-1$, we denote $\ogi\otimes\ogi$ also by
$\ogi$. Then $\ogi$ is an automorphism of the VOA structure on
$_{\FF}V_L$, and we may refer to it as a lift to $\Aut(_{\FF}V_L)$
of the $-1$ symmetry on $L$.
\begin{rmk}
It is a result of \cite{DonGriHohFVOAsMM} that all lifts of $-1$ are
conjugate under the action of $\Aut(_{\FF}V_L)$.
\end{rmk}
The space $_{\FF}V_L$ decomposes into eigenspaces for the action
of $\ogi$, and we express this decomposition as
$_{\FF}V_L=\,_{\FF}V_L^0\oplus\,_{\FF}V_L^1$ where $_{\FF}V_L^k$
is the eigenspace with eigenvalue $(-1)^k$ for the action of
$\ogi$.

\subsubsection{Real form for $_{\FF}V_L$}\label{RealFormLatSVOA}

The adjoint operators on $_{\FF}V_L$ determine an invariant
bilinear form, which is given by
\begin{gather}
    \langle u\mid v\rangle \vac
    ={\rm Res}_{z=0}\,z^{-1}(-1)^nY(e^{zL(1)}
        z^{-2L(0)}u,z^{-1})v
\end{gather}
when $u$ is in $(_{\FF}V_L)_{n-1/2}$ or $(_{\FF}V_L)_{n}$ for some
$n\in\ZZ$. Consider the case that $\FF=\RR$. Then the form
$\langle\cdot\mid\cdot\rangle$ is not positive definite on
$_{\RR}V_L$. In fact, the form is positive definite on the
subspace $_{\RR}V_L^0$, and is negative definite on $_{\RR}V_L^1$.
Suppose now that we view $_{\RR}V_L$ as embedded in $_{\CC}V_L$
curtesy of the the natural inclusion $\RR\hookrightarrow \CC$. Let
us set ${V}_L$ to be the $\RR$ subspace of $_{\CC}V_L$ described
by
\begin{gather}
    {V}_L=\,_{\RR}V_L^0\oplus \ii\,_{\RR}V_L^1
        =\{u+\ii v\mid u\in\,_{\RR}V_L^0,\,
            v\in\,_{\RR}V_L^1\}
\end{gather}
Then ${V}_L$ closes under the vertex operators $Y$ associated to
$_{\CC}V_L$, and the form $\langle\cdot\mid\cdot\rangle$ restricts
to be positive definite on ${V}_L$. From now on we write ${V}_L$
for the real VOA with positive definite bilinear form obtained in
this way, by restricting the form and vertex operator algebra
structure from $_{\CC}V_L$.

%%%%%%%%%%%%%%%%%%%%%%%%%%%%%%%%%%%%%%%%%%%%%%%%%%%%%%%%%%%

\section{Clifford algebras}\label{sec:cliffalgs}

The construction of SVOAs that we will use is based on the structure
of Clifford algebra modules. In this section we recall some basic
properties of Clifford algebras and we exhibit a construction of
modules over finite dimensional Clifford algebras using doubly-even
linear binary codes (see \S\ref{Sec:Notation}). In
\S\ref{sec:cliffalgs:SVOAs} we recall the construction of SVOA
module structure on modules over certain infinite dimensional
Clifford algebras.

\subsection{Clifford algebra structure}\label{sec:cliffalgs:struc}

In this section $\FF$ denotes either $\RR$ or $\CC$. For $\gt{u}$
an $\FF$-vector space with non-degenerate symmetric bilinear form
$\langle \cdot\,,\cdot\rangle$ we write $\Cl(\gt{u})$ for the
Clifford algebra over $\FF$ generated by $\gt{u}$. More precisely,
we set $\Cl(\gt{u})=T(\gt{u})/I(\gt{u})$ where $T(\gt{u})$ is the
tensor algebra of $\gt{u}$ over $\FF$ with unit denoted ${\bf 1}$,
and $I(\gt{u})$ is the ideal of $T(\gt{u})$ generated by all
expressions of the form $u\otimes u+\langle u,u\rangle{\bf 1}$ for
$u\in \gt{u}$. The natural algebra structure on $T(\gt{u})$
induces an associative algebra structure on $\Cl(\gt{u})$. The
vector space $\gt{u}$ embeds in $\Cl(\gt{u})$, and when it is
convenient we identify $\gt{u}$ with its image in $\Cl(\gt{u})$.
We also write $\alpha$ in place of $\alpha{\bf
1}+I(\gt{u})\in\Cl(\gt{u})$ for $\alpha\in\FF$ when no confusion
will arise. For $u\in \gt{u}$ we have the relation $u^2=-|u|^2$ in
$\Cl(\gt{u})$. Polarization of this identity yields
$uv+vu=-2\langle u,v\rangle$ for $u,v\in \gt{u}$.

The linear transformation on $\gt{u}$ which is $-1$ times the
identity map lifts naturally to $T(\gt{u})$ and preserves
$I(\gt{u})$, and hence induces an involution on $\Cl(\gt{u})$ which
we denote by $\ogi$. The map $\ogi$ is often referred to as the {\em
parity involution}. We have $\ogi(u_1\cdots u_k)=(-1)^ku_1\cdots
u_k$ for $u_1\cdots u_k\in\Cl(\gt{u})$ with $u_i\in \gt{u}$, and we
write $\Cl(\gt{u})=\Cl(\gt{u})^0\oplus \Cl(\gt{u})^1$ for the
decomposition into eigenspaces for $\ogi$. Define a bilinear form on
$\Cl(\gt{u})$, denoted $\langle\cdot\,,\cdot\rangle$, by setting
$\langle{\bf 1},{\bf 1}\rangle=1$, and requiring that for $u\in
\gt{u}$, the adjoint of left multiplication by $u$ is left
multiplication by $-u$. Then the restriction of
$\langle\cdot\,,\cdot\rangle$ to $\gt{u}$ agrees with the original
form on $\gt{u}$.
\begin{gather}
    \langle u,u\rangle=-\langle{\bf 1},u^2\rangle
        =-\langle{\bf 1},-|u|^2{\bf 1}\rangle
        =|u|^2
\end{gather}
More generally, the adjoint of $u=u_1\cdots u_k$ for $u_i\in \gt{u}$
is $(-1)^ku_k\cdots u_1$. The {\em main anti-automorphism} of
$\Cl(\gt{u})$ is the map we denote $\alpha$, which acts by sending
$u_1\cdots u_k$ to $u_k\cdots u_1$ for $u_i\in\gt{u}$.

\subsection{Spin groups}\label{sec:cliffalgs:spin}

Let us write $\Cl(\gt{u})^{\times}$ for the group of invertible
elements in $\Cl(\gt{u})$. For $x\in\Cl(\gt{u})^{\times}$ and
$a\in \Cl(\gt{u})$, we set $x(a)=xax^{-1}$. %Note that $u\in
%\gt{u}$ is invertible in $\Cl(\gt{u})$ with inverse $-u/|u|^2$ so
%long as $u\neq 0$.
We will define the Pinor and Spinor groups associated to $\gt{u}$
slightly differently according as $\gt{u}$ is real or complex: in
the case that $\gt{u}$ is real, we define the Pinor group
$\textsl{Pin}(\gt{u})$ to be the subgroup of $\Cl(\gt{u})^{\times}$
comprised of elements $x$ such that $x(\gt{u})\subset \gt{u}$ and
$\alpha(x)x=\pm 1$; in the case that $\gt{u}$ is complex we define
$\textsl{Pin}(\gt{u})$ to be the set of $x\in\Cl(\gt{u})^{\times}$
such that $x(\gt{u})\subset\gt{u}$ and $\alpha(x)x=1$. In both cases
we define the Spinor group by setting
$\Sp(\gt{u})=\textsl{Pin}(\gt{u})\cap\Cl(\gt{u})^0$.

Let $x\in \textsl{Pin}(\gt{u})$. Then we have $\langle
x(u),x(v)\rangle=\langle u,v\rangle$ for $u,v\in \gt{u}$, and thus
the map $x\mapsto x(\cdot)$, which has kernel $\pm {\bf 1}$,
realizes the Pinor group as a double cover of ${O}(\gt{u})$. (If
$u\in\gt{u}$ and $\langle u,u\rangle=1$, then $u(\cdot)$ is the
orthogonal transformation of $\gt{u}$ which is minus the reflection
in the hyperplane orthogonal to $u$.) The image of $\Sp(\gt{u})$
under the map $x\mapsto x(\cdot)$ is just the group $SO(\gt{u})$.

In the case that $\gt{u}$ is real with definite bilinear form, we
have $\alpha(x)x=1$ for all $x\in\Sp(\gt{u})$, and the group
$\Sp(\gt{u})$ is generated by the (well-defined) expressions
$\exp(\lambda e_i e_j)\in\Cl(\gt{u})^{\times}$ for $\lambda\in\RR$
and $\{e_i\}$ an orthonormal basis of $\gt{u}$. The Spinor group of
the complexified space $_{\CC}\gt{u}$ is then generated by the
$\exp(\lambda e_ie_j)$ for $\lambda\in\CC$.

\subsection{Clifford algebra modules}\label{sec:cliffalgs:mods}

For the remainder of this section we suppose that $\gt{u}$ is a
finite dimensional real vector space with positive definite
symmetric bilinear form, and also that $\mc{E}=\{e_i\}_{i\in\fset}$
is an orthonormal basis for $\gt{u}$, indexed by a finite set
$\fset$. For $S=(s_1,\ldots,s_k)\in\fset^{\times k}$, a $k$-tuple of
elements from $\fset$ for any $k$, we write $e_S$ for the element
$e_{s_1}e_{s_2}\cdots e_{s_k}$ in $\Cl(\gt{u})$. We suppose that
$\fset$ is equipped with some ordering. If $S=\{s_1,\ldots,s_k\}$ is
a subset of $\fset$ (that is, an unordered subset), we denote by
$\vec{S}$ the $k$-tuple given by $\vec{S}=(s_1,\ldots,s_k)$ just
when $s_1<\cdots <s_k$, so that $e_{\vec{S}}=e_{s_1}\cdots e_{s_k}$.
We then abuse notation to write $e_{S}$ for $e_{\vec{S}}$. In this
way we obtain an element $e_{S}$ in $\Cl(\gt{u})$ for any
$S\subset\fset$. (We set $e_{\emptyset}={\bf 1}$.) This
correspondence depends on the choice of ordering, but our discussion
will be invariant with respect to this choice. Note that
$e_{S}e_{R}=\pm e_{S+R}$ for any $S,R\subset\fset$, and the set
$\{e_S\mid S\subset\fset\}$ furnishes an orthonormal basis for
$\Cl(\gt{u})$.

We now obtain an $\FF_2^{\fset}$-grading on $\Cl(\gt{u})$ by
decreeing that for $S\subset\fset$, the homogeneous subspace of
$\Cl(\gt{u})$ with degree $S$ is just the $\FF$-span of the vector
$e_{S}$.
\begin{gather}
    \Cl(\gt{u})=\bigoplus_{S\subset\fset}
        \Cl(\gt{u})^S,\;
    \Cl(\gt{u})^S=\FF e_S,\;
    \Cl(\gt{u})^S\Cl(\gt{u})^R\subset\Cl(\gt{u})^{S+R}.
\end{gather}
Since this grading depends on the choice of orthonormal basis
$\mc{E}$, we will refer to it as the $\FF_2^{\mc{E}}$-grading, and
we refer to the homogeneous elements $e_S$ as
$\FF_2^{\mc{E}}$-homogeneous elements. A given subset of
$\Cl(\gt{u})$ is called $\FF_2^{\mc{E}}$-homogeneous if all of its
elements are $\FF_2^{\mc{E}}$-homogeneous.

Suppose that $E$ is an $\FF_2^{\mc{E}}$-homogeneous subgroup of
$\Sp(\gt{u})$. Then $E$ is a union of elements of the form $\pm e_C$
for $C\subset\fset$, and hence is finite, and has exponent four.
Furthermore, the set of $C$ for which $\pm e_C$ is in $E$ determines
a linear binary code on $\fset$. For $E$ an
$\FF_2^{\mc{E}}$-homogeneous subgroup of $\Sp(\gt{u})$, we write
$\mc{C}(E)$ for the associated linear binary code on $\fset$. The
following result is straightforward.
\begin{prop}
Suppose that $-1\notin E$ and $\mc{C}(E)$ is a doubly-even code.
Then the map $E\to \mc{C}(E)$ such that $\pm e_S\mapsto S$ is an
isomorphism of abelian groups, and the sub-algebra of $\Cl(\gt{u})$
generated by $E$ is naturally isomorphic to the group algebra $\RR
E$.
\end{prop}
Suppose now that $E$ is an $\FF_2^{\mc{E}}$-homogeneous subgroup of
$\Sp(\gt{u})$ such that $-1\notin E$ and $\mc{C}(E)$ is a self-dual
doubly-even code. (Note that this forces ${\rm dim}(\gt{u})$ to be a
multiple of eight.) Then we write $\Cm(\gt{u})_E$ for the
$\Cl(\gt{u})$-module defined by
$\Cm(\gt{u})_{E}=\Cl(\gt{u})\otimes_{\RR E}\RR_{1}$ where $\RR_1$
denotes a trivial $E$-module. Let us set $1_E={\bf 1}\otimes 1\in
\Cm(\gt{u})_E$. Then $\Cm(\gt{u})_E$ admits a bilinear form defined
so that $\langle 1_{E},1_E\rangle=1$, and the adjoint to left
multiplication by $u\in\gt{u}\hookrightarrow \Cl(\gt{u})$ is left
multiplication by $-u$.
\begin{prop}
The $\Cl(\gt{u})$-module $\Cm(\gt{u})_E$ is irreducible, and a
vector-space basis for $\Cm(\gt{u})_E$ is naturally indexed by the
elements of the co-code $\mc{C}(E)^{*}$. The bilinear form on
$\Cm(\gt{u})_E$ is non-degenerate.
\end{prop}
\begin{proof}
We have $e_{S+C}{1}=\pm e_{S}{1}$ for any $S\subset\fset$ when
$C\in\mc{C}(E)$. This shows that a basis for $\Cm(\gt{u})_E$ is
indexed by the elements of the co-code
$\mc{C}(E)^{*}=\mc{P}(\fset)/\mc{C}(E)$, and it follows that the
irreducible submodules of $\Cm(\gt{u})_E$ are indexed by the cosets
of $\mc{C}(E)$ in its dual code
$\mc{C}(E)^{\perp}=\{S\in\mc{P}(\fset)\mid |S\cap C|\equiv
0\pmod{2},\;\forall C\in\mc{C}(E) \}$. Since $\mc{C}(E)$ is
self-dual, $\Cm(\gt{u})_E$ is irreducible.
\end{proof}
Note that the vector $1_E\in\Cm(\gt{u})_E$ is such that $g1_E=1_E$
for all $g\in E$, and $\Cm(\gt{u})_E$ is spanned by the $a1_E$ for
$a\in \Cl(\gt{u})$. We have the following
\begin{prop}\label{CliffModEquivsB}
Suppose that $\Cm(\gt{u})_0$ is a non-trivial irreducible
$\Cl(\gt{u})$-module with a vector $1_0$ such that $g1_0=1_0$ for
all $g\in E$. Then $\Cm(\gt{u})_0$ is equivalent to
$\Cm(\gt{u})_E$, and a module equivalence is furnished by the map
$\phi:\Cm(\gt{u})_E\to\Cm(\gt{u})_0$ defined so that $\phi(a1_E)=
a1_0$ for $a\in\Cl(\gt{u})$.
\end{prop}
\begin{proof}
Recall our assumption that $\mc{C}(E)$ is self-dual so that $2d={\rm
dim}(\gt{u})$ is divisible by eight, and $|\mc{C}(E)|=2^{d}$. In
this case we have that $\Cl(\gt{u})$ is isomorphic to
$M_{2^{d}}(\RR)$ so that there is a unique non-trivial irreducible
$\Cl(\gt{u})$-module up to equivalence and it has dimension $2^{d}$.
It suffices then to show that $\phi$ is well defined, and this
follows from the universal mapping property of the induced module
$\Cm(\gt{u})_E$.
%We claim that $\phi$ is well defined. To this end, suppose that
%$a\in\Cl(\gt{u})$ and $a1_E=0$. We will show that $a1_0=0$. Let
%$\varphi_E$ be the map $\Cl(\gt{u})\to\Cm(\gt{u})_E$ given by
%$a\mapsto a1_E$. Then $\varphi_E$ is surjective, and $a1_E=0$ just
%when $a$ is in the kernel of $\varphi_E$. Let us set
%$\mc{C}=\mc{C}(E)$, and let $\mc{C}'$ be a set of representatives
%for the cosets of $\mc{C}$ in $\FF_2^{\fset}$; that is, let
%$\mc{C}'$ be a set of subsets of $\fset$ such that
%$\FF_2^{\fset}/\mc{C}=\cup_{S\in\mc{C}'}S+\mc{C}$. Then ${\rm
%Ker}(\varphi_E)$ contains the set $\{e_S(1-g)\mid S\in\mc{C}',
%g\in E\}$. Since $|\mc{C}'|=|E|=2^{d}$, this set has
%$2^{d}(2^{d}-1)=2^{2d}-2^{d}$ linearly independent vectors. Since
%$\varphi_E$ is surjective and $M_E$ has dimension $2^{d}$, we have
%accounted for all of ${\rm Ker}(\varphi_E)$.
%Define $\varphi_0$ to be the map $\Cl(\gt{u})\to \Cm(\gt{u})_0$
%given by $a\mapsto a1_0$. By hypothesis, any element of the form
%$e_S(1-g)$ for $S\in\mc{C}'$ and $g\in E$ lies in the kernel of
%$\varphi_0$, so the kernel of $\varphi_0$ contains the kernel of
%$\varphi_E$, and ${\rm dim}({\rm Ker}(\varphi_0))\geq
%2^{2d}-2^{d}$. By hypothesis, $\varphi_0$ is surjective, so ${\rm
%dim}(M_0)\leq 2^{d}$. A non-trivial module over $\Cl(\gt{u})$ has
%dimension $2^{d}$ and is unique since up to equivalence so ${\rm
%dim}(M_0)=2^{d}$ and $\varphi_0=\varphi_E$. This shows that
%$a1_0=0$ just when $a1_E=0$, and hence $\phi$ is well-defined and
%necessarily a $\Cl(\gt{u})$-module equivalence.
\end{proof}

We will sometimes wish to replace $\gt{u}$ with its complexification
$_{\CC}\gt{u}$ in the above, in which case we shall understand
$\Cm(_{\CC}\gt{u})_E$ to be the complexification
$_{\CC}\Cm(\gt{u})_E$ of $\Cm(\gt{u})_E$. Then $\Cm(_{\CC}\gt{u})_E$
is an irreducible module over $\Cl(_{\CC}\gt{u})$.

\subsection{Clifford module construction of
SVOAs}\label{sec:cliffalgs:SVOAs}

Let $\gt{u}$ be a finite dimensional vector space over $\RR$ with
positive definite symmetric bilinear form. In this section we review
the construction of SVOA structure on modules over certain infinite
dimensional Clifford algebras associated to $\gt{u}$. The
construction is quite standard and one may refer to \cite{FFR} for
example, for full details. Our setup is somewhat different from that
in \cite{FFR} in that we prefer to be able to work over $\RR$, and
we must therefore use an alternative construction of the canonically
twisted SVOA module, since a polarization of $\gt{u}$ does not exist
in this case. On the other hand, all one requires is an irreducible
module over the (finite dimensional) Clifford algebra $\Cl(\gt{u})$,
and the arguments of \cite{FFR} then go through with only cosmetic
changes.

So that the reader can translate between this section and the
exposition in \cite{FFR}, let us note that in \cite{FFR} they
consider the case that $\gt{a}$ say, is a complex vector space with
non-degenerate bilinear form $\langle\cdot \,,\cdot \rangle_0$, and
a polarization $\gt{a}=\gt{a}^+\oplus\gt{a}^-$ such that
$\gt{a}^{\pm}$ is spanned by vectors $a^{\pm}_i$, which satisfy
$\langle a^{\pm}_i, a^{\mp}_j\rangle_0=\delta_{ij}$. They consider
the Clifford algebra $\Cl_0(\gt{a})$ defined by
$\Cl_0(\gt{a})=T(\gt{a})/I_0(\gt{a})$ where $I_0(\gt{a})$ is the
ideal spanned by elements of the form $uv+vu=\langle u,v\rangle_0
{\bf 1}$. In this section we take $\gt{u}$ to be a real vector space
of even dimension with positive definite bilinear form and
orthonormal basis $\{e_i\}$. Let $d={\rm dim}(\gt{u})/2$, and in the
complexification $_{\CC}\gt{u}$ of $\gt{u}$ set
$a^{\pm}_j=\tfrac{1}{2}(ie_j\mp e_{j+d})$. Then we have an
identification of vector spaces $\gt{a}=\,_{\CC}\gt{u}$. We also
have $\langle a^{\pm}_i,a^{\mp}_j \rangle=-\tfrac{1}{2}\delta_{ij}$
so that $\{a^{\pm}_i,a^{\mp}_j\}=-2\langle
a^{\pm}_i,a^{\mp}_j\rangle = \delta_{ij}$ in $\Cl(_{\CC}\gt{u})$,
and setting $\langle \cdot\,,\cdot\rangle_0
=-2\langle\cdot\,,\cdot\rangle$ we have an identification of
algebras $\Cl_0(\gt{a})=\Cl(_{\CC}\gt{u})$.

\medskip

We now proceed with the construction. We assume that $\gt{u}$ admits
an orthonormal basis $\{e_i\}_{i\in \fset}$ indexed by a finite set
$\fset$. For simplicity we suppose that the dimension of $\gt{u}$ is
divisible by eight so that maximal self-orthogonal doubly-even codes
on $\fset$ are self-dual. Let $\hat{\gt{u}}$ and
$\hat{\gt{u}}_{\ogi}$ denote the infinite dimensional inner product
spaces described as follows.
\begin{gather}
    \hat{\gt{u}}=\coprod_{m\in\ZZ}\gt{u}\otimes t^{m+1/2},\quad
    \hat{\gt{u}}_{\ogi}=\coprod_{m\in\ZZ}\gt{u}\otimes t^m,\\
    \langle u\otimes t^r,v\otimes t^s\rangle
        =\langle u,v\rangle \delta_{r+s,0},
    \; \text{ for $u,v\in\gt{u}$ and $r,s\in\hZZ$.}
\end{gather}
We write $u(r)$ for $u\otimes t^r$ when $u\in\gt{u}$ and $r\in\hZZ$.
We consider the Clifford algebras $\Cl(\hat{\gt{u}})$ and
$\Cl(\hat{\gt{u}}_{\ogi})$. The inclusion of $\gt{u}$ in
$\hat{\gt{u}}_{\ogi}$ given by $u\mapsto u(0)$ induces an embedding
of algebras $\Cl(\gt{u})\hookrightarrow \Cl(\hat{\gt{u}}_{\ogi})$.
For $S=(i_1,\ldots,i_k)$ an ordered subset of $\fset$ we write
$e_S(r)$ for the element $e_{i_1}(r)\cdots e_{i_k}(r)$, which lies
in $\Cl(\hat{\gt{u}})$ or $\Cl(\hat{\gt{u}}_{\ogi})$ according as
$r$ is in $\ZZh$ or $\ZZ$. With this notation $e_S(0)$ coincides
with the image of $e_S$ under the embedding
$\Cl(\gt{u})\hookrightarrow \Cl(\hat{\gt{u}}_{\ogi})$. Let $E$ be an
$\FF_2^{\mc{E}}$-homogeneous subgroup of $\Cl(\gt{u})^{\times}$ such
that $\mc{C}(E)$ is a self-dual doubly-even code on $\fset$.

We write $\mc{B}(\hat{\gt{u}})$ for the subalgebra of
$\Cl(\hat{\gt{u}})$ generated by the $u(m+\tfrac{1}{2})$ for
$u\in\gt{u}$ and $m\in\ZZ_{\geq 0}$. We write
$\mc{B}(\hat{\gt{u}}_{\ogi})_E$ for the subalgebra of
$\Cl(\hat{\gt{u}}_{\ogi})$ generated by $E\subset\Cl(\gt{u})$, and
the $u(m)$ for $u\in\gt{u}$ and $m\in\ZZ_{> 0}$. Let $\RR_1$
denote a one-dimensional module for either $\mc{B}(\hat{\gt{u}})$
or $\mc{B}(\hat{\gt{u}}_{\ogi})_E$, spanned by a vector $1_E$,
such that $u(r)1_E=0$ whenever $r\in\hZZ_{>0}$, and such that
$g(0)1_{E}=1_{E}$ for $g\in E$. We write $A(\gt{u})$ (respectively
$A(\gt{u})_{\ogi,E}$) for the $\Cl(\hat{\gt{u}})$-module
(respectively $\Cl(\hat{\gt{u}}_{\ogi})$-module) induced from the
$\mc{B}(\hat{\gt{u}})$-module structure (respectively
$\mc{B}(\hat{\gt{u}}_{\ogi})_E$-module structure) on $\RR_{1}$.
\begin{gather}
    A(\gt{u})=\Cl(\hat{\gt{u}})
        \otimes_{\mc{B}(\hat{\gt{u}})}\RR_{1},\qquad
    A(\gt{u})_{\ogi,E}=\Cl(\hat{\gt{u}}_{\ogi})
        \otimes_{\mc{B}(\hat{\gt{u}}_{\ogi})_E}\RR_{1}.
\end{gather}
We write ${\bf 1}$ for $1\otimes 1_{E}\in A(\gt{u})$, and we write
${\bf 1}_{\ogi}$ for $1\otimes 1_{E}\in A(\gt{u})_{\ogi,E}$. When
no confusion will arise, we simply write $A(\gt{u})_{\ogi}$ in
place of $A(\gt{u})_{\ogi,E}$. %We often write $u_1(-r_1)\cdots
%u_k(-r_k)$ for $u_1(-r_1)\cdots u_k(-r_k){\bf 1}$ when the $r_i$
%are in $\ZZh_{>0}$, and similarly, we write $u_1(-m_1)\cdots
%u_k(-m_k)$ for $u_1(-m_1)\cdots u_k(-m_k){\bf 1}_{\ogi}$ when the
%$m_i$ are in $\ZZ_{\geq 0}$.

\medskip

The space $A(\gt{u})$ supports a structure of SVOA. In order to
define the vertex operators we require the notion of {\em fermionic
normal ordering} for elements in $\Cl(\hat{\gt{u}})$ and
$\Cl(\hat{\gt{u}}_{\ogi})$. The fermionic normal ordering on
$\Cl(\hat{\gt{u}})$ is the multi-linear operator defined so that for
$u_i\in\gt{u}$ and $r_i\in\ZZh$ we have
\begin{equation}
    :u_1(r_1)\cdots u_k(r_k):=
        {\rm sgn}(\sigma)u_{\sigma 1}(r_{\sigma 1})
        \cdots u_{\sigma k}(r_{\sigma k})
\end{equation}
where $\sigma$ is any permutation of the index set
$\{1,\ldots,k\}$ such that $r_{\sigma 1}\leq\cdots\leq r_{\sigma
k}$. For elements in $\Cl(\hat{\gt{u}}_{\ogi})$ the fermionc
normal ordering is defined in steps by first setting
\begin{equation}
    :u_1(0)\cdots u_k(0):=\frac{1}{k!}\sum_{\sigma\in S_k}
        {\rm sgn}(\sigma)u_{\sigma 1}(0)
        \cdots u_{\sigma k}(0)
\end{equation}
for $u_i\in \gt{u}$. Then in the situation that $n_i\in\ZZ$ are
such that $n_{i}\leq n_{i+1}$ for all $i$, and there are some $s$
and $t$ (with $1\leq s\leq t\leq k$) such that $n_j=0$ for $s\leq
j\leq t$, we set
\begin{equation}
    :u_1(n_1)\cdots u_k(n_k):=
    u_{1}(n_{1})\cdots u_{ s-1}(n_{s-1}):u_s(0)\cdots u_t(0):
    u_{t+1}(n_{t+1})\cdots u_{k}(n_k)
\end{equation}
Finally, for arbitrary $n_i\in\ZZ$ we set
\begin{equation}
    :u_1(n_1)\cdots u_k(n_k):=
        {\rm sgn}(\sigma):u_{\sigma 1}(n_{\sigma 1})
        \cdots u_{\sigma k}(n_{\sigma k}):
\end{equation}
where $\sigma$ is again any permutation of the index set
$\{1,\ldots,k\}$ such that $n_{\sigma 1}\leq\cdots\leq n_{\sigma
k}$, and we extend the definition multi-linearly to
$\Cl(\hat{\gt{u}}_{\ogi})$.

Now for $u\in\gt{u}$ we define the generating function, denoted
$u(z)$, of operators on $A(\gt{u})_{\Theta}=A(\gt{u})\oplus
A(\gt{u})_{\ogi}$ by $u(z)=\sum_{r\in\hZZ}u(r)z^{-r-1/2}$. Note
that $u(r)$ acts as $0$ on $A(\gt{u})$ if $r\in \ZZ$, and acts as
$0$ on $A(\gt{u})_{\ogi}$ if $r\in \ZZh$. To an element $a\in
A(\gt{u})$ of the form $a=u_{1}(-m_1-\tfrac{1}{2})\cdots
u_{k}(-m_k-\tfrac{1}{2}){\bf 1}$ for $u_i\in \gt{u}$ and $m_i\in
\ZZ_{\geq 0}$, we associate the operator valued power series
$\overline{Y}(a,z)$, given by
\begin{gather}
    \overline{Y}(a,z)=:D_z^{(m_1)}u_{i_1}(z)\cdots
        D_z^{(m_k)}u_{i_k}(z):
\end{gather}
We define the vertex operator correspondence
\begin{gather}
    Y(\cdot\,,z):A(\gt{u})\otimes A(\gt{u})_{\Theta}
        \to A(\gt{u})_{\Theta}((z^{1/2}))
\end{gather}
by setting $Y(a,z)b=\overline{Y}(a,z)b$ when $b\in A(\gt{u})$, and
by setting $Y(a,z)b=\overline{Y}(e^{\Delta_z}a,z)b$ when $b\in
A(\gt{u})_{\ogi}$, where $\Delta_z$ is the expression defined by
\begin{gather}
    \Delta_z=-\frac{1}{4}\sum_i\sum_{m,n\in\ZZ_{\geq 0}}C_{mn}
        e_i(m+\tfrac{1}{2})e_i(n+\tfrac{1}{2})z^{-m-n-1}\\
        C_{mn}=\frac{1}{2}\frac{(m-n)}{m+n+1}
            \binom{-\tfrac{1}{2}}{m}\binom{-\tfrac{1}{2}}{n}
\end{gather}
Set $\cas=-\tfrac{1}{4}\sum_i e_i(-\tfrac{3}{2})
e_i(-\tfrac{1}{2}){\bf 1} \in A(\gt{u})_2$. Then one has the
following
\begin{thm}[\cite{FFR}]
The map $Y$ defines a structure of self-dual rational SVOA on
$A(\gt{u})$ when restricted to $A(\gt{u})\otimes A(\gt{u})$. The
Virasoro element is given by $\cas$, and the rank is
$\tfrac{1}{2}\dim(\gt{u})$. The map $Y$ defines a structure of
$\ogi$-twisted $A(\gt{u})$-module on $A(\gt{u})_{\ogi}$ when
restricted to $A(\gt{u})\otimes A(\gt{u})_{\ogi}$.
\end{thm}
The definition of $Y(a,z)b$ for $b\in A(\gt{u})_{\ogi}$ appears
quite complicated, but all we need to know about these operators
is contained in the following
\begin{prop}\label{CliffTwOpsAWNTK}
Let $b\in A(\gt{u})_{\ogi}$.
\begin{enumerate}
\item   If $a={\bf 1}\in A(\gt{u})_0$ then $Y(a,z)b=b$.
\item   If $a\in A(\gt{u})_1$ then $\Delta_za=0$ so that
$Y(a,z)b=\overline{Y}(a,z)b$.
\item   If $a\in A(\gt{u})_2$ then $\Delta_za=0$ and $Y(a,z)b=
\overline{Y}(a,z)b$ unless $a$ has non-trivial projection onto the
span of the vectors $e_i(-\tfrac{3}{2})e_i(-\tfrac{1}{2})\vac$ for
$i\in \Omega$.
\item   For $a=e_i(-\tfrac{3}{2})e_i(-\tfrac{1}{2})\vac$ we have
$\Delta_za=-\tfrac{1}{4}z^{-2}$ and $\Delta_z^2a=0$ so that
$Y(a,z)b=\overline{Y}(a,z)b-\tfrac{1}{4}bz^{-2}$ in this case.
\end{enumerate}
\end{prop}
As a corollary of Proposition~\ref{CliffTwOpsAWNTK} we have that
$Y(\omega,z){\bf 1}_{\ogi}=\tfrac{1}{16}{\rm dim}(\gt{u}){\bf
1}_{\ogi} z^{-2}$, and consequently the $L(0)$ grading on
$A(\gt{u})_{\Theta}$ is given by
\begin{gather}
    A(\gt{u})=\coprod_{n\in\hZZ_{\geq 0}}A(\gt{u})_n,\quad
    A(\gt{u})_{\ogi}=\coprod_{n\in\hZZ_{\geq 0}}
        (A(\gt{u})_{\ogi})_{h+n},
\end{gather}
where $h=\tfrac{1}{16}{\rm dim}(\gt{u})$. Given a specific choice of
$E$, the embedding of $\Cl(\gt{u})$ in $\Cl(\hat{\gt{u}}_{\ogi})$
gives rise to an isomorphism of $\Cm(\gt{u})_E$ with
$(A(\gt{u})_{\ogi,E})_{h}$, and it will be convenient to consider
these spaces as identified. We may have occasion to replace $\gt{u}$
by its complexification $_{\CC}\gt{u}$ in the above; in such a
situation we shall understand $A(_{\CC}\gt{u})_{\Theta}$ to be the
complexified space $_{\CC}A(\gt{u})_{\Theta}$.

\section{Structure of $\afn$}\label{sec:strucafn}

\subsection{Construction}\label{CliffConst}

Suppose that $\Omega$ is some finite set with cardinality $24$ and
an arbitrary ordering. Let $\gt{l}$ be a $24$ dimensional vector
space over $\RR$ with positive definite bilinear form, and let
$\mc{E}=\{e_i\}_{i\in\Omega}$ be an orthonormal basis for
$\gt{l}$. The goal of this section is to show that the space
$\afn$ given by
\begin{gather}
    \afn=A(\gt{l})^0\oplus A(\gt{l})_{\ogi}^0
\end{gather}
admits a structure of self-dual rational $N=1$ SVOA.

Let $\mc{G}\subset\mc{P}(\Omega)$ be a copy of the Golay code. Let
$G$ be an $\FF_2^{\mc{E}}$-homogeneous subgroup of $\Sp(\gt{l})$
such that $G$ does not contain $-{\bf 1}\in\Sp(\gt{l})$ and the
associated code $\mc{C}(G)$ is $\mc{G}$. Curtesy of
\S\ref{sec:cliffalgs} we have the $\Cl(\gt{l})$-module
$\Cm(\gt{l})_G=\Cl(\gt{l})\otimes_{\RR G}\RR_1$, and this space can
be used to give an explicit realization of the $A(\gt{l})$-module
$A(\gt{l})_{\ogi}$. From now on we set
$A(\gt{l})_{\ogi}=A(\gt{l})_{\ogi,G}$. Notice that the odd parity
subspace of $\afn$ is precisely $A(\gt{l})_{\ogi}^0$. Since $\gt{l}$
has dimension $24$, the $L(0)$-homogeneous subspace of
$A(\gt{l})_{\ogi}$ with minimal degree is
$(A(\gt{l})_{\ogi})_{3/2}$, and this space is identified with the
$\Cl(\gt{l})$-module $\Cm(\gt{l})_G$. Note that ${\bf
1}_{\ogi}\leftrightarrow 1_{G}$ under this identification. Also, the
bilinear form on $\afn$ coincides with that on $\Cl(\gt{l})_G$ when
restricted to $(\afn)_{3/2}$, and in particular, is normalized so
that $\langle {\bf 1}_{\ogi}|{\bf 1}_{\ogi}\rangle=1$.

We require a vertex operator correspondence
$Y:\afn\otimes\afn\to\afn((z))$ and as yet this map is defined
only on $(\afn)_{\bar{0}}\otimes(\afn)_{\bar{0}}$ and on
$(\afn)_{\bar{0}}\otimes(\afn)_{\bar{1}}$. Such a map $Y$ must
satisfy skew-symmetry if it exists, so for $u\otimes v\in
(\afn)_{\bar{1}}\otimes(\afn)_{\bar{0}}$ we define $Y(u,z)v$ by
\begin{gather}\label{YdefAtwOnA}
  Y(u,z)v=e^{zL(-1)}Y(v,- z)u
\end{gather}
(since $|u||v|=0$ in this case). Suppose now that $u\otimes
v\in(\afn)_{\bar{1}}\otimes (\afn)_{\bar{1}}$. Motivated by
\S\ref{Sec:AdjOps} we define $Y(u,z)v$ by requiring that for any
$w\in (\afn)_{\bar{0}}$ we should have
\begin{gather}\label{AdjOfAtwOnAtw}
  \langle Y(u,z)v\mid w\rangle=
  (-1)^n\langle v\mid Y(e^{zL(1)}z^{-2L(0)}
  u,z^{-1})w\rangle
\end{gather}
whenever $u\in (\afn)_{n-1/2}$ for $n\in\ZZ$. Now the operator on
the right of (\ref{AdjOfAtwOnAtw}) is defined by (\ref{YdefAtwOnA}).
We can use this later expression to rewrite (\ref{AdjOfAtwOnAtw}) in
terms of the operator $Y$ defined on $(\afn)_{\bar{0}}\otimes\afn$
in \S\ref{sec:cliffalgs:SVOAs}, and doing so we obtain the following
convenient working definition for the operator $Y$ on
$(\afn)_{\bar{1}}\otimes (\afn)_{\bar{1}}$. For $u\in
(\afn)_{n-1/2}$ with $n\in\ZZ$, for $v\in(\afn)_{\bar{1}}$ and
$w\in(\afn)_{\bar{0}}$ we have
\begin{gather}\label{YdefAtwOnAtw}
    \langle Y(u,z)v\mid w\rangle
    =(-1)^n\langle e^{z^{-1}L(1)}v\mid
    Y(w,- z^{-1})e^{zL(1)} z^{-2L(0)}
    u\rangle
\end{gather}

\begin{prop}\label{AfnIsSVOA}
The map $Y:\afn\otimes\afn\to\afn((z))$ defines a structure of rank
$12$ self-dual rational SVOA on $\afn$.
\end{prop}
\begin{proof}
Let $_{\CC}\afn$ denote the complexification of $\afn$. Then
$(_{\CC}\afn)_{\bar{0}}$ is a simple VOA of rank $12$, and
$(_{\CC}\afn)_{\bar{1}}$ is an irreducible module over
$(_{\CC}\afn)_{\bar{0}}$. Let us write $Y_{\bar{k}\bar{l}}$ for the
restriction of $Y$ to $(_{\CC}\afn)_{\bar{k}}\otimes\,
(_{\CC}\afn)_{\bar{l}}$ for ${k},{l}\in\{{0},{1}\}$.

By the Boson-Fermion correspondence \cite{FreBF} (see also
\cite{DoMaBF}) we have that $(_{\CC}\afn)_{\bar{0}}$ is isomorphic
to a lattice VOA $_{\CC}V_{M_0}$ where $M_0$ is an even lattice of
type $D_{12}$. The irreducible modules over a lattice VOA
$_{\CC}V_L$ for $L$ an even lattice are known to be indexed by the
cosets of $L$ in its dual $L^*=\{u\in {\RR}\otimes_{\ZZ} L\mid
\langle u,L\rangle\subset\ZZ\}$ \cite{DonVAsLats}. In particular, a
lattice VOA is rational. Further, it is known that the fusion
algebra associated to the modules over $_{\CC}V_L$ coincides with
the group algebra of $L^*/L$ in the natural way. (One may refer to
\cite{DonLepGVAs} for a thorough treatment.) In the case that
$L=M_0$, there are exactly three non-trivial cosets, and for any one
of these $M_0+\mu$ say, the set $M_0\cup (M_0+\mu)$ forms an
integral lattice in $_{\RR}L=\RR\otimes_{\ZZ}L$, and in particular,
$M_0^*/M_0\cong\ZZ/2\times\ZZ/2$ has exponent two.

From this we conclude that $_{\CC}\afn$ is isomorphic to
$_{\CC}V_{M}$ for $M=M_0\cup (M_0+\mu)$ for some $\mu\in
M_0^*\setminus M_0$, that there exist unique up to scale
intertwiners of types $\binom{M_0+\mu}{M_0+\mu\, M_0}$ and
$\binom{M_0}{M_0+\mu\,M_0+\mu}$ for $_{\CC}V_{M_0}$, and that these
intertwiners are just those obtained by restricting the SVOA
structure on $_{\CC}V_{M}$. In particular, there is a unique
structure of rank $12$ SVOA on $_{\CC}\afn$. On the other hand it is
known from \cite{DonLiMasSmpCrt} for example, that the maps
$Y_{\bar{1}\bar{0}}$ and $Y_{\bar{1}\bar{1}}$ defined by equations
(\ref{YdefAtwOnA}) and (\ref{YdefAtwOnAtw}), respectively, yield
intertwining operators of types $\binom{\bar{1}}{\bar{1}\,\bar{0}}
=\binom{M_0+\mu}{M_0+\mu\, M_0}$ and
$\binom{\bar{0}}{\bar{1}\,\bar{1}} =\binom{M_0}{M_0+\mu\,M_0+\mu}$,
respectively, for $(_{\CC}\afn)_{\bar{0}}$. By uniqueness, they must
coincide with those inherited from the SVOA structure on
$_{\CC}V_{M}$ up to some scalar factors, and in any case the map $Y$
defined above furnishes $_{\CC}\afn$ with a structure of rational
rank $12$ SVOA. We have chosen scalars in such a way that $\afn$ is
a real form for $_{\CC}\afn$. One can check directly that $M$ is a
self-dual lattice (given that $(\afn)_{1/2}$ is trivial, $M$ must be
a copy of the $D_{12}^+$ lattice --- the unique self-dual lattice of
rank $12$ with no vectors of unit norm \cite[Ch.19]{CoS93}), and it
follows from the above that $_{\CC}V_{M}$ and $_{\CC}\afn$ are then
self-dual SVOAs. This completes the proof of the proposition.
\end{proof}
\begin{rmk}
The method used here to extend the vertex operator map from a VOA
to the sum of the VOA and a module over it was given earlier in
\cite{HuaXtnMoonVOA}.
%The proof of Proposition~\ref{AfnIsSVOA} shows that the SVOA
%structure on $\afn$ is essentially unique.
\end{rmk}
The following proposition gives a convenient criterion for when a
vector in $(\afn)_{3/2}$ is superconformal. In section
\S\ref{SecUniq} it will be shown that all such vectors are
equivalent up to the action of $\Sp(\gt{l})$.
\begin{prop}\label{Prop:CodeCriForSC}
Suppose that $t\in(\afn)_{3/2}$ is such that $\langle
t|t\rangle=8$ and $\langle e_C(0)t|t\rangle=0$ whenever
$C\subset\Omega$ has cardinality two or four. Then $t$ is a
superconformal vector for $\afn$.
\end{prop}
\begin{proof}
We should compute $t_{n}t$ for $n=0$, $n=1$ and $n=2$, and then
compare with the results of Proposition~\ref{prop:SConfCriterion}.
Using (\ref{YdefAtwOnAtw}) and recalling $t\in (\afn)_{2-1/2}$ we
find that for arbitrary $u\in\afn$ we have
\begin{gather}
    \begin{split}
    \langle u| Y(t,z)t\rangle
        =&\langle
            Y(u,-z^{-1})
            e^{zL(1)}z^{-2L(0)}t
            |e^{z^{-1}L(1)}t\rangle\\
        =&\langle Y(u,-z^{-1})t|
        t\rangle z^{-3}
    \end{split}
\end{gather}
Then for $t_nt$ we obtain
\begin{gather}\label{tauAontauA}
    \langle u|t_nt\rangle
        ={\rm Res}_{z=0}\langle
            Y(u,-z^{-1}) t|t\rangle z^{n-3}
        = \langle u_{-n+1}t|t\rangle
            (-1)^{n-2}
\end{gather}
For $L(0)$-homogeneous $u\in \afn$ the expression (\ref{tauAontauA})
is zero unless $u\in (\afn)_{-n+2}$. In order to determine $t_nt$,
we should compute $u_mt$ for various $u\in (\afn)^0$, and apply the
equation (\ref{tauAontauA}). To compute $u_mt$ we will use the
results of Proposition~ \ref{CliffTwOpsAWNTK}.

For $n=2$ the expression (\ref{tauAontauA}) is zero unless
$u\in(\afn)_0=\RR{\bf 1}$, in which case we obtain $\langle {\bf
1}| t_2t\rangle = \langle{\bf 1}_{-1}t|t \rangle$. Since $Y({\bf
1},z)={\bf 1}$ and $\langle {\bf 1}|{\bf 1}\rangle=1$, we find
that $t_2t=\langle t|t\rangle{\bf 1}$.

For the case that $n=1$ we should consider the $u$ in $(\afn)_1$.
Suppose $u=e_i(-\tfrac{1}{2})e_j(-\tfrac{1}{2})\vac$ for $i\neq j\in
\Omega$. Then
\begin{gather}
    Y(u,z)t =\overline{Y}(u,z)t
    =e_{ij}(0)tz^{-1}+\ldots
\end{gather}
By hypothesis we have that $e_{ij}(0)t$ is orthogonal to $t$, so
$t_1t=0$.

Finally, when $n=0$ we are concerned with $u_{1}t$ for $u\in
(\afn)_2$. The space $(\afn)_2$ is spanned by the vectors of the
form $e_i(-\tfrac{3}{2})e_j(-\tfrac{1}{2})\vac$ for $i,j\in \Omega$,
and also by the $e_C(-\tfrac{1}{2})\vac$ for
$C=\{i_1,i_2,i_3,i_4\}\subset\Omega$. For $u$ one of these vectors
we have $Y(u,z)t=\overline{Y}(u,z) t=e_{C}(0)t z^{-2}+\ldots$ for
$C\subset \Omega$ of size two or four unless
$u=e_{i}(-\tfrac{3}{2})e_i(-\tfrac{1}{2})\vac$. In the former cases,
$u_1t$ is orthogonal to $t$ by hypothesis, and in the later case we
have
\begin{gather}
    Y(u,z)t=
    \overline{Y}(u,z)t
        -\tfrac{1}{4}tz^{-2}%+\ldots
        =-\tfrac{1}{4}tz^{-2}+\ldots
\end{gather}
The expression (\ref{tauAontauA}) now reduces to $\langle
u|t_0t\rangle=-\tfrac{1}{4}\langle t|t\rangle=-2$, and we conclude
that $t_0t=-\tfrac{1}{2} \sum_{\Omega} e_{i}(-\tfrac{3}{2})
e_i(-\tfrac{1}{2})\vac$ since we have $\langle u|u\rangle=4$ for
$u=e_i(-\tfrac{3}{2})e_i(-\tfrac{1}{2})\vac$ for any $i$ in
$\Omega$.

We have verified that $t_2t=8$, $t_1t=0$ and $t_0t=2\cas$. Since the
rank of $\afn$ is $12$ and $8=\tfrac{2}{3}12$, an application of
Proposition~\ref{prop:SConfCriterion} confirms that $t$ is
superconformal for $\afn$.
\end{proof}

Set $\tau_A=\sqrt{8}{\bf 1}_{\ogi}\in (\afn)_{3/2}$. Then $\langle
\tau_A|\tau_A\rangle=8$, and we have
\begin{cor}\label{tauAIsSC}
The vector $\tau_A$ is a superconformal vector for $\afn$.
\end{cor}
\begin{proof}
That $\tau_A$ satisfies the hypotheses of Proposition~
\ref{Prop:CodeCriForSC} is a consequence of the fact that the Golay
code $\mc{G}=\mc{C}(G)$ has minimum weight $8$.
\end{proof}

We record the results of Proposition~\ref{AfnIsSVOA} and Corollary
\ref{tauAIsSC} in the following
\begin{thm}\label{ThmCnstafn}
The quadruple $(\afn,Y,\vac,\tau_A)$ is a self-dual rational $N=1$
SVOA.
\end{thm}

\subsection{Symmetries}

In this section we show that the automorphism group of the $N=1$
SVOA structure on $\afn$ is isomorphic to Conway's largest sporadic
group, $\Co_1$.

The operators $x_{0}$ for $x\in A(\gt{l})_1$ span a Lie algebra of
type $D_{12}$ in $\End( A(\gt{l})_{\Theta})$, and the exponentials
$\exp(x_{0})$ for $x\in A(\gt{l})_1$ generate a group $S$ say, which
acts as $A(\gt{l})^0$-module automorphisms of
$A(\gt{l})_{\Theta}=A(\gt{l})\oplus A(\gt{l})_{\ogi}$. This group
$S$ is isomorphic to the group $\Sp(\gt{l})$, and we may choose the
isomorphism so that $\exp(x_{0})$ in $S$ corresponds to
$\exp(\tfrac{1}{2}(ab-ba))$ in $\Sp(\gt{l})<\Cl(\gt{l})^{\times}$
when $x=a(-\tfrac{1}{2})b(-\tfrac{1}{2})\vac\in A(\gt{l})_1$ for
some $a,b\in\gt{l}$. The action of $S$ on $A(\gt{l})_{\Theta}$
commutes with the action of $\ogi$, and so preserves the subspace
$\afn=A(\gt{l})^0\oplus A(\gt{l})_{\ogi}^0$. The kernel of this
action is the group of order $2$ generated by $\ogi$. Taking the
complexification $_{\CC}\gt{l}$ in place of $\gt{l}$ in the above we
obtain an action of the complex Lie group $\Sp(_{\CC}\gt{l})$ on the
complexified SVOA $_{\CC}\afn=A(_{\CC}\gt{l})^0\oplus
A(_{\CC}\gt{l})^0_{\ogi}$. Let us write $_{\CC}S$ for this copy of
$\Sp(_{\CC}\gt{l})$ generated by exponentials $\exp(x_0)$ with $x$
in $A(_{\CC}\gt{l})_1$.

\begin{prop}
The group $_{\CC}S$ maps surjectively onto the group of SVOA
automorphisms of $_{\CC}\afn$.
\end{prop}
\begin{proof}
From the proof of Proposition~\ref{AfnIsSVOA} we may regard
$_{\CC}\afn$ as the (complex) lattice SVOA $_{\CC}V_{M}$ where $M$
is a lattice of type $D_{12}^+$. Let us write $M^0$ for the even
sublattice of $M$ (of type $D_{12}$), and $M^1$ for the unique coset
of $M^0$ in $M$; we may write $_{\CC}V_{M^0}\oplus\,_{\CC}V_{M^1}$
for the superspace decomposition of $_{\CC}V_M$. Let us write ${\sf
G}$ (not to be confused with the $G$ of \S\ref{CliffConst}) for the
group of SVOA automorphisms of $_{\CC}V_M$, and ${\sf G}^0$ for the
group of VOA automorphisms of $_{\CC}V_{M^0}$. Let ${\sf S}$ denote
the image of $_{\CC}S$ in ${\sf G}=\Aut_{{\rm SVOA}}(_{\CC}\afn)$.
We wish to show that ${\sf S}={\sf G}$.

Any element of ${\sf G}$ preserves the superspace structure on
$_{\CC}V_M$, so we have a natural map $\phi:{\sf G}\to {\sf G}^0$.
By a similar token, any element of ${\sf G}^0$ preserves the Lie
algebra structure on the degree $1$ subspace of $_{\CC}V_{M^0}$ (we
denote this Lie algebra by $\gt{g}$) so we have also a natural map
${\sf G}^0\to\Aut(\gt{g})$. In fact, this map is faithful and onto
since $_{\CC}V_{M^0}$ is generated by its subspace of degree $1$
elements, in the sense that we have
\begin{gather}
     _{\CC}V_{M^0}={\rm Span}_{\CC}\left\{
          x^1_{-n_1}x^2_{-n_2}\cdots x^r_{-n_r}\vac
          \mid \deg(x^i)=1,\; n_i\in\ZZ_{>0}\right\}
\end{gather}
so that any element $g\in\Aut(\gt{g})$ extends to an element of
$\Aut(_{\CC}V_{M^0})={\sf G}^0$ once we decree
\begin{gather}
     g:x^1_{-n_1}x^2_{-n_2}\cdots x^r_{-n_r}\vac\mapsto
          (gx^1)_{-n_1}(gx^2)_{-n_2}\cdots (gx^r)_{-n_r}\vac
\end{gather}
and any element of ${\sf G}^0$ that fixes $\gt{g}$ fixes all of
$_{\CC}V_{M^0}$. Thus we may identify $\Aut(\gt{g})$ with ${\sf
G}^0$.

We claim that $\phi({\sf S})=\phi({\sf G})$. In
\cite{DonNagAutsLattVOAs} it is proved that the automorphism group
of a lattice VOA (for an even positive definite lattice, such as
$M^0$) is generated by exponentials of zero-modes of degree $1$
elements and by lifts of automorphisms of the lattice. In our
situation this means ${\sf G}^0=\lab\phi({\sf S}),O(M^0)\rab$ where
$O(M^0)$ denotes the group of automorphisms of $_{\CC}V_{M^0}$
generated by lifts of elements of $\Aut(M^0)$. On the other hand, we
know that the group $\Inn(\gt{g})$ of inner automorphisms of
$\gt{g}$ (this is just our group $\phi({\sf S})$) has index $2$ in
$\Aut(\gt{g})$, since $\gt{g}$ is a simple complex Lie algebra of
type $D_{12}$. So there is some $x\in O(M^0)$ of order $2$, such
that ${\sf G}^0=\phi({\sf S})\cup x\phi({\sf S})$, and either
$\phi({\sf G})=\phi({\sf S})$, or $\phi({\sf G})={\sf G}^0$. Let
$\bar{x}$ denote the canonical image of $x$ in $\Aut(M^0)$. The
coset $x\phi({\sf S})$ corresponds to a so-called diagram
automorphism of $D_{12}$, and $\bar{x}$ acts non-trivially on the
coset space $(M^0)^*/M^0$ interchanging the two cosets with minimal
norm $3$ (one of which is $M^1$). In particular, $\bar{x}$ does not
preserve the lattice $M=M^0+M^1$, and thus $x$ cannot be extended to
an automorphism of $_{\CC}V_M$ (c.f. \cite[Lemma~
2.3]{DonNagAutsLattVOAs}). We conclude that $\phi({\sf S})=\phi({\sf
G})$.

Next we claim that $\ker(\phi)$ is contained in ${\sf S}$. For
suppose $g\in\ker(\phi)$. Then $g$ fixes all elements of $\gt{g}$,
and therefore commutes with the action of ${\sf S}$ on
$(_{\CC}V_M)_{3/2}=\Cm(_{\CC}\gt{l})_G^0$. The space
$\Cm(_{\CC}\gt{l})^0_G$ is irreducible for the action of ${\sf S}$,
so $g$ acts as scalar multiplication by $\zeta\in\CC$ say, on
$(_{\CC}V_M)_{3/2}$, and indeed, on all of $_{\CC}V_{M^1}$ (since
$_{\CC}V_{M^1}$ is generated by the action of $_{\CC}V_{M^0}$ on
$(_{\CC}V_M)_{3/2}$). Then for $x,y\in\, _{\CC}V_{M^1}$, we have
$g(x_{(n)}y)=(gx)_{(n)}(gy)=\zeta^2x_{(n)}y$, and also
$g(x_{(n)}y)=x_{(n)}y$ since $x_{(n)}y$ lies in $_{\CC}V_{M^0}$.  It
follows that $\zeta=\pm 1$ and $\ker(\phi)$ has order $2$. In the
non-trivial case that $g|_{_{\CC}V_{M^0}}=\Id$ and
$g|_{_{\CC}V_{M^1}}=-\Id$, we have that $g$ is realized by the image
of the element $-{\bf 1}\in S$ in ${\sf S}$. This proves the claim.

We have shown that $\phi({\sf S})=\phi({\sf G})$ and
$\ker(\phi)<{\sf S}$, and it follows that ${\sf S}={\sf G}$, which
is what we required.
\end{proof}
From now on it will be convenient to regard $\Sp(\gt{l})$ and
$\Sp(_{\CC}\gt{l})$ as groups of SVOA automorphisms of $\afn$ and
$_{\CC}\afn$, respectively.

Recall that the ordering on $\Omega$ is chosen so that
$e_{\Omega}\in\Sp(\gt{l})$ lies in $G$. We denote $e_{\Omega}$ also
by $\gt{z}$. Then the action of $e_{\Omega}$ on $_{\CC}\afn$ is
trivial, the kernel of the map $\Sp(_{\CC}\gt{l})\to\Aut_{{\rm
SVOA}}(_{\CC}\afn)$ is $\{1,\ogz\}$, and the full group of SVOA
automorphisms of $_{\CC}\afn$ is $\Sp(_{\CC}\gt{l})/\lab\ogz\rab$.
Any automorphism of $\afn$ extends to an automorphism of the
complexification $_{\CC}\afn$, so $\Sp(_{\CC}\gt{l})/\lab\ogz\rab$
contains the group of SVOA automorphisms of the real form $\afn$. On
the other hand, this latter group contains
$\Sp(\gt{l})/\lab\ogz\rab$ which is maximal compact in
$\Sp(_{\CC}\gt{l})/\lab\ogz\rab$. We conclude that
$\Sp(\gt{l})/\lab\ogz\rab$ is the full group of SVOA automorphisms
of the real SVOA $\afn$.

Let $F$ denote the subgroup of $\Sp(_{\CC}\gt{l})$ that fixes
$1_G\in\Cm(_{\CC}\gt{l})^0_G\leftrightarrow (_{\CC}\afn)_{3/2}$.
Then the full group of automorphisms of the $N=1$ SVOA structure on
$_{\CC}\afn$ is $F/\lab\ogz\rab$. Let us write $\gt{l}'$ for the
span of the vectors $u1_G\in\Cm(\gt{l})_G$ for $u\in \gt{l}$, and
$_{\CC}\gt{l}'$ for the complexification of $\gt{l}'$, regarded as a
subspace of $\Cm(_{\CC}\gt{l})_G$. Then $_{\CC}\gt{l}'$ has
dimension $24$, and $F$ embeds naturally in $SO(_{\CC}\gt{l}')$
since $xu1_G=xux^{-1}1_G=x(u)1_G$ for $u\in \,_{\CC}\gt{l}$ and
$x\in F$, and $x(\gt{u})\subset\gt{u}$ for $x\in\Sp(_{\CC}\gt{l})$.
We will now show
\begin{prop}\label{FHasC}
The group $F$ contains a group isomorphic to $\Co_0$.
\end{prop}
\begin{proof}
Recall that the natural map $\Sp(_{\CC}\gt{l})\to\SO(_{\CC}\gt{l})$
is denoted $x\mapsto x(\cdot)$. Recall that $G$ is an
$\FF_2^{\mc{E}}$-homogeneous lift of the Golay code $\mc{G}$ to
$\Sp(\gt{l})$ (see \S\ref{sec:cliffalgs:mods}), so that for each
$C\in\mc{G}$ there is a unique $g_C\in G$ such that $g_C(\cdot)$ is
$-1$ on $e_i$ for $i\in C$, and $+1$ on $e_i$ otherwise. Let $C_0$
be a subgroup of $\SO(\gt{l})$ isomorphic to $\Co_0$ such that $C_0$
contains $g(\cdot)$ for each $g\in G<\Sp(\gt{l})$. The Golay code
construction of the Leech lattice given in \cite{ConLctXcptGps} for
example, shows that this is possible.

Let $\hat{C}$ be the preimage of $C_0$ in $\Sp(\gt{l})$. The group
$\Co_0$ has trivial Schur multiplier \cite{ATLAS} so there exists a
group $C'$ in $\Sp(\gt{l})$, a subgroup of index $2$ in $\hat{C}$,
such that the map $x\mapsto x(\cdot)$ restricts to an isomorphism of
$C'$ with $C_0<\SO(\gt{l})$. Set $\hat{G}=\{g_C,-g_C\mid C\in
\mc{G}\}$, and set $G'=\hat{G}\cap C'$. Then we have
$G'=\{\gamma_Cg_C\mid C\in G\}$ where $C\mapsto\gamma_C$ is a map
$\mc{G}\to\{\pm 1\}$ such that $\gamma_C\gamma_D=\gamma_{C+D}$. In
particular, $C\to \gamma_C$ is a homomorphism, and there must be
some $S\subset\Omega$ such that we have $\gamma_C=(-1)^{\lab
C,S\rab}$ for all $C\in\mc{G}$. In other words, we have
$G'=\{e_Sg{e_S}^{-1}\mid g\in G\}$.

Set $C=\{{e_S}^{-1}xe_S\mid x\in C'\}$. Then $C$ is isomorphic to
$\Co_0$, and contains $G$. In particular, $C$ contains the central
element $\ogz$, and $C/\lab\ogz\rab$ must be isomorphic to $\Co_1$.
The space $\Cm(\gt{l})_G^0$ is then a $C/\lab\gt{z}\rab$-module of
dimension $2048$ and since the only $\Co_1$ irreducibles with
dimension less than $2048$ have dimension $1$, $276$, $299$ and
$1771$ \cite{ATLAS}, the space $\Cm(\gt{l})_G^0$ must have a fixed
point $t$ say, for the action of $C$. We may assume that $t$ has
unit norm. Since $G$ is contained in $C$, the vector $t$ is also
invariant for $G$, and this forces $t=1_G$ by Proposition~3.3. We
conclude that $C$ is a subgroup of $F$ isomorphic to $\Co_0$, and
this completes the proof.
\end{proof}
The proof of Proposition~\ref{FHasC} shows that a copy of $\Co_0$
may be found even inside the intersection $F\cap\Sp(\gt{l})$.
\begin{prop}
The group $F$ is finite.
\end{prop}
\begin{proof}
$F$ is a subgroup of the algebraic group
$\Sp(_{\CC}\gt{l})\cong\Sp_{24}(\CC)$. The condition that a subspace
be stabilized by a linear transformation is polynomial, so $F$ too
is algebraic, since it is by definition the stabilizer of a subspace
in a representation of $\Sp(_{\CC}\gt{l})$. At the same time, $F$ is
a subgroup of $\SO(_{\CC}\gt{l}')\cong \SO_{24}(\CC)$ containing the
algebraic group $\Co_0$ by Proposition~\ref{FHasC}. Since the latter
group acts irreducibly on $_{\CC}\gt{l}'$, we conclude that $F$ is
reductive. We now check if there is any non-trivial semisimple
complex algebraic group or algebraic torus that can occur as a
factor of the connected component of the identity in $F$. Any such
group would have a non-trivial Lie algebra $\gt{k}$ say, with an
embedding in the Lie algebra $\gt{g}$ of $\Sp(_{\CC}\gt{l})$, which
we may identify with the degree one subspace of $_{\CC}\afn$
(equipped with the bracket $[x,y]=x_0y$). Now for all $x\in\gt{k}$
we have $\exp(x_0)1_G=1_G$, and this implies $x_01_G=0$ for some
non-trivial $x\in\gt{k}$. We claim that if $x\in\gt{g}$ satisfies
$x_01_G=0$ then $x=0$. For consider the map
$\gt{g}\to\Cm(_{\CC}\gt{l})_G$ given by $x\mapsto x_01_G$, and write
$\gt{g}'$ for the image of $\gt{g}$ under this map. Then the
dimension of $\gt{g}'$ is at most $276$. On the other hand $\gt{g}'$
contains the span of the vectors $\{e_{ij}1_G\}$ for $i<j$, and
since the Golay code has minimum weight $8$ these vectors are
linearly independent, and we see that $\gt{g}'$ has dimension not
less than $276$. It follows that the map $x\mapsto x_01_G$ is a
linear isomorphism from $\gt{g}$ to $\gt{g}'$, and in particular,
the kernel is trivial. This verifies the claim. We conclude that
$\dim(F)=0$, whence $F$ is finite.
\end{proof}
We have shown that $F$ is a finite subgroup of $\Sp(_{\CC}\gt{l})$
such that $F\cap\Sp(\gt{l})$ contains a copy of $\Co_0$. Our last
main task for this section is to show that $F$ itself is isomorphic
to $\Co_0$. With this in mind, we offer the following proposition,
the proof of which owes much to the methods used in Theorems 5.6 and
6.5 of \cite{NebRaiSloInvtsCliffGps}. In particular, we utilize the
notion of primitive matrix group: a group $G\leq\GL(V)$, for $V$ a
vector space, is said to be {\em primitive} if there is no
non-trivial decomposition $V=V_1\oplus\cdots\oplus V_k$ into
subspaces permuted by the the action of $G$. Note that if $N$ is
normal in $G$, then $G$ permutes the isotypic components of the
restricted module $V|_N$, so that $V|_N$ must be multiple copies of
a single irreducible representation for $N$ in the case that $G$ is
primitive.
\begin{prop}\label{CMaxSubjFin} The group $\Co_0$ is a
maximal subgroup of $\SO_{24}(\CC)$ subject to being finite.
\end{prop}
\begin{proof}
Any compact subgroup of $\SO_n(\CC)$ is realizable over $\RR$ (c.f.
\cite[\S13.2]{SeRep}), so it suffices to show that $\Co_0$ is a
maximal finite subgroup of $\SO_{24}(\RR)$. Let $V$ denote a real
vector space of dimension $24$, equipped with a non-degenerate
symmetric bilinear form. Suppose that $F$ is a finite subgroup of
$\SO(V)$ properly containing a copy $C$ of the group $\Co_0$. Then
$C$ is not normal in $F$, since $C$ acts absolutely irreducibly on
$V$, and if $C$ were normal in $F$ then $F/C$ would embed in the
outer automorphism group of $C$, which is trivial \cite{ATLAS}. Let
us write $Z$ for the center $\{\pm\Id\}$ of $C$ (and $F$).

We claim that $F$ has no non-trivial normal $p$-subgroups for $p$
odd, and the only non-trivial normal $2$-subgroup is $Z$. For
suppose $N$ is a normal $p$-subgroup of $F$ for some prime $p$. Then
$C_F(N)\cap C$ (we write $C_F(N)$ for the centralizer in $F$ of $N$)
is normal in $C$ and contains $Z$, so that $C_F(N)\cap C$ is either
$Z$ or $C$. In the former case we have that $C/Z\cong \Co_1$ is a
subgroup of $\Aut(N)$. In the latter case, $N$ is centralized by
(the absolutely irreducible action of) $C$, and hence must consist
of scalar matrices. It follows that $N$ is trivial unless $p=2$, in
which case $N$ is either trivial or $N=Z$. We suppose then that $N$
is a normal $p$-subgroup of $F$ such that $\Aut(N)$ contains a copy
of $\Co_1$. The group $C$ acts primitively on $V$, and even on the
complexification $_{\CC}V=\CC\otimes_{\RR}V$, and hence so does $F$.
It follows that $_{\CC}V|_N$ is an isotypic module for $N$; i.e.
several copies of a single irreducible module $M$ say, for $N$.
Since $N$ is by definition a subgroup of $\SO_{24}(\RR)$, it follows
that $M$ is the complexification of an irreducible $N$-module
$_{\RR}M$ say, defined over $\RR$, and that the action of $N$ on
$_{\RR}M$ is faithful. Any $p$-group has non-trivial center, and a
central subgroup of a group acts by scalar multiplications on any
irreducible module for that group. We conclude that $p$ is not odd,
since a $p$-group for odd $p$ has central elements which must act as
multiplication by primitive (and non-real) $p$-th roots of unity. We
see also that any abelian normal subgroup of $F$ is cyclic. The
irreducible representations of a $p$-group are of $p$-power order,
so $\deg(_{\RR}M)$ is a power of $2$ dividing $24$. Without loss of
generality, we suppose $\deg(_{\RR}M)=8$. Note that $N$ can contain
no noncyclic characteristic abelian subgroups, since such a subgroup
would be a noncyclic normal abelian $p$-subgroup of $F$. A $p$-group
with this property is said to be of {\em symplectic type}, and such
groups are classified by a Theorem of P. Hall (c.f.
\cite[(23.9)]{AscFGT}). In particular, there are no $2$-groups of
symplectic type that both embed in $\SO_8(\RR)$ and admit a
non-trivial action by $\Co_1$ as automorphisms. The claim follows.

Now we seek a contradiction. Any finite group is realizable over a
cyclotomic number field (c.f. \cite[Thm 24]{SeRep}). In fact we may
assume that $F$ is a subgroup of $SO_{24}(K)$ where $K$ is a totally
real abelian number field (e.g. $K=\QQ(\zeta_m+\zeta_m^{-1})$ for
$m$ such that $x^m=1$ for all $x\in F$ and $\zeta_m=\exp(2\pi\ii/m)$
--- c.f. \cite[Prop 5.6]{DreIndctnStrucThmsOrth}). Let $K$ be a
minimal such field, and let $R$ be the ring of integers in $K$. Then
$F$ preserves an $RC$-lattice, and such a lattice is of the form
$I\otimes_{\ZZ}\LL$ for $\LL$ a copy of the Leech lattice and $I$ a
fractional ideal of $R$, since any $C$ invariant lattice in
$\QQ^{24}$ is isometric to $\LL$ (c.f. \cite{TieGIRsFGs}). It
follows that $F$ preserves the lattice $R\otimes_{\ZZ}\LL$. (To see
this, note that if $F$ preserves a lattice $L$, then it also
preserves $aL$ for any $a\in K$. Also, if $F$ preserves lattices
$L_1$ and $L_2$, then it preserves the sum $L_1+L_2$. Now take
$L=I\otimes_{\ZZ}\LL$ for $I$ a fractional ideal of $R$, and let
$\{y_i\}$ be a set of generators for the inverse fractional ideal.
Then $R=\sum y_iI$, and $\sum y_iL=R\otimes_{\ZZ}\LL$ is also
invariant for $F$.) We may now regard $F$ as a group of matrices
with entries in the ring $R$. If $K=\QQ$ then $F=C$ and we are done.
If not, then there is some rational prime $p$ that ramifies in $K$.
Let $\gt{p}$ be a prime ideal of $R$ lying above $p$, and let
$\Gamma_{\gt{p}}$ denote the subgroup of $\Gal(K/\QQ)$ consisting of
automorphisms $\sigma$ such that $\sigma(a)\equiv a \pmod{\gt{p}}$
for all $a\in R$. (This is the {\em first inertia group}. It
stabilizes $\gt{p}$, and has order equal to the ramification index
of $\gt{p}$ over $p$
--- c.f. \cite[III:4]{FroTayANT}) The Galois group of $K$ over $\QQ$
acts on $F$ by acting componentwise on matrices. Let
$F_{\gt{p}}=\{g\in F\mid g\equiv\Id\pmod{\gt{p}}\}$, let $\sigma$ be
a non-trivial element of $\Gamma_{\gt{p}}$, and let $\phi:F\to
F_{\gt{p}}$ be the map defined by $\phi(g)=g^{-1}\sigma(g)$ for
$g\in F$. The group $F_{\gt{p}}$ is a normal $p$-subgroup of $F$,
and we have shown that such a group is trivial except possibly in
the case that $p=2$. If $F_{\gt{p}}$ is trivial then $\sigma$ fixes
$F$, and this contradicts the minimality of $K$. So suppose $p=2$
and $F_{\gt{p}}$ is the group $Z=\{\pm\Id\}$. Then $\phi$ is in fact
a group homomorphism $F\to Z$ (since $g^{-1}\sigma(g)$ is now
central for all $g\in F$). Since the image of $\phi$ is abelian, the
derived subgroup $F^{(1)}$ of $F$ lies in the kernel of $\phi$, and
is thus fixed by $\sigma$. If $F=F^{(1)}$ then $F$ is realizable
over the subfield of $K$ fixed by $\sigma$, contradicting the
minimality of $K$. So $F$ properly contains $F^{(1)}$, and the
argument thus far shows that any finite subgroup of $\SO_{24}(\RR)$
properly containing $C$, properly contains its own derived subgroup.
Consider now the descending chain $F\geq F^{(1)}\geq
F^{(2)}\geq\cdots$ where $F^{(k+1)}$ is the derived subgroup of
$F^{(k)}$. Each term contains $C$ since $F>C$ and $C=C^{(1)}$, and
thus each containment $F^{(k)}\geq F^{(k+1)}$ is proper unless
$F^{(k)}=C$. Since $F$ is finite, not all containments are proper,
and thus we have $F^{(k)}=C$ for some $k$. Then $C$ is a
characteristic subgroup of $F$, and in particular, is normal in $F$,
and this is again a contradiction.

We conclude that $\Co_0$ is a maximal subgroup of $\SO_{24}(\CC)$
subject to being finite.
\end{proof}

We have established the following
\begin{thm}\label{Thm:PtStabIsCo0} The subgroup
of $\Sp(_{\CC}\gt{l})$ fixing $1_G\in \Cm(_{\CC}\gt{l})_{G}$ is
isomorphic to $\Co_0$.
\end{thm}
Recall that $\Aut(_{\CC}\afn) =F/\lab \gt{z}\rab$. The group
$\langle\gt{z}\rangle$ is the center of $F$, and thus
$F/\lab\gt{z}\rab$ is isomorphic to $\Co_1$. By construction this
copy of $\Co_1$ is contained in $\Aut(\afn)$. We have therefore
established
\begin{thm}\label{ThmSymms}
There are isomorphisms of groups
$\Aut(_{\CC}\afn)\cong\Aut(\afn)\cong\Co_1$.
\end{thm}

\begin{rmk}
It was noted in the Introduction that an action of $\Co_1$ on the
SVOA underlying $\afn$ was considered earlier in \cite{BorRybMMIII}.
In fact, an action of $\Co_0$ (the perfect double cover of $\Co_1$)
on the SVOA underlying $\afn$ was also considered in
\cite{BorRybMMIII}, and in our setting, this action arises naturally
by considering the action of the group $F$ on the object
$A^{f\flat}$ given by
\begin{gather}
    A^{f\flat}=A(\gt{l})^0\oplus A(\gt{l})_{\ogi}^1
\end{gather}
where we realize $A(\gt{l})_{\ogi}$ as $A(\gt{l})_{G,\ogi}$ with $G$
as in \S\ref{CliffConst}. We have seen that the group $F$ is
isomorphic to the quasi-simple group $\Co_0$, and in contrast to the
situation with $\afn$ the central element of $F$ acts non-trivially
on $A^{f\flat}$. The same method used in \S\ref{CliffConst} shows
that $A^{f\flat}$ has a unique structure of SVOA, and also that
$\afn$ and $A^{f\flat}$ are isomorphic, as SVOAs. There is however
no $\Co_0$ invariant vector in the degree $3/2$ subspace of
$A^{f\flat}$, and hence no $\Co_0$ invariant $N=1$ structure on
$A^{f\flat}$.
\end{rmk}

\section{Uniqueness}\label{SecUniq}

In this section we prove a uniqueness result for $\afn$. In the
first subsection we verify that any nice rational $N=1$ SVOA
satisfying
\begin{itemize}
\item   self-dual
\item   rank $12$
\item   no small elements
\end{itemize}
is isomorphic to $_{\CC}\afn$ as an SVOA. To do this we first
recall the modularity results for trace functions on VOAs due to
Zhu (see \cite{ZhuPhd}, \cite{ZhuModInv}), and their extension to
the SVOA case given in \cite{HohnPhD}. Then we make use of some
techniques from \cite{DonMasEfctCC} and \cite{DonMasHlmVOA},
replacing VOA concepts with their SVOA analogues as necessary. The
guiding principle that we adopt from these two papers is that one
may use modular invariance results for a VOA $V$ to deduce
properties about the Lie algebra structure on $V_1$, the degree
one subspace of $V$.

In the second subsection we show that the $N=1$ structure on
$\afn$ is unique in the sense that if $\tau\in(\afn)_{3/2}$ is a
superconformal vector then there is some SVOA automorphism of
$\afn$ mapping $\tau$ to $\tau_A$.

\subsection{SVOA structure}

Recall from \S\ref{sec:SVOAstruc:SVOAs} and \S\ref{sec:SVOAMods} the
definitions of niceness and rationality for an SVOA. Recall also
from \S\ref{sec:SVOAMods} that a rational SVOA has finitely many
irreducible modules up to equivalence.

\subsubsection{Theta group}\label{sec:thetagp}

Let $\Gamma=\SL(2,\ZZ)$ and recall that the modular group
$\bar{\Gamma}=\Gamma/\{\pm 1\}=\PSL(2,\ZZ)$ acts faithfully on the
upper half plane $\hh=\{\sigma+\ii t\mid t>0\}\subset\CC$, with the
action generated by modular transformations $S$ and $T$ where
$S:\tau\mapsto-1/\tau$ and $T:\tau\mapsto \tau+1$. We identify
$\bar{\Gamma}$ with its image in the isometry group of $\hh$ and set
$\bar{\Gamma}_{\theta}=\langle S,T^2\rangle$. The compactification
of the quotient space $\bar{\Gamma} \backslash \hh$ is topologically
a sphere, and the same is true for $\bar{\Gamma}_{\theta} \backslash
\hh$. The space $\bar{\Gamma}_{\theta} \backslash \hh$ has two
cusps, with representatives $1$ and $\infty$. There is a unique
holomorphic function on $\hh$ that is invariant under
$\bar{\Gamma}_{\ogi}$, has a $q$ expansion of the form
$q^{-1/2}+a+bq^{1/2}+cq+\ldots$, and vanishes at $1$. We denote this
function by $J_{\ogi}(\tau)$ since it is an analogue of the $J$
function, which generates the field of functions on the compactified
curve $\bar{\Gamma}\backslash \hh$. The function $J_{\ogi}$
furnishes a bijective map from the compactification of
$\bar{\Gamma}_{\ogi}\backslash \hh$ to the Riemann sphere
$\CC\cup\{\infty\}$. One has the following expression for
$J_{\ogi}(\tau)$.
\begin{gather}
    \begin{split}
    J_{\ogi}(\tau)&=\frac{\eta(\tau)^{48}}
                    {\eta(\tau/2)^{24}\eta(2\tau)^{24}}\\
            &=q^{-1/2}+24+276q^{1/2}+2048q+11202q^{3/2}
                +49152q^{2}+          \cdots
    \end{split}
\end{gather}
To see the behavior of $J_{\ogi}$ at $1$, note that $TS\tau\to 1$ as
$\tau\to \infty$. For $J_{\ogi}|_{TS}$ we have
\begin{gather}
    \begin{split}
     J_{\ogi}|_{TS}=J_{\ogi}(-1/\tau+1)
          &=-2^{12}\frac{\eta(2\tau)^{24}}
                    {\eta(\tau)^{24}}\\
            &=-(4096q+98304q^{2}+1228800q^{3}+\ldots)
    \end{split}
\end{gather}
confirming that $J_{\ogi}$ vanishes as $\tau\to 1$. We write
$\Gamma_{\ogi}$ for the preimage of $\bar{\Gamma}_{\ogi}$ in
$\Gamma$.

%%%%%%%%%%%%%%%%%%%%%%%%%%%%%%%%%%%%%%%%%%%%%%%%%%%%%%%%
%%%%%%%%%%%%%%%%%%%%%%%%%%%%%%%%%%%%%%%%%%%%%%%%%%%%%%%%
\subsubsection{Modular Invariance}\label{sec:uniq:SVOA:modinv}

Suppose now that $(V,Y,\vac,\cas)$ is a nice rational VOA. Following
\cite{ZhuModInv} we may define a new VOA structure on the space $V$
as follows. We define the genus one VOA associated to $V$ to be the
four-tuple $(V,Y[\;],\vac,\tilde{\cas})$ where
$\tilde{\cas}=(2\pi\ii)^2(\cas-c/24)$ and the linear map
$Y[\;]:V\otimes V\to V((z))$ is defined so that
\begin{gather}
    Y[u,z]=\sum_{n\in\ZZ}u_{[n]}z^{-n-1}
        =Y(u,e^{2\pi\ii z}-1)e^{{\rm deg}(u)2\pi\ii z}
\end{gather}
for $u$ an $L(0)$-homogeneous element in $V$. The object thus
defined is again a VOA and is isomorphic to $(V,Y,\vac,\cas)$
\cite{ZhuModInv}. In particular the coefficients of
$Y[\tilde{\cas},z]$ define a representation of the Virasoro
algebra with central charge $c$. We write
\begin{gather}
    L[z]=Y[\tilde{\cas},z]=\sum_{\ZZ}L[n]z^{-n-2}
\end{gather}
and for $n\in\ZZ$, we set $V_{[n]}=\{u\in V\mid L[0]u=nu\}$. Note
that for $u$ an $L(0)$-homogeneous element in $V$ we have
\begin{gather}
    \begin{split}
    Y[u,z]=
        &\sum_{n}u_n(e^{2\pi\ii z}-1)^{-n-1}
        e^{{\rm deg}(u)2\pi\ii z}\\
        =&\sum_{n}u_n
        (2\pi\ii z+\tfrac{1}{2}(2\pi\ii z)^2+\ldots)^{-n-1}
        e^{{\rm deg}(u)2\pi\ii z}\\
        =&\sum_{n}
            (2\pi\ii)^{-n-1}u_nz^{-n-1}
            (1+\tfrac{1}{2}2\pi\ii z+\ldots)^{-n-1}
        e^{{\rm deg}(u)2\pi\ii z}
    \end{split}
\end{gather}
and in particular, $u_{[n]}=(2\pi\ii)^{-n-1}u_n
+\sum_{k>0}c_ku_{n+k}$ for some constants $c_k\in\CC$. If $u,v\in
V_1$ then $u_1v=\langle u|v\rangle\vac$ and $u_nv=0$ for $n>1$ so
we have the following
\begin{lem}\label{NRu1vlemma}
Let $V$ be a nice rational VOA and let $u,v\in V_1$. Then
$u_{[1]}v=-(4\pi^2)^{-1}\langle u|v\rangle\vac$.
\end{lem}

Recall the Eisenstein series $G_2(\tau)$ given by
\begin{gather}
    G_2(\tau)=\frac{\pi^2}{3}+\sum_{m\neq 0}\sum_n
        \frac{1}{(m\tau+n)^2}
\end{gather}
The function $G_2(\tau)$ has a $q$ expansion which may be expressed
in the form
\begin{gather}
    G_2(\tau)=
    \frac{\pi^2}{3}-8\pi^2
        \sum_{n=1}^{\infty}\sigma_1(n)q^n
\end{gather}
where $\sigma_1(n)$ is the sum of the divisors of $n$. We denote
the corresponding formal power series (that is, element of
$\CC[[q]]$) by $\tilde{G}_2(q)$.

We define a linear function $o(\cdot):V\to{\rm End}(V)$ by setting
$o(u)=u_{{\rm deg}(u)-1}$ for $L(0)$-homogeneous $u\in V$. The
following result is a special case of Proposition~4.3.5 in
\cite{ZhuModInv}.
\begin{prop}[\cite{ZhuModInv}]\label{NR_ZhuProp}
Let $V$ be a nice rational SVOA and $M$ a finitely generated
$V$-module. Then for $u,v\in V_1$ we have
\begin{gather}
    {\sf tr}|_Mo(u)o(v)q^{L(0)}=
    {\sf tr}|_Mo(u_{[-1]}v)q^{L(0)}
    -\tilde{G}_{2}(q)
    {\sf tr}|_M o(u_{[1]}v)q^{L(0)}
\end{gather}
\end{prop}

Let $V$ be a nice rational SVOA and let $M$ be a finitely generated
$V$-module. For an $n$-tuple $(u_1,\ldots,u_n)$ of
$L(0)$-homogeneous elements in $V$ we define the following formal
series.
\begin{gather}
\begin{split}
    &\tilde{F}_M((u_1,x_1),\ldots,(u_n,x_n);q)\\
        &=x_1^{{\rm deg}(u_1)}\cdots x_n^{{\rm deg}(u_n)}
        {\sf tr}|_MY(u_1,x_1)\cdots Y(u_n,x_n)q^{L(0)}
\end{split}
\end{gather}
We extend the definition of $\tilde{F}_M$ to arbitrary $n$-tuples
of elements from $V$ by linearity. As in \cite[Th
4.2.1]{ZhuModInv} One can show that this series $\tilde{F}_M$
converges to a holomorphic function in the domain
\begin{gather}
    \{ (x_1,\ldots,x_n,q)\mid 1>|x_1|>\ldots>|x_n|>|q|\}
\end{gather}
and can be continuously extended to be meromorphic in the domain
\begin{gather}
    \{ (x_1,\ldots,x_n,q)\mid x_i\neq 0,\,|q|<1\}
\end{gather}
We denote the meromorphic function so obtained by $F_M$. We
substitute variables $x_i$ with $e^{2\pi\ii z_i}$ and $q$ with
$e^{2\pi\ii \tau}$, and we set
\begin{gather}
    T_M((u_1,z_1),\ldots,(u_n,z_n);\tau)
        =q^{-c/24}F_M((u_1,x_1),\ldots,(u_n,x_n);q)
\end{gather}
Following \cite{ZhuModInv} and \cite{HohnPhD} we call
$T_M((u_1,z_1), \ldots, (u_n,z_n); \tau)$ the $n$-point
correlation function on the torus with parameter $\tau$ for the
operators $Y(u_i,z_i)$ and the module $M$.
\begin{prop}[\cite{ZhuModInv},\cite{HohnPhD}]
The function $T_M((u_1,z_1),\ldots,(u_n,z_n);\tau)$ is doubly
periodic in each variable $z_i$ with periods $1$ and $2\tau$, and
possible singularities only at the points $z_i=z_j+k+l\tau$ for
$i\neq j$, $k,l\in\ZZ$. For $j\in\{1,\ldots,n\}$ and $u_j$ an
$L(0)$-homogeneous element of $V$ we have
\begin{gather}
    T_M(\ldots,(u_j,z_j+\tau), \ldots; \tau)=
    (-1)^{p(u_j)}T_M(\ldots,(u_j,z_j), \ldots; \tau)
\end{gather}
For a permutation $\sigma\in S_n$ we have
\begin{gather}
    T_M((u_1,z_1),\ldots,(u_n,z_n);\tau)=(-1)^w
    T_M((u_{\sigma(1)},z_{\sigma(1)}),\ldots,
        (u_{\sigma(n)},z_{\sigma(n)});\tau)
\end{gather}
where $w$ is the number of permutations of the elements $u_i$ that
lie in $V_{\bar{1}}$.
\end{prop}

We denote by $T_M$ the mapping defined on the set
$\bigcup_{n=1}^{\infty}((V\times \CC)^n\times\hh)$ that sends
$((u_1,z_1),\ldots,(u_n,z_n);\tau)$ to
$T_M((u_1,z_1),\ldots,(u_n,z_n);\tau)$.

Suppose now that $\{M^1,\ldots, M^r\}$ is a complete list of
irreducible $V$-modules. The superconformal block on the torus
associated to the SVOA $V$ is the $\CC$-vector space spanned by
the $r$ mappings $T_{M^i}$. We denote it by ${\rm SB}_V$.

The following result is an analogue for SVOAs of a celebrated
theorem due to Zhu concerning the modularity properties of $n$-point
correlation functions on the torus associated to the vertex
operators on VOAs. As is indicated in \cite{HohnPhD}, this analogue
may be proven in a manner directly analogous to that of the VOA
version given in \cite{ZhuModInv}, and one should use the SVOA
analogues of Zhu algebras defined in \cite{KacWanSVOAs}.

\begin{thm}[\cite{ZhuModInv},\cite{KacWanSVOAs},\cite{HohnPhD}]
\label{Thm:thetagpinv}
Let $V$ be a nice rational SVOA and suppose that $\{M^1,\ldots,
M^r\}$ is a complete list of irreducible $V$-modules. Then the
superconformal block on the torus associated to $V$ is
$r$-dimensional and the functions $T_{M^i}$ form a basis. Moreover,
there exists a representation $\rho$ of $\Gamma_{\ogi}$ on ${\rm
SB}_V$ such that for Virasoro highest weight vectors
$\{u_1,\ldots,u_n\}\in V$ the $n$-point correlation functions on the
torus for the operators $Y(u_i,z_i)$ and the modules $M^i$ satisfy
the following transformation property
\begin{gather}
    \begin{split}
    &T_{M^i}\left(
        \left(u_1,\frac{z_1}{c\tau+d}\right),\ldots,
        \left(u_n,\frac{z_n}{c\tau+d}\right);
        \frac{a\tau+b}{c\tau+d}\right)\\
        &\qquad=
    (c\tau+d)^{\sum_k {\rm deg}(u_k)}
        \sum_j\rho(A)_{ij}
        T_{M^j}((u_1,z_1),\ldots,(u_n,z_n);\tau)
    \end{split}
\end{gather}
where $A=\binom{a\;b}{c\;d}$ is an element of $\Gamma_{\ogi}$ and
$(\rho(A)_{ij})$ is the matrix representing $\rho(A)\in{\rm
End}({\rm SB}_V)$ with respect to the basis $\{T_{M^i}\}$.
\end{thm}
In the case that $n=1$ the function $T_M((u,z);\tau)$ is elliptic in
the variable $z$ and without poles, and is therefore constant with
respect to $z$. We may therefore set $T_M(u;\tau)= T_M((u,z);\tau)$.
Note that $T_M(u;\tau)={\sf tr}|_Mo(u)q^{L(0)-c/24}$, and in
particular $T_M(\vac;\tau)={\sf tr}|_Mq^{L(0)-c/24}$.
\begin{cor}\label{SDNR_wtkModFrm}
Let $V$ be a self-dual nice rational SVOA. Then ${\rm SB}_V$ is one
dimensional, and for $u\in V_{\bar{0}}$ a Virasoro highest weight
vector with ${\rm deg}(u)=k$, the function ${\sf
tr}|_Mo(u)q^{L(0)-c/24}$ is a weight $k$ modular form on
$\bar{\Gamma}_{\ogi}$, possibly with character.
\end{cor}
\begin{rmk}
Corollary \ref{SDNR_wtkModFrm} now appears as a special case of the
more general Theorem~3 of \cite{DonZhaMdltyOrbSVOA} which
incorporates also $g$-twisted SVOA modules for $g$ in a finite group
of automorphisms of a suitable SVOA $V$.
\end{rmk}

We apply Corollary \ref{SDNR_wtkModFrm} immediately with $u=\vac$ in
order to determine the character of an SVOA satisfying our
hypotheses.
\begin{prop}\label{Prop:VChar}
Suppose that $V$ is a self-dual nice rational SVOA of rank $12$.
Then we have
\begin{gather}
\begin{split}
     {\sf tr}|_Vq^{L(0)-c/24}&=J_{\ogi}(\tau)
               +\dim(V_{1/2})-24\\
          &=q^{-1/2}+\dim(V_{1/2})+276q^{1/2}
               +2048q+11202q^{3/2}+\ldots
\end{split}
\end{gather}
for the character of $V$.
\end{prop}
\begin{proof}
Let us set $f(\tau)={\sf tr}|_Vq^{L(0)-c/24}$. By hypothesis,
$f(\tau)$ admits a Fourier expansion of the form
$q^{-1/2}+\sum_{n\geq 0}f_nq^{n/2}$ with all the $f_n$ non-negative
integers. By Corollary \ref{SDNR_wtkModFrm} we know that $f(\tau)$
is holomorphic on $\hh$ and is invariant for the action of
$\Gamma_{\ogi}$. It follows that $f(\tau)=P(J_{\ogi})/Q(J_{\ogi})$
for some polynomials $P(X),Q(X)\in \CC[X]$, with
$\deg(P)=\deg(Q)+1$, and we may assume that $P$ and $Q$ are both
monic and have no common factors. The function $J_{\ogi}(\tau)$ is a
surjective map from $\hh$ to $\CC\setminus\{0\}$, so that for
$f(\tau)$ to be holomorphic we must have $Q(X)=X^m$ for some $m$.
Then we have
\begin{gather}\label{eqn:fisLauJ}
     f(\tau)=J_{\ogi}+a_{m}+a_{m-1}J_{\ogi}^{-1}+\cdots+
          a_{0}J_{\ogi}^{-m}
\end{gather}
for $P(X)=X^{m+1}+a_mX^m+\cdots+a_0$ with $a_0\neq 0$ unless
possibly if $m=0$. We claim that $m=0$ and
$a_m=a_0=\dim(V_{1/2})-24$. Certainly, $a_m=\dim(V_{1/2})-24$, since
the first two terms of (\ref{eqn:fisLauJ}) determine the first two
Fourier coefficients of $f(\tau)$. Let us write
$J_{\ogi}(\tau)^{-d}=\sum_n r_{-d}(n)q^{n/2}$. Then the sequence
$\{r_{-d}(n)\}_n$ alternates in sign when $d$ is positive, as can be
seen from the following identity.
\begin{gather}
     J_{\ogi}(\tau)^{-d}
          =\frac{\eta({\tau}/{2})^{24d}\eta(2\tau)^{24d}}{\eta(\tau)^{48d}}
          =q^{d/2}\prod\frac{1}{(1+q^{n+1/2})^{24d}}
\end{gather}
%To see this, note that $J_{\ogi}(\tau)^{-d}$ is the supercharacter
%of the Weyl-module SVOA associated to a vector space of dimension
%$24d$ with non-degenerate antisymmetric form --- c.f.
%\cite{WeiSympPVOA}.
We see also from this that for $d>0$, the value of $|r_{-d}(n)|$ is
the number of partitions of $n-d$ into odd parts with $24d$ colors.
The asymptotic behavior of such functions is described by Theorem~1
of \cite{MeiAsympPtns}, and we will quote this result presently. The
value of $r_1(n)$ is the number of partitions of $n+1$ into odd
parts of $24$ colors without replacement, and for the asymptotics of
this function we refer to Proposition~1 of
\cite{HwaLimitThmsIntPtns}. The result is that we have
\begin{gather}
     r_1(n)\sim C_1\frac{e^{2\pi\sqrt{n}}}{n^{3/4}}\quad\text{and}
          \quad
     |r_{-d}(n)|\sim
     C_{-d}\frac{e^{2\pi\sqrt{2d}\sqrt{n}}}{n^{3/4}}\quad
          \text{for $d>0$,}
\end{gather}
for some constants $C_k$. Evidently, the growth of the $r_{-m}(n)$
outstrips that of the $r_{-d}(n)$ for $1\geq -d\geq-m+1$ if $m>0$.
In particular, the $f_n$ can be all non-negative integers only if
$m=0$.
\end{proof}
As demonstrated in \cite{DonZhaMdltyOrbSVOA}, one may recover a
modular invariance under the full modular group for the trace
functions associated to an SVOA by considering canonically twisted
modules together with untwisted modules. The following result is a
special case of Theorem~1 of \cite{DonZhaMdltyOrbSVOA} where we take
$V$ to be self-dual, and $G$ to be the group of SVOA automorphisms
generated by the canonical automorphism of $V$. Recall from
\S\ref{sec:SVOATwMods} that if $V$ is a self-dual rational
$C_2$-cofinite SVOA, then $V_{\sigma}$ denotes the unique
$\sigma$-stable $\sigma$-twisted $V$-module, and $\sigma$-stable
here means that $V_{\sigma}$ admits a compatible action by $\sigma$.
\begin{prop}[\cite{DonZhaMdltyOrbSVOA}]\label{prop:FullModInvSVOA}
Let $V$ be a self-dual rational $C_2$-cofinite SVOA. Let $w\in V$
such that $w\in V_{[k]}$ for some $k$. Then for $\gamma\in\Gamma$,
we have
\begin{gather}
     {\sf tr}|_{V}o(w)q^{L(0)-c/24}|_{\gamma}
     =(c\tau+d)^k\rho(\gamma)
     {\sf tr}|_{W}o(w)\sigma^{1+b+d}q^{L(0)-c/24}
\end{gather}
for some $\rho(\gamma)\in\CC$ independent of $w$, where
$W=V_{\sigma}$ if $\sigma^{1+a+c}=\sigma$, and $W=V$ otherwise, and
$\gamma$ is the matrix
     $\left(\begin{array}{cc}
          a & b \\
          c & d
     \end{array}\right)$.
\end{prop}
In Proposition~\ref{Prop:VChar} we determined that the character of
a self-dual nice rational SVOA of rank $12$ is
$J_{\ogi}(\tau)+\dim(V_{1/2})-24$. Applying Proposition~
\ref{prop:FullModInvSVOA} with $\gamma=TS$ (so that
$(a,b,c,d)=(1,-1,1,0)$) we find that the character ${\sf
tr}|_{V_{\sigma}}q^{L(0)-c/24}$ of the canonically twisted module
$V_{\sigma}$ over such an SVOA is just $\alpha
(J_{\ogi}|_{TS}+\dim(V_{1/2})-24)$ for some $\alpha\in\CC$.
Recalling from \S\ref{sec:thetagp} that the $q$ expansion of
$J_{\ogi}|_{TS}$ involves only positive integer powers of $q$, we
have
\begin{prop}\label{prop:twmodbnd}
Let $V$ be a self-dual nice rational SVOA of rank $12$ with
$V_{1/2}=0$. Then $(V_{\sigma})_n$ vanishes unless $n>0$ and
$n\in\ZZh$.
\end{prop}

\subsubsection{Structure of $V_1$}

Our plan now is to study the structure of the Lie algebra on $V_1$
for $V$ satisfying suitable hypotheses. By the end of this section,
knowledge of $V_1$ will determine $V$ uniquely under the conditions
we consider. Much of our method follows that employed in certain
sections of \cite{DonMasEfctCC} and \cite{DonMasHlmVOA}, and is a
manifestation of the principle established there that modular
invariance for an SVOA $V$ can be used to make strong conclusions
about the structure of $V_1$.

\begin{prop}\label{Prop:N=1NiceBiFm}
Suppose that $V$ is a nice $N=1$ SVOA with $V_{1/2}=0$. Then $V$ has
a unique non-degenerate invariant bilinear form.
\end{prop}
\begin{proof}
By the results of \cite{SchVASStgs}, the space of invariant bilinear
forms on $V$ is in natural correspondence with the space
$V_0/L(1)V_1$. We have that $V_0$ is one dimensional by hypothesis,
so we require to show that $L(1)V_1=0$. From the commutation
relations of the Neveu--Schwarz superalgebra we have
$G(\tfrac{1}{2})^2=L(1)$ so that $L(1)V_1\subset
G(\tfrac{1}{2})V_{1/2}$. Since $V_{1/2}=0$ the result follows.
\end{proof}

From now we assume $V$ to be a nice rational $N=1$ SVOA with
$V_{1/2}=0$. Then the import of Proposition~\ref{Prop:N=1NiceBiFm}
is that $V_{\bar{0}}$ is a strongly rational VOA in the sense of
\cite{DonMasEfctCC}. In particular Theorem~1 of \cite{DonMasEfctCC}
yields the following
\begin{thm}[\cite{DonMasEfctCC}]\label{Thm:V1Red}
The Lie algebra $V_1$ is reductive.
\end{thm}
Since $V_1$ is reductive, the Lie rank of $V_1$ is well defined.
Suppose in addition now that $V$ is self-dual. We will follow the
technique used to prove Theorem~2 in \cite{DonMasEfctCC} to
establish the following
\begin{thm}\label{Thm:RnkBndsLiernk}
The Lie rank of $V_1$ is bounded above by ${\rm rank}(V)$.
\end{thm}
\begin{proof}
We set $c={\rm rank}(V)$. Let $\gt{h}$ be a maximal abelian
subalgebra of $V_1$ consisting of semisimple elements. The Lie rank
of $V_1$ is the dimension of $\gt{h}$ and we denote this value by
$l$. The bilinear form $\lab\cdot|\cdot\rab$ restricts to be
non-degenerate on $\gt{h}$, and thus the vertex operators $Y(h,z)$
for $h\in\gt{h}$ generate an affine Lie algebra $\hat{\gt{h}}$, and
we can decompose $V$ as
\begin{gather}\label{Frm:TensDecompV}
    V=M(1)\otimes\Omega_V
\end{gather}
where $M(1)\simeq S(h_{-m}\mid h\in\gt{h},\,m>0)$ is the Heisenberg
VOA of rank $l$ associated to the space $\gt{h}$, and $\Omega_V$ is
the vacuum space consisting of vectors $u\in V$ such that $h_mu=0$
for all $h\in\gt{h}$ and $m>0$. Both factors on the right hand side
of (\ref{Frm:TensDecompV}) are invariant under the action of $L(0)$,
so that (\ref{Frm:TensDecompV}) holds even as a decomposition of
$L(0)$-graded spaces. Taking the trace of $q^{L(0)-c/24}$ on each
side of (\ref{Frm:TensDecompV}), multiplying both sides by
$\eta(q)^c$ and noting that ${\sf
tr}|_{M(1)}q^{L(0)}=q^{l/24}\eta(q)^{-l}$ we have
\begin{gather}\label{Frm:HolmMFrm}
    \eta(q)^c{\sf tr}|_Vq^{L(0)-c/24}
        =q^{(l-c)/24}\eta(q)^{c-l}
        {\sf tr}|_{\Omega_V}q^{L(0)}
\end{gather}
The expression on the left hand side of (\ref{Frm:HolmMFrm}) is a
holomorphic modular form on $\Gamma_{\ogi}$ (of weight $c/2$). It is
a classical result that the Fourier coefficients $r(n)$ say, of such
a function satisfy $r(n)= O(n^3)$. On the other hand, the Fourier
coefficients of $\eta(q)^{-s}$ grow like
$n^{-s/4-3/4}\exp(\pi\sqrt{2s/3}\sqrt{n})$ whenever $s>0$ (c.f.
\cite[Thm 1]{MeiAsympPtns}), so we must have $c-l\geq 0$. This is
what we required to show.
\end{proof}
\begin{rmk}
Certainly one expects Theorem~\ref{Thm:RnkBndsLiernk} to hold not
just for the case that $V$ is self-dual, but the present result is
sufficiently strong for our interests.
\end{rmk}

The following proposition is an analogue for self-dual rational
SVOAs of rank $12$ of Corollary 2.3 in \cite{DonMasHlmVOA}, and we
repeat here the method of proof used there.
\begin{prop}\label{SDNR_KFrmProp}
Suppose that $V$ is a self-dual nice rational $N=1$ SVOA of rank
$12$ with $V_{1/2}=0$. Then the Killing form on $V_1$ satisfies
$\kappa(\cdot\,,\cdot)=44\lab\cdot|\cdot\rab$.
\end{prop}
\begin{proof}
By Proposition~\ref{Prop:VChar} we have ${\sf
tr}|_Vq^{L(0)-c/24}=J_{\ogi}(\tau)-24$, and in particular,
$\dim(V_1)=276$. Now we apply Proposition~\ref{NR_ZhuProp} to $V$
with $M=V$ and use Lemma~\ref{NRu1vlemma} to rewrite the conclusion
as
\begin{gather}\label{PreKillingFrmEqn}
    {\sf tr}|_Vo(u)o(v)q^{L(0)-1/2}=
    {\sf tr}|_Vo(u_{[-1]}v)q^{L(0)-1/2}
    +\frac{\langle u|v\rangle}{4\pi^2}\tilde{G}_{2}(q)
    {\sf tr}|_V q^{L(0)-1/2}
\end{gather}
The leading term of the second summand on the right hand side of
(\ref{PreKillingFrmEqn}) is therefore $\tfrac{1}{12}\langle
u|v\rangle q^{-1/2}$, but the leading term on the left hand side of
(\ref{PreKillingFrmEqn}) is $\kappa(u,v)q^{1/2}$ where
$\kappa(\cdot\,,\cdot)$ is the Killing form on $V_1$. We conclude
that the leading term of the first summand on the right hand side of
(\ref{PreKillingFrmEqn}) is $-\tfrac{1}{12}\langle u|v\rangle
q^{-1/2}$. For $u,v\in V_1$ set $X_{u,v}(\tau)={\sf
tr}|_Vo(u_{[-1]}v)q^{L(0)-c/24}$. We claim that
$X_{u,v}(\tau)=\tfrac{1}{6}\lab u|v\rab qD_q{\sf
tr}|_Vq^{L(0)-c/24}$. This is certainly true if $\lab u|v\rab=0$.
Suppose not then, and note firstly that $X_{u,v}(\tau)$ is a weight
$2$ modular form for $\bar{\Gamma}_{\ogi}$ by Corollary
\ref{SDNR_wtkModFrm}, and secondly that $X_{u,v}(\tau)
=-\tfrac{1}{12}\lab u|v\rab q^{-1/2}+0+aq^{1/2}+\ldots$ for some
$a\in\CC$ by the above. Applying Proposition~
\ref{prop:FullModInvSVOA} with $w=u_{[-1]}v$ and $\gamma=TS$, we
find
\begin{gather}
     {X_{u,v}(-1/\tau+1)}{\tau^{-2}}
     =\alpha\,{\sf tr}|_{V_{\sigma}}
          o(u_{[-1]}v)q^{L(0)-c/24}%\\
\end{gather}
for some $\alpha$ in $\CC$, and this $q$ series belongs to
$\CC[[q]]$ by Proposition~\ref{prop:twmodbnd}. From this we note
thirdly, that $X_{u,v}(\tau)$ is holomorphic at the cusp represented
by $1$. We claim that these three properties determine
$X_{u,v}(\tau)$ uniquely, for if $X'(\tau)$ is another such
function, then $Z(\tau)=X_{u,v}(\tau)-X'(\tau)$ is a weight two
modular form for $\bar{\Gamma}_{\ogi}$ that is holomorphic at both
cusps, and vanishes at $\infty$. The space of weight $2$ modular
forms that are holomorphic at both cusps is spanned by the theta
function of the lattice $\ZZ^4$ (c.f. \cite{RanMdlrFrmsFns}), but
this function does not vanish at $\infty$, and the claim follows. It
is easy to check that $\tfrac{1}{6}\lab u|v\rab qD_q{\sf
tr}|_Vq^{L(0)-c/24}$ satisfies the three properties of
$X_{u,v}(\tau)$ so we may rewrite (\ref{PreKillingFrmEqn}) as
follows.
\begin{gather}
    \begin{split}
    {\sf tr}|_Vo(u)o(v)q^{L(0)-1/2}
    =&
    \frac{\langle u|v\rangle}{6}
        qD_q{\sf tr}|_Vq^{L(0)-1/2}
        +\frac{\langle u|v\rangle}{4\pi^2}\tilde{G}_2(q)
            {\sf tr}|_Vq^{L(0)-1/2}\\
    =&
    \frac{\langle u|v\rangle}{6}
        (-\tfrac{1}{2}q^{-1/2}+\tfrac{1}{2}276q^{1/2}+\ldots)\\
    &+\frac{\langle u|v\rangle}{4\pi^2}
        (\tfrac{\pi^2}{3}-8\pi^2q+\ldots)
        (q^{-1/2}+276q^{1/2}+\ldots)
    \end{split}
\end{gather}
Equating the coefficients of $q^{1/2}$ on each side we have
$\kappa(\cdot\,,\cdot)= 44\langle\cdot|\cdot\rangle$. %Since the
%bilinear form on $V$ is non-degenerate, the Killing form on $V_1$ is
%non-degenerate. This completes the proof.
\end{proof}

\begin{thm}\label{UniqSVOA}
Let $V$ be a self-dual nice rational $N=1$ SVOA with rank $12$ and
$V_{1/2}=0$. Then $V$ is isomorphic to $_{\CC}\afn$ as an SVOA.
\end{thm}
\begin{proof}
By Proposition~\ref{Prop:N=1NiceBiFm} the bilinear form defined by
the adjoint operators is non-de\-generate, and Theorems
\ref{Thm:V1Red} and \ref{Thm:RnkBndsLiernk} show that $V_1$ is a
reductive Lie algebra with Lie rank bounded above by $12$. From
Proposition~\ref{SDNR_KFrmProp} we find that $V_1$ is of dimension
$276$ and the Killing form $\kappa(\cdot\,,\cdot)$ on $V_1$
satisfies $\kappa(\cdot\,,\cdot)=44\lab\cdot|\cdot\rab$. In
particular, the Killing form is non-degenerate, and $V_1$ is a
semi-simple Lie algebra.

Suppose then that $\gt{g}$ is a simple component of $V_1$ with level
$k$ and dual coxeter number $h$. By the main theorem of
\cite{DonMasItgbtyVOAs} we have that $k$ is an integer. Suppose that
$(\cdot\,,\cdot)$ is the bilinear form on $\gt{g}$ normalized so
that $(\alpha,\alpha)=2$ for a long root $\alpha$. Then we have
$\lab u|v\rab=k(u,v)$ for $u,v \in\gt{g}$, and thus also
$\kappa(u,v)=44k(u,v)$ for $u,v \in\gt{g}$. Taking $u=v=\alpha$ we
obtain $h/k=22$ since $\kappa(\alpha,\alpha)=4h$. This argument is
independent of the choice of simple component and so the ratio $h/k$
must hold for each simple component. By inspection the only
possibility then is that $\gt{g}$ is of type $D_{12}$ with level
$k=1$.

Thus we find that $V_1$ is a semisimple Lie algebra of type
$D_{12}$, and the VOA $V_{\bar{0}}$ is isomorphic to the lattice VOA
$V_{M_0}$ for $M_0$ a copy of the $D_{12}$ lattice. It follows then
that the SVOA $V$ is isomorphic to a lattice VOA $V_M$ for some
positive definite integral lattice $M=M_0\cup M_1$. Since $V$ is
self-dual of rank $12$, $M$ is self-dual of rank $12$, and the fact
that $V_{1/2}=0$ implies that $M$ has no vectors of unit length.
There is one such lattice up to isomorphism; namely, the lattice
$D_{12}^+$. We conclude that any self-dual nice rational SVOA with
rank $12$ and trivial degree $1/2$ subspace is isomorphic to
$V_{M}$, where $M$ is a copy of the lattice $D_{12}^+$. From the
proof of Proposition~\ref{AfnIsSVOA} we see that the SVOA
$_{\CC}\afn$ is also such an object, and this completes the proof of
the theorem.
\end{proof}

\begin{rmk}
We sketch here an alternative approach to Theorem~\ref{UniqSVOA}
that was described to us by Gerald H\"ohn. Suppose that $V$ is as in
the statement of Theorem~\ref{UniqSVOA}, and let us write $U_0$ for
a copy of the lattice VOA associated to the lattice of type $D_4$.
This VOA has three irreducible modules beyond itself; we pick one of
them and denote it $U_1$. We then set $W=U_0\otimes
V_{\bar{0}}\oplus U_1\otimes V_{\bar{1}}$, so that $W$ is a module
for the VOA $U_0\otimes V_{\bar{0}}$ with only integral weights.
Using knowledge of the fusion of modules for $V_{\bar{0}}$ and $U_0$
it can be shown that the VOA structure extends in a unique way from
$U_0\otimes V_{\bar{0}}$ to the whole space. Then $W$ is a self-dual
VOA of rank $16$, and one may invoke Theorem~2 of
\cite{DonMasHlmVOA} to conclude that $W$ is a lattice VOA $W_L$ for
$L$ one of the two self-dual lattices of rank $16$; namely,
$E_8\oplus E_8$ or $D_{16}^+$. One then shows that the $D_4$ lattice
VOA $U_0$ can only be embedded in $W_L$ in such a way that $V_0$
must also be a lattice VOA, and the lattice must be of type
$D_{12}$.
\end{rmk}

\subsection{$N=1$ structure}\label{sec:Uniq_N=1Struc}

We now wish to demonstrate that the $N=1$ structure on $\afn$ is
unique. More precisely, suppose that $t\in \Cm(\gt{l})_G^0$
satisfies $\langle e_St,t\rangle=0$ for any $S\subset\Omega$ with
$0<w(S)\leq 4$. Citing Proposition~\ref{Prop:CodeCriForSC} as
justification, we call such a vector superconformal. We wish to show
that if $t$ is superconformal with $|t|=1$ then $t=x1_G$ for some
$x\in\Sp(\gt{l})$. This will be achieved in Theorem~\ref{ThmUniq},
after we establish a few preliminary lemmas. For the benefit of the
reader we now include a few words about the idea behind the proof of
this theorem.

\subsubsection{Strategy}\label{Sec:UniqStrat}

Our strategy is the following. Suppose that $t\in\Cm(\gt{l})_G^0$ is
a superconformal vector of unit norm, and define a function
$f_t:\Sp(\gt{l})\to [-1,1]$ by $f_t(x)=\langle x1_G,t\rangle$. Since
the bilinear form on $\Cm(\gt{l})_G^0$ is non-degenerate we are done
as soon as we find an $x\in \Sp(\gt{l})$ such that $f_t(x)=1$. We
will show that for any $x\in \Sp(\gt{l})$ with $f_t(x)<1$ there
exists some $x'\in \Sp(\gt{l})$ such that $f_t(x')>f_t(x)$. The
function $f$ is certainly continuous and $\Sp(\gt{l})$ is compact,
so showing that $f_t$ can always be made closer to $1$ suffices to
show that $f_t$ attains the value $1$. In fact, since $xt$ is
superconformal for $x\in \Sp(\gt{l})$ whenever $t$ is, it suffices
to show only that for any superconformal $t$ there is some $x\in
\Sp(\gt{l})$ such that $f_{t}(x)>f_t({\bf 1})$ where ${\bf 1}$
denotes the identity in $\Sp(\gt{l})$. Thus the following results up
to and including Theorem~\ref{ThmUniq} are dedicated to showing that
for any superconformal $t$ in $\Cm(\gt{l})_G^0$ there is some
$x\in\Sp(\gt{l})$ such that $\lab x1_G,t\rab>\lab 1_G,t\rab$.

We next recall some facts about the Golay co-code, and then
introduce some useful notation and terminology before presenting
Propositions \ref{Prop:dbevcntn} and \ref{Prop:dbevcntnlift}, which
will be the main tools we use to implement the stated strategy. A
superconformal vector $t$ may be regarded as an element of the unit
ball in $2048$ dimensional space, and Proposition~
\ref{Prop:dbevcntn} provides a way of regarding $t$ as an element of
the unit ball in $2048/|\Gamma|$ dimensional space via a kind of
linearization over the cosets of certain subgroups
$\Gamma<\mc{G}^*$. The Proposition~\ref{Prop:dbevcntnlift} is a
generalization of this result which arises essentially because
$\mc{G}^*$ has many distinct lifts to $\FF_2^{\Omega}$. We sometimes
refer to Propositions \ref{Prop:dbevcntn} and
\ref{Prop:dbevcntnlift} as the coset contraction results.

\subsubsection{Golay co-code}\label{Sec:Golayco-code}

Recall that the Golay co-code is the space
$\mc{G}^*=\FF_2^{\Omega}/\mc{G}$, and recall the co-weight function
$w^*$ on $\FF_2^{\Omega}$ from \S\ref{Sec:Notation}. We write
$X\mapsto \bar{X}$ for the canonical map
$\FF_2^{\Omega}\to\mc{G}^*$. %The weight function
%$w:\FF_2^{\Omega}\to\{0,\ldots,24\}$ induces a function $w^*$ on
%$\mc{G}^*$ defined so that $w^*(\bar{X})$ is the minimum weight
%among all code words in $\bar{X}=X+\mc{G}$. We call $w^*$ the
%co-weight function, and lift $w^*$ to a function on $\FF_2^{\Omega}$
%by defining $w^*(X)=w^*(\bar{X})$ for $X\in\FF_2^{\Omega}$.
The Golay code corrects three errors and the range of the co-weight
function is the set $\{0,1,2,3,4\}$. Let $X\in\FF_2^{\Omega}$. We
say that $X$ and $\bar{X}$ are co-even if $2|w^*(X)$, and we say
that $X$ and $\bar{X}$ are doubly co-even if $4|w^*(X)$.

If $w^*(X)=2$ then there is a unique pair of points $i,j\in\Omega$
such that $X+\mc{G}=\{ij\}+\mc{G}$. If $w(X)=w^*(X)=4$ then there
are exactly five weight $8$ words (octads) in $\mc{G}$ containing
$X$, and we have $w(Y)=w^*(Y)=4$ and $\bar{X}=\bar{Y}$ just when
$X+Y$ is one of these. Thus for $\bar{X}\in\mc{G}^*$ with
$w^*(X)=4$, the six sets of cardinality four in $\FF_2^{\Omega}$
that lift $\bar{X}$ constitute a partition of $\Omega$ into six
disjoint four sets. Such a partition is called a sextet, and the
four sets in a given sextet are called tetrads.

Let $\bar{T},\bar{T}'\in\mc{G}^*$ with co-weight $4$, and let
$S=\{T_i\}$ and $S'=\{T_i'\}$ be the sextets determined by $\bar{T}$
and $\bar{T}'$ respectively. Then there are essentially four
different ways that the sextets $S$ and $S'$ can overlap.
\begin{enumerate}
\item   $|T_i\cap T_j'|\in\{0,4\}$ for all $i,j$.
\item   $|T_i\cap T_j'|\in\{0,1,3\}$ for all $i,j$.
\item   $|T_i\cap T_j'|\in\{0,1,2\}$ for all $i,j$.
\item   $|T_i\cap T_j'|\in\{0,2\}$ for all $i,j$.
\end{enumerate}
The first case is just the case that $S$ and $S'$ are the same
sextet. The second case is the case that $\bar{T}+\bar{T}'$ has
co-weight $2$. In the third and fourth cases $\bar{T}+\bar{T}'$ has
co-weight $4$, and the last case is distinguished by the property
that any lift of the set $\{\bar{T},\bar{T}'\}$ is isotropic. We
refer to the sextets $S$ and $S'$ as commuting sextets when
$|T_i\cap T_j'|$ is even for all $i$ and $j$, and we refer to them
as non-commuting otherwise. In this setting, we refer to tetrads
$T_0$ and $T_0'$ as commuting when the corresponding sextets
$\{T_i\}$ and $\{T_i'\}$ are commuting, and we refer to them as
non-commuting otherwise.

\subsubsection{Co-code lifts}

Let $\Delta$ be a subset of $\FF_2^{\Omega}$ such that the map
$\FF_2^{\Omega}\to \mc{G}^*$ induces a bijection
$\Delta\leftrightarrow\bar{\Delta}=\{\bar{X}\mid X\in\Delta\}$. We
then call $\Delta$ a lift of $\bar{\Delta}$ to $\FF_2^{\Omega}$. If
further we have that $w(X)=w^*(X)$ for all $X\in\Delta$, we say that
$\Delta$ is a balanced lift of $\bar{\Delta}$. Recall that the space
$\FF_2^{\Omega}$ admits a bilinear form
$\FF_2^{\Omega}\times\FF_2^{\Omega}\to\FF_2$ defined so that $\lab
X,Y\rab=|X\cap Y|\pmod{2}$. In some cases a subset
$\bar{\Delta}<\mc{G}^*$ has a lift $\Delta\subset\FF_2^{\Omega}$
such that $\lab X,Y\rab\equiv 0$ for any $X,Y\in\Delta$, and we call
such a lift isotropic. Suppose that $\Delta$ is doubly co-even and
balanced. We then say that $\Delta$ is commuting if the sextets
determined by any pair of tetrads in $\Delta$ are commuting in the
sense of \S\ref{Sec:Golayco-code}.

Suppose now that $\Sigma$ is a balanced lift of $(\mc{G}^*)^0$,
the even part of $\mc{G}^*$. Then the set $\Sigma$ is in natural
bijective correspondence with $(\mc{G}^*)^0$, and the group
structure on the later may be lifted via this correspondence so as
to define a group structure on the former. We denote this group
operation with $\dotplus$, so that $X\dotplus Y=Z$ just when
$X+Y+Z\in\mc{G}$. Note that the bilinear form on $\FF_2^{\Omega}$
is not bilinear with respect to $\dotplus$, so that in general
$\lab A\dotplus B,X\rab\neq \lab A,X\rab+\lab B,X\rab$ for
example. A multiplicative $2$-cocycle $\sigma$ with values in
$\{\pm 1\}$ is defined on the group $\Sigma=(\Sigma,\dotplus)$ by
requiring that $e_Xe_Y1_G=\sigma(X,Y)e_{X\dotplus Y}1_G$ for
$X,Y\in\Sigma$.
\begin{prop}\label{prop:2cocycRHsymm}
We have
\begin{gather}
    \sigma(X,X)=(-1)^{w^*(X)/2}\\
    \sigma(X,Y)=(-1)^{\lab X,Y\rab}\sigma(Y,X)\\
    \sigma(X,Y)=\sigma(X,X)\sigma(X,X\dotplus Y)
\end{gather}
for all $X,Y\in\Sigma$. In particular
$\sigma(X,Y)=\sigma(X,X\dotplus Y)$ for all $Y\in\Sigma$ just when
$X$ is doubly co-even.
\end{prop}
\begin{proof}
For any $X,Y\in\Sigma$ we have $e_Xe_Y1_G=(-1)^{\lab
X,Y\rab}e_Ye_X1_G$, and this implies $\sigma(X,Y)=(-1)^{\lab
X,Y\rab}\sigma(Y,X)$. Left multiplying both sides of $e_Xe_Y1_G
=\sigma(X,Y)e_{X\dotplus Y}1_G$ by $e_X$ yields $\sigma(X,Y)
e_Xe_{X\dotplus Y}1_G =e_X^2e_Y1_G$. On the other hand
$e_Xe_{X\dotplus Y}1_G=\sigma(X,X\dotplus Y)e_Y1_G$. Since
$e_X^2=\sigma(X,X)$, this verifies the claim.
\end{proof}
When $\Sigma$ is a lift of $(\mc{G}^*)^0$, the set
$\{e_X1_G|X\in\Sigma\}$ constitutes an orthonormal basis for
$\Cm(\gt{l})_G^0$. We then have $t=\sum_{\Sigma}t_Xe_X1_G$ for
unique $t_X\in\RR$ such that $\sum t_X^2=1$ when $t$ is a unit
vector in $\Cm(\gt{l})_G^0$. Note that $f_t({\bf
1})=t_{\emptyset}$. We write ${\rm supp}(t)$ for the set of
$\bar{X}\in\mc{G}^*$ such that $t_X\neq 0$.

One way to obtain a balanced lift of $(\mc{G}^{*})^0$ is the
following. Choose an element in $\Omega$ and denote it by
$\infty$. Let $\Delta=\Delta_0\cup\Delta_2\cup\Delta_4$ where
$\Delta_0$ contains just the emptyset, $\Delta_2$ is the set of
pairs of elements from $\Omega$, and $\Delta_4$ is the set of
subsets of $\Omega$ of size four containing $\infty$. Then
$\Delta$ is a balanced lift of $(\mc{G}^*)^0$. A doubly co-even
subgroup $\Gamma<\Delta$ is then isotropic just when it is
commuting.

\subsubsection{Coset contraction}

Let $\Sigma$ be a balanced lift of $(\mc{G}^*)^0$, suppose that
$\bar{\Gamma}$ is a subgroup of $\mc{G}^*$, and suppose that
$W\in\Sigma$ is chosen so that the coset $\bar{W}+\bar{\Gamma}$ of
$\bar{\Gamma}$ in $\mc{G}^*$ is doubly co-even. Suppose also that
the corresponding lift $W\dotplus \Gamma\subset\Sigma$ obtained by
restriction from $\Sigma$ is isotropic. Then the $2$-cocycle
$\sigma$ is symmetric on $W\dotplus \Gamma$. Let $\chi:\Sigma\to
\pm 1$ be a function such that
\begin{gather}\label{form:chicond}
    \chi(A\dotplus W)\sigma(A\dotplus W,Z)=\chi(Z)\chi(A\dotplus
        W\dotplus Z),\quad \forall A\in\Gamma,\;
            Z\in\Sigma.
\end{gather}
Then by Proposition~\ref{prop:2cocycRHsymm} we have
$\sigma(A\dotplus W ,Z)=\sigma(A\dotplus W,A\dotplus W\dotplus Z)$
so that the invariance under swapping $Z$ with $A\dotplus W\dotplus
Z$ is evident for both sides of the expression (\ref{form:chicond}).
The assumption that $W\dotplus \Gamma$ be isotropic ensures that
(\ref{form:chicond}) can be satisfied for $Z$ in $W\dotplus \Gamma$.

In the case that $W\dotplus \Gamma=\Gamma$ the condition
(\ref{form:chicond}) implies that the restriction $\chi|_{\Gamma}$
of $\chi$ to $\Gamma$ is a
$1$-cocycle with %values in $\{\pm 1\}$, and
coboundary $\sigma|_{\Gamma\times\Gamma}$ since we have
$\sigma(A\dotplus B,A) =\sigma(A,B)$ for $A,B \in\Gamma$ when
$\Gamma$ is isotropic. The values of $\chi$ on a coset $Z\dotplus
\Gamma$ of $\Gamma$ in $\Sigma$ are determined by those on
$\Gamma$ together with $\chi(Z)$ since $\chi(Z\dotplus
A)=\chi(A)\sigma(A,Z)/\chi(Z)$. Any two such functions
$\chi:\Sigma\to\pm 1$ therefore differ by a single element of
$\Gamma^*={\rm Hom}(\Gamma,\pm 1)$ on each coset of $\Gamma$ in
$\Sigma$. Recall that the maps $\bar{X}\mapsto(-1)^{\lab D,X\rab}$
for $D\in\mc{G}$ exhaust the homomorphisms $\mc{G}^*\to\pm 1$.
Indeed, $(-1)^{\lab D,X\rab}$ is independent of the choice of lift
$X$ for $\bar{X}$, and thus any homomorphism $(\Sigma,\dotplus)\to
\pm 1$ is of the form $X\mapsto (-1)^{\lab D,X\rab}$ for some
$D\in\mc{G}$.

To a function $\chi$ satisfying (\ref{form:chicond}) and to any
given $Z\in \Sigma$ we associate the element $u_{\chi,Z}
=\sum_{A\in\Gamma} \chi(Z\dotplus A) e_{Z\dotplus A}$ in
$\Cl(\gt{l})$, and for ease of notation we set
$u_{\chi}=u_{\chi,\emptyset}$. The following proposition is our
main tool for classifying superconformal vectors in
$\Cm(\gt{l})_G^0$.
\begin{prop}\label{Prop:dbevcntn}
Let $t=\sum_{\Sigma}t_Xe_X1_G$ be superconformal with $|t|=1$. Let
$\bar{\Gamma}$ be a subgroup of $\mc{G}^*$ and suppose that
$\bar{W}\in\mc{G}^*$ is chosen so that $\bar{W}+ \bar{\Gamma}$ is
doubly co-even. Suppose also that $W\dotplus\Gamma \subset\Sigma$
is an isotropic lift of $\bar{W}+\bar{\Gamma}$. Then for
$\chi:\Sigma\to \pm 1$ satisfying (\ref{form:chicond}) and for
${\mathsf T}$ any transversal of $\Gamma$ in $\Sigma$ we have
\begin{gather}\label{eqn:cosetlinz}
    \sum_{Z\in{\mathsf T}}
        \lab u_{\chi,Z}1_G,t\rab
        \lab u_{\chi,W\dotplus Z}1_G,t\rab
        =\lab u_{\chi,W}t,t\rab
        =\begin{cases} 1&\text{if $W\in\Gamma$,}\\
            0&\text{if $W\notin\Gamma$.}
            \end{cases}
\end{gather}
In particular, for the case that ${\rm supp}(t)$ is a doubly
co-even group with an isotropic lift $\Gamma$ we have
$\sum_{A\in\Gamma} \chi(A)t_A=\pm 1$ for any $1$-cocycle
$\chi:\Gamma\to \pm 1$ with coboundary
$\sigma|_{\Gamma\times\Gamma}$.
\end{prop}
\begin{proof}
Let us consider the expression $\lab u_{\chi}t,t\rab$. Since $t$
is superconformal, we have $\lab e_Xt,t\rab=0$ whenever $w(X)<8$,
so that $\lab u_{\chi}t,t\rab=\chi(\emptyset)\lab t,t\rab=1$. On
the other hand we have
\begin{gather}
    \begin{split}
        u_{\chi}t&=\sum_{A\in\Gamma, Z\in\Sigma}
            \chi(A)t_Ze_Ae_Z1_G\\
        &=\sum_{A\in\Gamma,Z\in \Sigma}
                \chi(A)\sigma(A,Z)
                t_Ze_{Z\dotplus A}1_G\\
        &=\sum_{A\in\Gamma,Z\in\Sigma}
            \chi(Z)\chi(Z\dotplus A)
                t_Ze_{Z\dotplus A}1_G
    \end{split}
\end{gather}
and then $1=\lab u_{\chi}t,t\rab =\sum \chi(Z)\chi(Z\dotplus A)
t_Zt_{Z\dotplus A}$. From the fact that $\lab
u_{\chi,Z}1_G,t\rab=\sum_{A\in\Gamma}\chi(Z\dotplus A)t_{Z\dotplus
A}$ we see that the left hand side of (\ref{eqn:cosetlinz})
coincides with $\lab u_{\chi}t,t\rab$, and the equation
(\ref{eqn:cosetlinz}) follows. This handles the case that
$W\in\Gamma$, and the case that $W\notin\Gamma$ is similar. The
case that all $t_X$ vanish for $X\notin\Gamma$ then yields
$(\sum_{\Sigma}\chi(X)t_X)^2=1$, and the last part follows from
this.
\end{proof}

Suppose that $\Sigma$ and $\Sigma'$ are balanced lifts of
$\mc{G}^*$ to $\FF_2^{\Omega}$. We denote the group operations on
$\Sigma$ and $\Sigma'$ both by $\dotplus$. There is a
correspondence $\Sigma\leftrightarrow \Sigma'$ such that we have
$X\leftrightarrow X'$ if and only if
$\bar{X}=\bar{X}'\in\mc{G}^*$. We then have $X\dotplus Y=Z$ in
$\Sigma$ just when $X'\dotplus Y'=Z'$ in $\Sigma'$. For $X\in
\Sigma$ we have $\lab X,X'\rab=0$ and $X+X'\in\mc{G}$. Define
$\varrho:\mc{G}^*\to\pm 1$ so that $\varrho_{\bar{X}}e_Xe_{X'}1_G
=1_G$. We abuse notation to write $\varrho_X$ for
$\varrho_{\bar{X}}$ whenever $X\in\Sigma$. Let $\sigma'$ denote
the multiplicative $2$-cocycle on $(\Sigma',\dotplus)$ such that
$e_{X'}e_{Y'}1_G=\sigma'(X',Y')e_{Z'}1_G$ for $Z'=X'\dotplus Y'$.
The following lemma gives the relationship between $\sigma$ and
$\sigma'$.
\begin{lem}\label{Lem:sigrelsig'}
For $X,Y\in\Sigma$ and $Z=X\dotplus Y$ we have
\begin{gather}
    \sigma'(X',Y')=(-1)^{\lab X+X',Y\rab}
        \varrho_X\varrho_Y\varrho_Z
        \sigma(X,Y)
\end{gather}
\end{lem}
As before we assume that $W\dotplus \Gamma$ is doubly co-even and
isotropic, and we assume now the same for $W'\dotplus\Gamma'$. We
also assume that $\lab X',Y\rab=\lab X,Y'\rab$ for all $X,Y\in
W\dotplus \Gamma$. Given $\chi:\Sigma\to\pm 1$ satisfying
(\ref{form:chicond}), we define $\chi':\Sigma'\to\pm 1$ by setting
$\chi'(X')=\psi(X)\chi(X)$ for $X'\in\Sigma'$ where
$\psi:\Sigma\to \pm 1$ is chosen to satisfy
\begin{gather}\label{Form:LiftTrans}
    \psi(X)
    (-1)^{\lab X+X',Z\rab}
    \varrho_{X}\varrho_Z
    \varrho_{X\dotplus Z}
    =\psi(Z)\psi(X\dotplus Z),\quad
    \forall X\in W\dotplus\Gamma,  Z\in\Sigma.
\end{gather}
Then from Lemma~\ref{Lem:sigrelsig'} we obtain that
$\chi'(X')\sigma'(X',Z') =\chi'(Z')\chi'(X'\dotplus Z')$ whenever
$X'\in W'\dotplus \Gamma'$ and $Z'\in\Sigma'$. Now we may apply
Proposition~\ref{Prop:dbevcntn} with $\Sigma'$ in place of $\Sigma$,
and $\chi'$ in place of $\chi$, and we obtain the following
generalization of that proposition.
\begin{prop}\label{Prop:dbevcntnlift}
Let $\Sigma$, $\Sigma'$ be balanced lifts of $\mc{G}^*$, and
suppose that $\Gamma<\Sigma$ and $W\in\Sigma$ are chosen so that
both $W\dotplus \Gamma$ and $W'\dotplus \Gamma'$ are doubly
co-even and isotropic. Suppose also that $\lab X',Y\rab=\lab
X,Y'\rab$ for all $X,Y\in W\dotplus\Gamma$. Then for
$\chi:\Sigma\to\pm 1$ satisfying (\ref{form:chicond}), for
$\psi:\Sigma\to \pm 1$ satisfying (\ref{Form:LiftTrans}), and for
${\mathsf T}$ a transversal of $\Gamma$ in $\Sigma$, we have
\begin{gather}
    \sum_{Z\in{\mathsf T}}
        \lab u_{\chi'',Z}1_G,t\rab
        \lab u_{\chi'',W\dotplus Z}1_G,t\rab
        =\begin{cases}
            1&\text{if $W\in\Gamma$,}\\
            0&\text{if $W\notin\Gamma$.}
        \end{cases}
\end{gather}
where $\chi'':\Sigma\to \pm 1$ is given by
$\chi''(X)=\varrho_X\psi(X)\chi(X)$, and $u_{\chi'',Z}$ is defined
by $u_{\chi'',Z} =\sum_{A\in\Gamma}\chi''(A\dotplus Z)e_{A\dotplus
Z}$ for $Z\in\Sigma$.
\end{prop}

\subsubsection{Superconformal vectors}

We now embark upon the task of realizing the strategy summarized in
\S\ref{Sec:UniqStrat}; that is, the task of showing that for any
superconformal $t\in\Cm(\gt{l})_G^0$ there is some $x$ in
$\Sp(\gt{l})$ such that $f_t(x)>t_{\emptyset}$. In practice, we
treat all possible unit vectors $t$ on a case by case basis using
the coset contraction results to narrow down the possibilities for
the coefficients $t_X$ of $t=\sum_{\Sigma}t_Xe_X1_G$ (given a
balanced lift $\Sigma$ of $(\mc{G}^*)^0$ to $\FF_2^{\Omega}$) that
can make $t$ superconformal. In the course of doing so we find
superconformal vectors in the $\Sp(\gt{l})$ orbit of $1_G$ other
than those of the form $\pm e_X1_G$ for $X\in\FF_2^{\Omega}$, or
$\exp(re_X)1_G$ for $r\in\RR$ and $w(X)=2$, and ultimately we find
that $t$ either has projection on one of these vectors exceeding
$t_{\emptyset}$ or cannot be superconformal.

Suppose that $|t_X|>t_{\emptyset}$ for some $X\in\Sigma$. Then we
have $\lab x1_G,t\rab=\sigma(X,X)t_X$ for $x=e_X$. Multiplying by
$-1$ if necessary, we have $f_t(x)>t_{\emptyset}$. Thus from now
on we may suppose that $t_{\emptyset}\geq|t_X|$ for all
$X\in\Sigma$.

Suppose that $t_X\neq 0$ for some $X$ with $w^*(X)=2$. Setting
$\exp(re_X)=\cos(r)+\sin(r)e_X$ we have $\lab\exp(re_X)1_G,t\rab=
(1-\tfrac{1}{2}r^2)t_{\emptyset}-rt_X+o(r^2)$ so that
$f_t(x)>t_{\emptyset}$ for $x=\exp(re_X)$ and suitably chosen $r$.
From now on we assume that $t_X=0$ whenever $w^*(X)=2$. That is we
may assume that ${\rm supp}(t)$ is a doubly co-even subset of
$\mc{G}^*$.

Suppose that ${\rm supp}(t)$ is contained in a doubly co-even
subgroup $\bar{\Gamma}$ of $\mc{G}^*$ and suppose that
$\bar{\Gamma}$ has a balanced isotropic lift $\Gamma$. We may assume
that $\Sigma$ is a balanced lift of $(\mc{G}^*)^0$ containing
$\Gamma$. Since it is useful, we now state the following result,
which is obtained by direct application of
Proposition~\ref{Prop:dbevcntnlift} to our present situation.
\begin{prop}\label{Prop:dbevsupplift}
Let $\bar{\Gamma}$ be a doubly co-even subgroup of $\mc{G}^*$ and
suppose that $\Gamma$ and $\Gamma'$ are balanced isotropic lifts
of $\bar{\Gamma}$ such that $\lab A',B\rab=\lab A,B'\rab$ for all
$A,B\in \Gamma$. Then for $\chi:\Gamma\to\pm 1$ a $1$-cocycle with
coboundary $\sigma|_{\Gamma\times\Gamma}$, and for $\psi:\Gamma\to
\pm 1$ satisfying
\begin{gather}
    \psi(A)
    (-1)^{\lab A+A',B\rab}
    \varrho_{A}\varrho_B
    \varrho_{A\dotplus B}
    =\psi(B)\psi(A\dotplus B)
\end{gather}
for all $A,B\in\Gamma$ we have $\lab u_{\chi''}1_G,t\rab^2=1$
where $\chi'':\Gamma\to \pm 1$ is given by
$\chi''(A)=\varrho_A\psi(A)\chi(A)$, and
$u_{\chi''}=\sum_{A\in\Gamma}\chi''(A)e_{A}$.
\end{prop}
The requirement that $\bar{\Gamma}$ be doubly co-even is quite
strong. Any maximal doubly co-even subgroup of $\mc{G}^*$ has order
$16$ or $32$, and every doubly co-even subgroup of order $16$ or
less has an isotropic lift. A convenient way to generate doubly
co-even subgroups of $\mc{G}^*$ is the following. Choose a weight
$12$ word in $\mc{G}$ (a dodecad), and partition the $12$ non-zero
coordinates into six pairs $\{A_i\}$. Then the $15$ elements
$\bar{A}_i+\bar{A}_j \in\mc{G}^*$ are the non-trivial elements in a
doubly co-even subgroup $\bar{\Gamma}$ say, of $\mc{G}^*$ of order
$16$. Furthermore, the set $\{\emptyset,A_i+A_j\}$ furnishes a
balanced isotropic lift of $\Gamma$. For some partitions there is a
sextet $S=\{T_i\}$ such that each pair $A_i$ is contained in a
tetrad of $S$. In this case, the addition of one of the $T_i$
extends $\Gamma$ to be a balanced isotropic lift of a doubly co-even
$32$ group in $\mc{G}^*$.

Suppose then that ${\rm supp}(t)$ is contained in a two group
$\bar{\Gamma}$ with balanced lift $\Gamma=\{\emptyset,A\}$. Then
since $t$ is a unit vector we have $t_{\emptyset}^2+t_{A}^2=1$. On
the other hand from Proposition~\ref{Prop:dbevsupplift} we have
$(t_{\emptyset}\pm t_A)^2= 1$ and thus $t_{\emptyset}+t_A=\pm 1$ and
$t_{\emptyset}-t_A=\pm 1$. The only solutions are $t_{\emptyset}=\pm
1$ and $t_{A}=\pm 1$. Since we have assumed $t_{\emptyset}\geq
|t_X|$ for all $X\in\Sigma$, we have $t_{\emptyset}=1$ and $t=1_G$.

Suppose now that ${\rm supp}(t)$ is contained in a four group
$\bar{\Gamma}$ with balanced isotropic lift
$\Gamma=\{\emptyset,A,B,C\}$. Similar to before we have
$t_{\emptyset}^2+t_A^2+t_B^2+t_C^2=1$. A function
$\chi_0:\Gamma\to\pm 1$ suitable for an application of Proposition~
\ref{Prop:dbevsupplift} may be given arbitrary values on generators
of $\Gamma$, and then the remaining values are determined by
$\sigma$. For example, we may set
$\chi_0(\emptyset)=\chi_0(A)=\chi_0(B)=1$ and
$\chi_0(C)=\sigma(A,B)\chi_0(A)\chi_0(B)=\sigma(A,B)$. Any other
suitable $1$-cocycle $\chi$ differs from $\chi_0$ by some element of
$\Gamma^*$, and Proposition~\ref{Prop:dbevsupplift} (with
$\Sigma'=\Sigma$) now yields
\begin{gather}\label{Form:dbevsuppcntnon4}
    t_{\emptyset}+\mu(A)t_A+\mu(B)t_B
        +\sigma(A,B)\mu(C)t_C=\pm 1
\end{gather}
where $\mu$ is any homomorphism $\mu:\Gamma\to\pm 1$. Summing over
(\ref{Form:dbevsuppcntnon4}) for various choices of $\mu$ we find
that $t_X\pm t_Y\in\{0,\pm 1\}$ for each $X,Y\in\Gamma$, and then
$t_X\in\{0,\pm \tfrac{1}{2},\pm 1\}$ for each $X\in\Gamma$. Thus if
one of the $t_X$ vanishes then the remaining $t_Y$ must lie in
$\{0,\pm 1\}$, and then $t=\pm t_Y$ for some $Y\in\Gamma$. If one of
the $t_X$ is $\pm \tfrac{1}{2}$ then the remaining $t_Y$ must lie in
$\pm \tfrac{1}{2}$ and there is some restriction on the signs: any
two solutions differ by an element of $\Gamma^*$, and one solution
is given by $s=\tfrac{1}{2}(1+e_A+e_B-\sigma(A,B)e_C)1_G$. Taking
$\Sigma'$ different from $\Sigma$ we can say more: if there is some
lift $\Gamma'$ of $\bar{\Gamma}$ such that $\lab X',Y\rab\neq 0$ for
some $X,Y\in\Gamma$, then we obtain another equation like
(\ref{Form:dbevsuppcntnon4}) with signs not differing by an element
of $\Gamma^*$, and this extra restriction is enough to rule out the
possibility that any $t_X$ has norm $\tfrac{1}{2}$. Such a lift
$\Gamma'$ exists just when two of the sextets corresponding to
$\{\bar{A},\bar{B},\bar{C}\}$ are non-commuting (see
\S\ref{Sec:Golayco-code}). We are left with the question of whether
of not a vector $s\in\Cm(\gt{l})_G^0$ of the form $s =\tfrac{1}{2}(1
+\mu(A)e_A+\mu(B)e_B-\sigma(A,B)\mu(C)e_C)1_G$ with $\mu\in\Gamma^*$
is in the $\Sp(\gt{l})$ orbit containing $1_G$ given that the
sextets in $\bar{\Gamma}$ are commuting. The answer is affirmative,
as the following lemma demonstrates.
\begin{lem}\label{Lem:sc4gp}
Suppose ${\Gamma}=\{\emptyset,A,B,C\}$ is a balanced lift of a
commuting four group in $\mc{G}^*$, and
\begin{gather}
    s =\frac{1}{2}\left(1 +\mu(A)e_A +\mu(B)e_B -\sigma(A,B)
    \mu(C)e_C\right)1_G
\end{gather}
for some $\mu\in\Gamma^*$. Then $s$ is in the $\Sp(\gt{l})$ orbit
containing $1_G$.
\end{lem}
\begin{proof}
For $S\in\mc{G}$ let $g_S=\pm e_S$ with the sign chosen so that
$g_S\in G$. For any given $S\in\mc{G}$ either $\lab S,X\rab=0$ for
all $X\in\Gamma$, or there is a unique non-trivial $X\in\Gamma$ such
that $\lab S,X\rab=0$. We define a group $G'=\{g_S'\mid
S\in\mc{G}\}< \Sp(\gt{l})$ by setting $g_S'=s_Xe_Xg_S$ when $X$ is
the unique non-trivial element of $\Gamma$ such that $\lab
S,X\rab=0$, and setting $g_S'=g_S$ when $\lab S,X\rab =0$ for all
$X\in\Gamma$. Then a simple computation shows that $g_S's=s$ for all
$S\in\mc{G}$. The group $G'$ is an $\FF_2^{\mc{E}}$-homogeneous
subgroup of $\Sp(\gt{l})$ whose associated code is doubly-even
self-dual and has no short roots. In other words, $G'$ is a lift of
a Golay code on $\Omega$. Noting that both $G$ and $G'$ contain the
volume element $e_{\Omega}$, It follows from the uniqueness of the
Golay code that there is some coordinate permutation in
$\Sp(\gt{l})$ that sends $s$ to $1_G$.
\end{proof}
The method illustrated above for the case that $\bar{\Gamma}$ has
order four is a model for the cases of higher order. For this
reason we will summarize only the results for the higher order
cases that we need, and refrain from burdening the reader with all
details.
\begin{lem}\label{Lem:sc8gp}
In the case that ${\rm supp}(t)$ is contained in an eight group
$\bar{\Gamma}$ with balanced isotropic lift $\Gamma= \lab
A,B,C\rab$, either ${\rm supp}(t)$ is contained in a commutative
four group, or $\Gamma$ is totally commutative and $t$ is of the
form
\begin{gather}
    t=\frac{1}{4}\left(3
    +\mu(A)e_A+\mu(B)e_B+\mu(C)e_C
    +\mu(A\dotplus B)\sigma(A,B)e_{A\dotplus B}
    +\ldots\right)1_G
\end{gather}
for some $\mu\in\Gamma^*$, and $t$ belongs to the $\Sp(\gt{l})$
orbit containing $1_G$.
\end{lem}
\begin{lem}\label{Lem:sc16gp}
In the case that ${\rm supp}(t)$ is contained in a $16$ group
$\bar{\Gamma}$ with balanced isotropic lift $\Gamma$, either ${\rm
supp}(t)$ is contained in some commutative eight group, or there
is a dodecad in $\mc{G}$ and a partition $P$ of its non-trivial
coordinates into six pairs $P=\{A_0,\ldots,A_5\}$ such that
$\Gamma=\{\emptyset,A_{ij}\}$ where we write $A_{ij}$ for the
tetrad $A_i+A_j\in\FF_2^{\Omega}$. We may assume that
$e_{A_0}e_{A_1}\cdots e_{A_5}\in G$. Then
$\sigma(A_{ij},A_{ik})=-1$ and $\sigma(A_{ij},A_{kl})=1$ for
distinct $i,j,k,l$. The vector $t$ is then of the form
\begin{gather}
    t=\frac{1}{4}
    \left( 1+\mu(A_{01})e_{A_{01}}
            +\ldots
            +\mu(A_{04})e_{A_{04}}
            -\mu(A_{12})e_{A_{12}}
            -\cdots-\mu(A_{45})e_{A_{45}}
        \right)1_G
\end{gather}
for some $\mu\in\Gamma^*$, and $t$ belongs to the $\Sp(\gt{l})$
orbit containing $1_G$.
\end{lem}
\begin{lem}\label{Lem:sc32gp}
In the case that ${\rm supp}(t)$ is contained in a $32$ group with
balanced isotropic lift, either ${\rm supp}(t)$ is contained in
some $16$ group with isotropic lift, or $f_t(x)>t_{\emptyset}$ for
$x1_G=s$ a superconformal vector with ${\rm supp}(s)$ contained in
a doubly co-even $16$ group with isotropic lift.
\end{lem}

We must now treat the case that ${\rm supp}(t)$ is not contained
in any doubly co-even group with isotropic lift.  We remind here
that we assume $t_{\emptyset}\geq |t_X|$ for all $X\in\Sigma$, and
$t_X=0$ whenever $w^*(X)=2$. We claim that for such $t$, either
$f_t(x)>t_{\emptyset}$ for some $x1_G=s$ with $s$ supported on a
doubly co-even group with isotropic lift as given in Lemmata
\ref{Lem:sc4gp} through \ref{Lem:sc16gp}, or $t$ is not
superconformal.

So let us assume that $f_t(x)\leq t_{\emptyset}$ for any
$x\in\Sp(\gt{l})$ such that $x1_G=s$ for $s$ one of the
superconformal vectors appearing in Lemmata \ref{Lem:sc4gp} through
\ref{Lem:sc16gp}. This condition amounts to putting upper bounds on
the moduli of the coefficients $t_X$ that are non-zero. For example,
take $s$ and $\Gamma$ as in Lemma~\ref{Lem:sc4gp}. Then $\lab
s,t\rab\leq t_{\emptyset}$ for all $\mu\in\Gamma^*$ is equivalent to
the inequalities
\begin{gather}
    0\leq \frac{1}{2}\left(
        t_{\emptyset}+\mu(A)t_A+\mu(B)t_B
        +\mu(C)\sigma(A,B)t_C\right),\quad
        \forall\mu\in\Gamma^*,
\end{gather}
which in turn imply that the smallest of $|t_A|$, $|t_B|$ and
$|t_C|$ is not greater than $\tfrac{1}{4}$, given that all are
non-zero. Also, one can construct elements $x\in\Sp(\gt{l})$ of the
form $x=\exp(\theta_1e_{X_1})\cdots\exp(\theta_ke_{X_k})$ for
$w(X_i)=2$ such that $f_t(x)>t_{\emptyset}$ so long as not all the
non-vanishing $t_X$ in $t$ are too small. On the other hand,
Proposition~\ref{Prop:dbevcntn} applied in the case that $W\notin
\Gamma$ can be used to show that the non-vanishing of some co-weight
$4$ coefficients $t_X$ implies the non-vanishing of others. The
simplest result of this kind is the following.
\begin{lem}\label{Lem:FourGpCond}
If $t_A\neq 0$ then there is some $B\in\Sigma$ such that
$t_Bt_{A\dotplus B}\neq 0$.
\end{lem}
\begin{proof}
We take $\Gamma=\{\emptyset\}$ and $W=X$ in Proposition~
\ref{Prop:dbevcntn}. We then obtain
\begin{gather}
    0=\langle e_{A}t,t\rangle
    =\sum_{X\in\Sigma}\left\langle
        t_Xe_Ae_X1_G,t\right\rangle
    =2t_{\emptyset}t_A+
    \sum_{X\in\Sigma\setminus\emptyset,A}
    \sigma(A,X)t_Xt_{A\dotplus X}
\end{gather}
and this implies the claim.
\end{proof}
We require to show then that the coefficients $t_X$ for $w^*(X)=4$
cannot all be too small. Since some non-vanishing co-weight $4$
components implies more non-vanishing co-weight $4$ components, let
us consider the extreme case that $t_X\neq 0$ for all $X\in\Sigma$
with $w^*(X)=4$. Suppose that $t_Z=\varepsilon$ is the greatest
among these (we may assume $t_Z>0$) so that $\varepsilon\geq|t_X|$
for all $X\in\Sigma$ with $w^*(X)=4$. Then since $\sum t_X^2=1$ we
have $t_{\emptyset}>t_{\emptyset}^2> 1-N\varepsilon^2$ where
$N=1771$ is the number of co-weight $4$ elements in $\mc{G}^*$. On
the other hand Proposition~\ref{Prop:dbevcntn} with
$\Gamma=\emptyset$ and $W=Z$ yields $0=\lab e_Zt,t\rab$, and we then
have
\begin{gather}
    0=\lab e_Zt,t\rab> 2(1-N\varepsilon^2)\varepsilon
        -M\varepsilon^2
\end{gather}
where $M$ is the number of $X\in\Sigma$ with $w^*(X)=4$ such that
$w^*(Z\dotplus X)=4$. We have $2(1-N\varepsilon^2)\varepsilon
-M\varepsilon^2>0$ (that is, a contradiction) just when
$2>M\varepsilon+2N\varepsilon^2$, so that $\varepsilon$ cannot be
smaller than $1/\sqrt{2N}$ for example. In this way we find that any
$t$ which does not satisfy $f_t(x)>t_{\emptyset}$ for some
superconformal $x1_G=s$ already constructed is not superconformal.
That is, we have the following

\begin{thm}\label{Thm:UniqOrb}
The superconformal vectors in $\Cm(\gt{l})_G$ form a single orbit
under the action of $\Sp(\gt{l})$. This orbit contains $1_G$.
\end{thm}

From Theorems \ref{UniqSVOA} and \ref{Thm:UniqOrb} we deduce the
following characterization of $\afn$.
\begin{thm}\label{ThmUniq}
Let $V$ be a self-dual nice rational $N=1$ SVOA with rank $12$ and
$V_{1/2}=0$. Then $V$ is isomorphic to $_{\CC}\afn$ as an $N=1$
SVOA.
\end{thm}

For $x\in\Sp_{24}(\RR)$ let us write ${\rm tr}|_{24}x$ for the trace
of $x$ in the representation of $\Sp_{24}(\RR)$ on $\RR^{24}$. We
have ${\rm tr}|_{24}(e_X)=24-2n$ when $w(X)=n$. Combining Theorems
\ref{Thm:PtStabIsCo0} and \ref{Thm:UniqOrb} we obtain the following
characterization of the group $\Co_0$.
\begin{thm}\label{Thm:Co0Chrztn}
Let $M$ be a spin module for $\Sp_{24}(\RR)$ and let $t\in M$ such
that $\lab xt,t\rab=0$ whenever $x\in\Sp_{24}(\RR)$ is an involution
with ${\rm tr}|_{24}x=16$. Then the subgroup of $\Sp_{24}(\RR)$
fixing $t$ is isomorphic to $\Co_0$.
\end{thm}

\section{Structure of $\vfn$}\label{LatConst}

In this section we summarize the construction of $\vfn$, mentioned
in the introduction, and we indicate how to construct an explicit
isomorphism with $\afn$.

\subsection{Lattice $N=1$ SVOAs}

There is a standard construction which assigns an $N=1$ SVOA to a
positive definite integral lattice, and we summarize that
construction now. Suppose that $L$ is a positive definite integral
lattice. Let $\FF$ be $\RR$ or $\CC$, and recall from
\S\ref{sec:SVOAstruc:LattSVOAs}, the SVOA $_{\FF}V_L$ associated to
$L$ via the standard construction. Let $_{\FF}\gt{a}
=\FF\otimes_{\ZZ}L$, and let us denote $A(_{\FF}\gt{a})$ by
$_{\FF}A_L$. Define $_{\FF}V^f_L$ to be the tensor product of SVOAs
\begin{gather}
    _{\FF}V^f_L=\,_{\FF}A_L\otimes\,_{\FF}V_L
\end{gather}
For $\FF=\CC$ the SVOA $_{\CC}V^f_L$ admits a natural structure of
$N=1$ SVOA. To define the superconformal element, let $h_i$ be an
orthonormal basis of $_{\CC}\gt{h}=\CC\otimes_{\ZZ} L$ and let $e_i$
be the corresponding elements of $_{\CC}\gt{a}$ under the
identification $_{\CC}\gt{a}=\,_{\CC}\gt{h}=\CC\otimes_{\ZZ} L$.
Then we set $\tau$ to be the element in $(_{\CC}V^f_L)_{3/2}$ given
by
\begin{equation}
    \tau=\frac{\ii}{\sqrt{8}}\sum_ie_i(-\tfrac{1}{2})h_i(-1)\vac
    %\in (V^f_L)_{3/2}
\end{equation}
where we suppress the tensor product from our notation. From
\cite{KapOrlCTorus}, \cite{SchVASStgs} for example, we have the
following
\begin{thm}
The element $\tau$ is a superconformal vector for $_{\CC}V^f_L$. In
particular, $_{\CC}V^f_L$ admits a natural structure of $N=1$ SVOA.
\end{thm}
Just as in \S\ref{RealFormLatSVOA} we can obtain a real form
${V}_{L}^f$ for $_{\CC}V^f_L$ by setting ${V}^f_L=
{_{\RR}A_L}\otimes {V}_L$. Noting that $\ii h_i(-1)\in\ii
{_{\RR}V^1_L}\subset {V}_L$ we see that $\tau\in{V}^f_L$, and the
$N=1$ structure on $_{\CC}V_L^f$ restricts so as to furnish an $N=1$
structure on ${V}_L^f$.

For simplicity, let us suppose that the rank of $L$ is even. By
the Boson-Fermion correspondence \cite{FreBF}, \cite{DoMaBF}, we
have an isomorphism of SVOAs $_{\CC}A_L\cong{_{\CC}V_{\ZZ^n}}$
where $n={\rm rank}(L)/2$, so that $_{\CC}A_L$ is self-dual as an
SVOA. The tensor product $_{\CC}V^f_L$ is therefore isomorphic to
a lattice SVOA $_{\CC}V_{\ZZ^n\oplus L}$. The lattice $\ZZ^n\oplus
L$ is self-dual just when $L$ is, so we conclude that the $N=1$
SVOA associated to any self-dual lattice is self-dual as an SVOA.
More generally, the irreducible $_{\CC}V^f_L$ modules are indexed
by the cosets of $L$ in its dual.

\subsection{The case that $L=E_8$}

From now on we take $L$ be a lattice of $E_8$ type so that
$_{\CC}V_L^f$ is a realization of the $N=1$ SVOA associated to the
$E_8$ lattice. Since $L$ is a self-dual lattice, $_{\CC}V^f_L$ is
a self-dual $N=1$ SVOA. The idea is that the $N=1$ SVOA
$_{\CC}\vfn$ should be a $\ZZ/2$-orbifold of $_{\CC}V^f_L$. More
particularly, we wish to define the space underlying $_{\CC}\vfn$
to be $(_{\CC} V^f_L)^0 \oplus (_{\CC}V^f_L)_{\ogi}^0$ where
$\ogi$ is a suitably chosen involution on $_{\CC}V_L^f$, the space
$(_{\CC}V^f_L)_{\ogi}$ is an $\ogi$-twisted $_{\CC}V_L^f$-module,
and the superscripts outside the brackets indicate $\ogi$-fixed
points. In order to construct $_{\CC}\vfn$ we must therefore
specify the involution $\ogi$, and construct a $\ogi$-twisted
module. This will be the objective of the next two subsections.

\subsection{Twisting}

The SVOA $_{\CC}V^f_L$ admits an automorphism $\ogi$ given by
$\ogi=\ogi_f\otimes\ogi_b$ where $\ogi_f$ denotes the parity
involution on $_{\CC}A_L$, and $\ogi_b$ denotes a lift of $-1$ on
$L$ to $\Aut(_{\CC}V_L)$. Observe that both $\ogi_f$ and $\ogi_b$
may be regarded as a lift of $-1$ on $L$. We have
$\ogi(\scas_V)=\scas_V$, so $\ogi$ is an automorphism of the $N=1$
structure on $_{\CC}V^f_L$. Also, the real form of $_{\CC}V^f_L$ is
just
\begin{gather}
    {V}^f_L=
        \left\{u+\ii v\mid u\in(_{\RR}V^f_L)^0,\;
        v\in(_{\RR}V^f_L)^1\right\}\subset\,_{\CC}V^f_L
\end{gather}
where $_{\FF}V^f_L=(_{\FF}V^f_L)^0\oplus (_{\FF}V^f_L)^1$
indicates the decomposition into $\ogi$-eigenspaces.

An $\ogi$-twisted $_{\CC}V^f_L$-module $(_{\CC}V^f_L)_{\ogi}$ is
of the form
\begin{gather}
    (_{\CC}V^f_L)_{\ogi}=(_{\CC}A_L)_{\ogi_f}
        \otimes(_{\CC}V_L)_{\ogi_b}
\end{gather}
where $(_{\CC}A_L)_{\ogi_f}=A(_{\CC}\gt{a})_{\ogi_f}$ is a
canonically twisted $_{\CC}A_L$-module and may be constructed as in
\S\ref{sec:cliffalgs:SVOAs}, and $(_{\CC}V_L)_{\ogi_b}$ is a
$\ogi_b$-twisted module over $_{\CC}V_L$. There is a well known
method for constructing $\ogi_b$-twisted $_{\CC}V_L$ modules for
$\ogi_b$ a lift of $-1$ on $L$, and one may refer to \cite{FLM} for
a thorough treatment. It turns out that for the case we are
interested in there is a simpler approach using only modules over
lattice VOAs, and this in turn can be viewed from the point of view
of Clifford module SVOAs. Such an approach is convenient for our
purpose.

Recall that the lattice $L$ contains a sublattice of the form
$K\oplus K$ where $K$ is a lattice of $D_4$ type. Let $K^*$ denote
the dual lattice to $K$, and let $K_{\gamma}$ for
$\gamma\in\Gamma=\{0,1,\w,\wc\}$ be an enumeration of the cosets of
$K$ in $K^*$. We decree that $K_0=K$. The remaining cosets
$K_{\gamma}$ for $\gamma\neq 0$ are permuted by automorphisms of
$K^*$ preserving $K$, and for this reason it is natural to regard
$\Gamma$ as a copy of the field of order $4$. We may assume that the
lattice $L$ decomposes as $L=\bigcup_{\Gamma} K_{\gamma}\oplus
K_{\gamma}$ into cosets of $K\oplus K$. Then the VOA $_{\CC}V_L$ has
a decomposition
\begin{gather}\label{KDecompVL}
    _{\CC}V_L=\bigoplus_{\gamma\in\Gamma}
        \,_{\CC}V_{K_{\gamma}}
        \otimes\,_{\CC}V_{K_{\gamma}}
\end{gather}
The VOA $_{\CC}V_K$, being a VOA of $D_4$ type, may be realized
using Clifford module SVOAs, and similarly for its modules
$_{\CC}V_{K_{\gamma}}$. In fact we may take $_{\CC}V_K$ to be a
copy of $A(_{\CC}\gt{a})^0$, and then $\bigoplus_{\Gamma}\,
_{\CC}V_{K_{\gamma}}$ is isomorphic as an
$A(_{\CC}\gt{a})^0$-module to the space $A(_{\CC}\gt{a})\oplus
A(_{\CC}\gt{a})_{\ogi_f}$.

The corresponding construction of $_{\CC}V_L$, using
$A(_{\CC}\gt{a})^0$-modules in place of $_{\CC}V_K$-modules, was
achieved in \cite{FFR}. Indeed, they provided more than this,
describing certain twisted modules over $_{\CC}V_L$ and proving that
the direct sum of these twisted and untwisted structures may be
equipped with a certain generalization of VOA structure; namely
para-VOA structure. It turns out that the twisting involutions
considered in \cite{FFR} are conjugate to the involution $\ogi_b$
under the action of $\Aut(_{\CC}V_L)\cong E_8(\CC)$, and in
particular, we may use one of these in place of $\ogi_b$. Before
describing precisely the involution we will use, we set up some new
notation, and recall the relevant results of \cite{FFR}.

\subsection{Clifford construction of $_{\CC}V_L$}

Recall that $_{\FF}\gt{a}=\FF\otimes_{\ZZ}L$ for $\FF=\RR$ or $\CC$.
The extended Hamming code is the unique up to equivalence
doubly-even self-dual code of length $8$. Let $\Pi$ be some set of
cardinality $8$, and let $\mc{H}$ be a copy of the extended Hamming
code, which we regard at once as a subset of $\mc{P}(\Pi)$, and as a
subspace of $\FF_2^{\Pi}$. Let $\mc{E}=\{e_i\}_{i\in\Pi}$ be an
orthonormal basis for ${_{\RR}\gt{a}}\subset {_{\CC}\gt{a}}$, and
let $H$ be an $\FF_2^{\mc{E}}$-homogeneous lift of $\mc{H}$ to
$\Sp(_{\RR}\gt{a})$. We may then realize an $\ogi_f$-twisted
$A(_{\FF}\gt{a})$-module explicitly by setting
$A(_{\FF}\gt{a})_{\ogi_f}= A(_{\FF}\gt{a})_{\ogi_f,H}$. We define
$_{\FF}U_0$ to be the VOA $A(_{\FF}\gt{a})^0$ and we enumerate the
$_{\FF}U_0$-modules $_{\FF}U_{\gamma}$ for $\gamma\in\Gamma$ by
setting
\begin{gather}
    _{\FF}U_0=A(_{\FF}\gt{a})^0,\quad
    _{\FF}U_1=A(_{\FF}\gt{a})^1,\quad
    _{\FF}U_{\w}=A(_{\FF}\gt{a})_{\ogi_f}^0,\quad
    _{\FF}U_{\wc}=A(_{\FF}\gt{a})_{\ogi_f}^1.
\end{gather}
Then $_{\CC}U_0$ is isomorphic to the $D_4$ lattice VOA, and the
$_{\CC}U_{\gamma}$ are its irreducible modules. We set
$_{\FF}U=\bigoplus_{\gamma\in\Gamma} {_{\FF}U_{\gamma}}$. From
\cite{FFR} we have $_{\CC}U_0$-module intertwining operators
$I_{\gamma\delta}: {_{\CC}U_{\gamma}}\otimes {_{\CC}U_{\delta}}
\to {_{\CC}U_{\gamma+\delta}}((z^{1/2}))$ such that the map
$I=(I_{\gamma\delta}): {_{\CC}U}\otimes {_{\CC}U} \to
{_{\CC}U}((z^{1/2}))$ furnishes $_{\CC}U$ with a structure of
para-VOA. We refer the reader to \cite{FFR} for detailed
information about para-VOAs, and we note here that the restriction
of $I$ to $_{\CC}U_0\oplus {_{\CC}U_{\gamma}}$ equips that space
with a structure of SVOA for any $\gamma\neq 0$.

For $k\in\{1,2,3\}$ let $_{\FF}\gt{a}^k$ be a copy of the space
$_{\FF}\gt{a}$ with orthonormal basis
$\mc{E}^k=\{e_i^k\}_{i\in\Pi}$, and let $_{\FF}U^k$ be a copy of
the space $_{\FF}U$. Suppose we define spaces $_{\FF}W_L$ and
$_{\FF}W_L'$ by setting
\begin{gather}
    _{\FF}{W}_L =\bigoplus_{\gamma\in\Gamma}
        {_{\FF}U^2_{\gamma}} \otimes
        {_{\FF}U^3_{\gamma}},\qquad
    _{\FF}W_L'  =\bigoplus_{\gamma\in\Gamma}
        {_{\FF}U^2_{\gamma}}\otimes
        {_{\FF}U^3_{\gamma+\omega}}.
\end{gather}
The main result from \cite{FFR} that we will use is that the
para-VOA structure on $_{\CC}U$ induces a VOA structure on
$_{\CC}W_L$ isomorphic to $_{\CC}V_L$, and induces a structure of
$\ogi_b'$-twisted $_{\CC}W_L$-module on $_{\CC}W_L'$, where
$\ogi_b'=1\otimes\ogi_f$. Furthermore, we may assume that the
isomorphism is chosen so that the action of $\ogi_b'$ on
$_{\CC}W_L$ corresponds to that of $\ogi_b$ on $_{\CC}V_L$.
Indeed, a Cartan subalgebra of $_{\CC}({W}_L)_1$ is spanned by the
$\ii e^2_i(-\tfrac{1}{2}) e^3_i(-\tfrac{1}{2})$ for ${i\in\Pi}$,
and $\ogi_b'$ acts as $-1$ on this space.

From now on we will regard the VOAs $_{\CC}W_L$ and $_{\CC}V_L$ as
identified via some VOA isomorphism such that $\ogi_b$ corresponds
to $\ogi_b'$, and we will write $_{\CC}V_L$ in place of
$_{\CC}W_L$ and $\ogi_b$ in place of $\ogi_b'$. Then for a
$\ogi_b$-twisted $_{\CC}V_L$-module we may take
$(_{\CC}V_L)_{\ogi_b}={_{\CC}W_L'}$. Note that $\ogi_b$ acts
naturally on the $\ogi_b$-twisted module $(_{\CC}V_L)_{\ogi_b}$. A
real form ${V}_L$ for $_{\CC}V_L$ may be described by
${V}_L={_{\RR}W_L}=\bigoplus_{\Gamma} {_{\RR}U_{\gamma}^2} \otimes
{_{\RR}U_{\gamma}^3}$.

We may now express the spaces $_{\CC}V_L^f$ and
$(_{\CC}V_L^f)_{\ogi}$ in the following way as sums of tensor
products of the $_{\CC}U^k_{\gamma}$.
\begin{gather}
    _{\CC}{V}^f_L
        =({_{\CC}U^1_0}\oplus {_{\CC}U^1_{1}})\otimes
        \left(\bigoplus_{\Gamma} {_{\CC}U^2_{\gamma}}
            \otimes {_{\CC}U^3_{\gamma}}
        \right)\\
    ({_{\CC}{V}^f_L})_{\ogi}
        =({_{\CC}U^1_{\omega}}\oplus
            {_{\CC}U^1_{\bar{\omega}}})\otimes
        \left(\bigoplus_{\Gamma} {_{\CC}U^2_{\gamma}}
        \otimes {_{\CC}U^3_{\gamma+\omega}}\right)
\end{gather}
We obtain real forms ${V}_L^f$ and $({V}_L^f)_{\ogi}$ by replacing
$\CC$ with $\RR$ in the subscripts of all the
$_{\CC}U^k_{\gamma}$. Note that the super-conformal element
$\tau_V\in{_{\CC}V^f_L}$ may now be written in the following form.
\begin{gather}
    \tau_V=-\frac{1}{\sqrt{8}}\sum_{\Pi}
        e^1_i(-\tfrac{1}{2})
        e^2_i(-\tfrac{1}{2})
        e^3_i(-\tfrac{1}{2})\vac
\end{gather}
\begin{rmk}
One can see that the Clifford module SVOA $A(\gt{u})$ has an $N=1$
structure whenever ${\rm dim}(\gt{u})$ is divisible by $3$.
\end{rmk}

We now define the space $_{\CC}V^{f\natural}$ and its real form
$V^{f\natural}$ as follows.
\begin{gather}
    _{\CC}V^{f\natural}=(_{\CC}V^f_L)^0\oplus
        (_{\CC}V^f_L)_{\ogi}^0,\quad
    V^{f\natural}=({V}^f_L)^0\oplus
        ({V}^f_L)^0_{\ogi}.
\end{gather}
Then in terms of the $_{\FF}U^k_{\gamma}$ we have
\begin{gather}
        _{\CC}\vfn=
            \bigoplus_{
        \substack{\gamma_k\in\{0,\omega\}\\
        \sum\gamma_k=0}}{_{\CC}U_{\gamma_1\gamma_2\gamma_3}}
        \oplus
        \bigoplus_{
        \substack{\gamma_k\in\{1,\bar{\omega}\}\\
        \sum\gamma_k=1}}{_{\CC}U_{\gamma_1\gamma_2\gamma_3}}
\end{gather}
where we use an abbreviated notation to write
${_{\CC}U_{\gamma_1\gamma_2\gamma_3}}$ for
${_{\CC}U^1_{\gamma_1}}\otimes {_{\CC}U^2_{\gamma_2}} \otimes
{_{\CC}U^3_{\gamma_3}}$, and there is a similar expression for the
real form $\vfn$ obtained by replacing $\CC$ with $\RR$ in the
subscripts of the $_{\CC}U^k_{\gamma}$. By a similar argument to
that used in \cite{FFR} to equip $_{\CC}W_L$ with VOA structure
via the para-VOA structure on $_{\CC}U$, one may also equip the
spaces $_{\CC}V^{f\natural}$ and $V^{f\natural}$ with $N=1$ SVOA
structure. Since our main focus is to study the $N=1$ SVOA $\vfn$
via its realization $\afn$, we omit a verification of this claim
and proceed directly to the task of indicating how one may arrive
at an $N=1$ SVOA isomorphism between $\vfn$ and $\afn$.

\subsection{Isomorphism}

We will concentrate on finding a correspondence between the real
$N=1$ SVOAs $\vfn$ and $\afn$. Recall that $\vfn$ may be described
as follows.
\begin{gather}\label{vfnExpInUs}
    \begin{split}
        \vfn=\bigoplus_{
        \substack{\gamma_k\in\{0,\omega\}\\
        \sum\gamma_k=0}}
        {_{\RR}U_{\gamma_1\gamma_2\gamma_3}}
        \oplus
        \bigoplus_{
        \substack{\gamma_k\in\{1,\bar{\omega}\}\\
        \sum\gamma_k=1}}
        {_{\RR}U_{\gamma_1\gamma_2\gamma_3}}
    \end{split}
\end{gather}
On the other hand, the space underlying $\afn$ is described as
$A(\gt{l})^0\oplus A(\gt{l})_{\ogi}^0$ where $\gt{l}$ is real
vector space of dimension $24$, and in particular, for a suitable
identification of $\gt{l}$ with $\bigoplus{_{\RR}\gt{a}^k}$, we
may identify $_{\RR}U_{000}=A(_{\RR}\gt{a}^1)^0\otimes
A(_{\RR}\gt{a}^2)^0\otimes A(_{\RR}\gt{a}^3)^0$ with a sub space
of $A(\bigoplus{_{\RR}\gt{a}^k})^0=A(\gt{l})^0$. As a sum of
modules over this subVOA $\afn$ admits the following description.
\begin{gather}\label{afnExpInUs}
    \begin{split}
        \bigoplus_{
        \substack{\gamma_k\in\{0,1\}\\
        \sum\gamma_k=0}}
        {_{\RR}U_{\gamma_1\gamma_2\gamma_3}}
        \oplus
        \bigoplus_{
        \substack{\gamma_k\in\{\omega,\bar{\omega}\}\\
        \sum\gamma_k=\omega}}
        {_{\RR}U_{\gamma_1\gamma_2\gamma_3}}
    \end{split}
\end{gather}
Thus it is evident that our method of constructing $\vfn$ has
almost delivered us an isomorphism with $\afn$ already. We require
to find some way of interchanging $1$ with $\w$ in the subscripts
on the right hand side of (\ref{vfnExpInUs}), and to do so we will
invoke the results of \cite{FFR} once more. It is well known that
the type $D_4$ Lie algebra admits an $S_3$ group of outer
automorphisms that has the effect of permuting transitively the
three inequivalent irreducible non-adjoint $D_4$ modules. As shown
in \cite{FFR} this action extends to the corresponding VOA
modules, and applying the outer automorphism that preserves the
spaces $U_0$ and $U_{\bar{\omega}}$, and interchanges $U_1$ with
$U_{\w}$ simultaneously to each tensor factor on the right hand
side of (\ref{vfnExpInUs}), we obtain an isomorphism of $\vfn$
with an $N=1$ SVOA ${\vfn}'$ whose underlying
$_{\RR}U_{000}$-module structure is as in (\ref{afnExpInUs}).
\begin{gather}
        {\vfn}'=\bigoplus_{
        \substack{\gamma_k\in\{0,1\}\\
        \sum\gamma_k=0}}
        {_{\RR}U_{\gamma_1\gamma_2\gamma_3}}
        \oplus
        \bigoplus_{
        \substack{\gamma_k\in\{\omega,\bar{\omega}\}\\
        \sum\gamma_k=\omega}}
        {_{\RR}U_{\gamma_1\gamma_2\gamma_3}}
\end{gather}
We have seen that the spaces $(\vfn')_{\bar{0}}$ and
$(\afn)_{\bar{0}}$ are isomorphic VOAs due to the fact that
$A(\gt{l})=A(\bigoplus_k {_{\RR}\gt{a}^k})$ and $\bigotimes_k
A(_{\RR}\gt{a}^k)$ are naturally isomorphic. Similarly,
$(\vfn')_{\bar{1}}$ is naturally isomorphic to a canonically-twisted
$A(\gt{l})$-module, and the same is true for $(\afn)_{\bar{1}}$ by
construction. The difference between $(\vfn')_{\bar{1}}$ and
$(\afn)_{\bar{1}}$ is that the former may be naturally identified
with the $A(\gt{l})^0$-module $A(\gt{l})_{\ogi,\tilde{H}}^0$ where
$\tilde{H}$ is an $\FF_{2}^{\mc{E}}$-homogeneous lift of a direct
sum of three copies of the Hamming code $\mc{H}^{\oplus 3}$, and the
later is realized as a the $A(\gt{l})^0$-module
$A(\gt{l})_{\ogi,G}^0$ for $G$ a lift of the Golay code $\mc{G}$.
Canonically twisted modules over $A(\gt{l})$ are unique up to
isomorphism, so we can be assured that $(\vfn')_{\bar{1}}$ and
$(\afn)_{\bar{1}}$ are isomorphic as $A(\gt{l})^0$-modules, and the
proof of Theorem~\ref{AfnIsSVOA} shows that the SVOA structures on
$\vfn'$ and $\afn$ are essentially unique.

What is perhaps not so clear is whether or not the $N=1$
structures on $\vfn'$ and $\afn$ coincide. Recall the following
description of the superconformal element in $\vfn$.
\begin{gather}
    \tau_V=-\frac{1}{\sqrt{8}}\sum_{\Pi}
        e^1_i(-\tfrac{1}{2})
        e^2_i(-\tfrac{1}{2})
        e^3_i(-\tfrac{1}{2})\vac\in {_{\RR}U_{111}}
\end{gather}
Since $\mc{H}$ may be defined as a quadratic residue code, it is
convenient to take $\Pi$ to be the points of the projective line
over $\FF_7$ so that $\Pi=\{\infty,0,1,2,3,4,5,6\}$ say. Then we
may choose the isomorphism $\vfn\to\vfn'$ in such a way that the
image $\tau_V'$ of $\tau_V$ in $\vfn'$ is the following.
\begin{gather}
    \tau_V'=-\frac{1}{\sqrt{8}}\sum_{\Pi}
        e^1_{i\infty}
        e^2_{i\infty}
        e^3_{i\infty}1_{\tilde{H}}
            \in (A(\gt{l})_{\ogi,\tilde{H}}^0)_{3/2}
\end{gather}
On the other hand, the superconformal element in $\afn$ is given by
$\tau_A=1_G\in (A(\gt{l})_{\ogi,G}^0)_{3/2}$. The spaces
$(A(\gt{l})_{\ogi,\tilde{H}})_{3/2}$ and
$(A(\gt{l})_{\ogi,G})_{3/2}$ are different realizations of the spin
module over $\Cl(\gt{l})$, and we may assume that $\tilde{H}$ and
$G$ are $\FF_2^{\mc{E}}$-homogeneous subgroups of $\Sp(\gt{l})$. In
the notation of \S\ref{sec:cliffalgs:mods} the spaces
$(A(\gt{l})_{\ogi,\tilde{H}})_{3/2}$ and
$(A(\gt{l})_{\ogi,G})_{3/2}$ are identified with
$\Cm(\gt{l})_{\tilde{H}}$ and $\Cm(\gt{l})_{G}$, respectively. These
two modules are equivalent and irreducible as ${\rm
Cliff}(\gt{l})$-modules, and in particular, the action of ${\rm
Cliff}(\gt{l})$ on any non-zero vector generates the entire module
in each case. The vector $\tau_A=1_G$ in $\Cm(\gt{l})_{G}$ is
determined by the property that $g1_G=1_G$ for any $g\in G$. It is a
remarkable fact that the correspondence between $\gt{l}$ and
$\bigoplus_k{_{\RR}\gt{a}}^k$ may be chosen in such a way that
$\tau_V$ also satisfies the property $g\tau_V=\tau_V$ for any $g\in
G$. Consequently we obtain explicit ${\rm Cliff}(\gt{l})$-module and
$\Sp(\gt{l})$-module equivalences between $\Cm(\gt{l})_{\tilde{H}}$
and $\Cm(\gt{l})_{G}$ such that $\tau_V'$ corresponds to $\tau_A$.
Using this isomorphism we can construct an explicit
$A(\gt{l})$-module equivalence between $A(\gt{l})_{\ogi,\tilde{H}}$
and $A(\gt{l})_{\ogi,G}$. This is the last piece of information
needed to construct an isomorphism of $N=1$ SVOAs $\vfn'\to\afn$,
and consequently we obtain
\begin{thm}\label{ThmIsom}
There is an isomorphism of $N=1$ SVOAs $\vfn \xrightarrow{\sim}
\afn$.
\end{thm}

\section{McKay--Thompson series}\label{sec:MTseries}

In this section we consider the McKay--Thompson series associated to
elements of $\Co_1$ acting on %$\vfn$ or
$\afn$. %Using $\vfn$ one may
%derive expression for elements in a certain maximal subgroup of
%$\Co_1$, and using $\afn$ we can derive explicit expressions for any
%element of $\Co_1$.

%\subsection{Series via $\vfn$}
%
%The group $\Co_1$ contains a maximal subgroup $C^f$ say, of the
%shape $2^{1+8}.(W_{E_8}'/\{\pm 1\})$. Here $W_{E_8}'/\{\pm 1\}$ is the
%quotient by $\pm 1$ of the derived subgroup of the automorphism
%group of the $E_8$ lattice $L$. Similar to the case of $\vn$ and the
%group $C=2^{1+24}.\Co_1$ which was considered in \cite[\S10.4]{FLM}
%it is possible to realize an action of $C^f$ on
%$\vfn=(V_L^f)^0\oplus (V^f_L)^0_{\ogi}$ in a reasonably explicit
%manner. Furthermore, one may derive expressions for the traces of
%elements in $C^f$ acting on $\vfn$ using just the same method as was
%used in \cite[\S10.5]{FLM}.
%
%Since we are able to compute the trace associated to any element of
%$\Co_1$ using $\afn$, we proceed directly to a description of how.

\subsection{Series via $\afn$}

Let $g\in \Co_0$, and suppose that $g^m=1$. Then there are unique
integers $p_k$ for $k|m$ such that for $\det(g-x)$, the
characteristic polynomial of $g$, we have
$\det(g-x)=\prod_{k|m}(1-x^k)^{p_k}$. This data can be expressed
using a kind of formal permutation notation as $\prod_{k|m}k^{p_k}$,
and this expression is called the {\em frame shape} for $g$. Recall
$\eta(\tau)$, the Dedekind eta function (\ref{Dedetafun}), and let
$\phi(\tau)$ be the function on the upper half plane given by
\begin{gather}
    \phi(\tau)=\frac{\eta(\tau/2)}{\eta(\tau)}
        =q^{-1/48}\prod_{n=0}^{\infty}(1-q^{n+1/2})
\end{gather}
For $g\in \Co_0$ with frame shape $\prod k^{p_k}$, we set
\begin{gather}
    \phi_g(\tau)=\prod_{k|m}\phi(k\tau)^{p_k},\quad
    \eta_g(\tau)=\prod_{k|m}\eta(k\tau)^{p_k}.
\end{gather}
The group $\Co_0$ has unique up to equivalence irreducible
representations of dimensions $1$, $24$, $276$, $2024$ and $1771$
\cite{ATLAS}. With $N$ any one of these numbers, we write $\chi_N$
for the trace function $\chi_N:\Co_0\to\CC$ on an irreducible
$\Co_0$-module of dimension $N$. Let us also write $\chi_{G}$ for
the trace function on the $\Co_0$-module $\Cm(\gt{l})_G$. (Recall
$\Cm(\gt{l})_G^0=(\afn)_{3/2}$.) We have $\chi_{G}=\chi_1
+\chi_{24}+\chi_{276} +\chi_{2024} +\chi_{1771}$, and the following
\begin{thm}\label{ThmChars}
For $\bar{g}\in \Aut(\afn)$, let $\pm g$ be the preimages of
$\bar{g}$ in $SO_{24}(\RR)$. Then we have
\begin{gather}\label{FrmChars}
    \mathsf{tr}|_{\afn}gq^{L(0)-c/24}=
        \frac{1}{2}\left(\phi_g(\tau)+\phi_{-g}(\tau)\right)
        +\frac{1}{2}(\chi_{G}(g)\eta_{-g}(\tau)
            +\chi_{G}(-g)\eta_{g}(\tau))
\end{gather}
\end{thm}
\begin{proof}
Suppose that $g\in \Co_0$ is of order $m$ and has frame shape
$\prod_{k|m}k^{p_k}$. Then $g^{-1}$ has the same frame shape. Let
$\{f_i\}_{i=1}^{24}$ be a basis for
$_{\CC}\gt{l}=\CC\otimes_{\RR}\gt{l}$ consisting of eigenvectors of
$g$ with eigenvalues $\{\xi_i\}_{i=1}^{24}$. Then we have
\begin{gather}\label{charpolyfm}
    \det(g-x)=\prod_i(\xi_i-x)= \prod_{k|m}(1-x^k)^{p_k}
\end{gather}
and we note also that $\sum_{k|m}kp_k=24$.

Recall that $\afn$ may be described as $\afn=A(\gt{l})^0\oplus
A(\gt{l})^0_{\ogi}$ where $A(\gt{l})$ is the Clifford module SVOA
associated to a $24$-dimensional inner product space $\gt{l}$, and
$A(\gt{l})_{\ogi}$ is a canonically twisted $A(\gt{l})$-module. It
is not hard to derive the following expressions for the trace of $g$
on the complexified spaces $_{\CC}A(\gt{l})$ and
$_{\CC}A(\gt{l})_{\ogi}$.
\begin{gather}
    {\sf tr}|_{_{\CC}A(\gt{l})}(-g)q^{L(0)-c/24}
        =q^{-1/2}\prod_{n\geq 0}\prod_i
            (1-\xi_iq^{n+1/2})\\
    {\sf tr}|_{_{\CC}A(\gt{l})_{\ogi}}(-g)q^{L(0)-c/24}
        =\chi_{G}(-g)q\prod_{n\geq 0}\prod_i
            (1-\xi_iq^{n+1})
\end{gather}
Substituting $q^r$ for $x$ in (\ref{charpolyfm}) and using the
fact that $\prod_i\xi_i=1$ we obtain
\begin{gather}
    \prod_i(1-\xi_iq^r)=\prod_{k|m}(1-(q^{k})^r)^{p_k}.
\end{gather}
Then for ${\sf tr}|_{_{\CC}A(\gt{l})}(-g)q^{L(0)-c/24}$ for example,
we have
\begin{gather}
    \begin{split}
    {\sf tr}|_{_{\CC}A(\gt{l})}(-g)q^{L(0)-c/24}
        &=q^{-1/2}\prod_{n\geq 0}\prod_{i}
            (1-\xi_iq^{n+1/2})\\
        &=\prod_{k|m}\left(
            q^{-kp_k/48}\prod_{n\geq 0}
            (1-(q^k)^{n+1/2})^{p_k}\right)
        =\phi_{g}(\tau)
    \end{split}
\end{gather}
and similarly, we obtain ${\sf tr}|_{_{\CC} A(\gt{l})_{\ogi}}
(-g)q^{L(0)-c/24}=\chi_G(-g)\eta_{g}(\tau)$. To compute the traces
of $\bar{g}\in\Co_1=\Co_0/\{\pm 1\}$ on $A(\gt{l})^0$ we should
average over the traces of its preimages $g$ and $-g$ on
$A(\gt{l})$, and similarly for the trace of $\bar{g}$ on
$A(\gt{l})^0_{\ogi}$. This completes the verification of
(\ref{FrmChars}).
\end{proof}

\section*{Acknowledgement}

The author is grateful to Richard Borcherds, John Conway, Gerald
H\"ohn, Atsushi Matsuo, Kiyokazu Nagatomo, Marcus Rosellen and
Olivier Schiffmann for interesting and useful discussions. The
author is also grateful to Jia-Chen Fu for lending a patient and
critical ear to many ideas. The author thanks Gerald H\"ohn for
filling a gap in the original treatment of Theorem~\ref{UniqSVOA},
and is extremely grateful to the referees, for suggesting many
improvements upon earlier versions. Finally, the author wishes to
thank his advisor Igor Frenkel for suggesting this project, and for
providing invaluable guidance and encouragement throughout its
completion.

% ------------------------------------------------------------------------
%\bibliographystyle{alpha}
%\bibliography{vaco1}
\newcommand{\etalchar}[1]{$^{#1}$}

% ------------------------------------------------------------------------
\end{document}